\documentclass[preprint]{elsarticle}

 \journal{arXiv}

\usepackage{amsmath}
\usepackage{amssymb}
\usepackage{amsthm}
\usepackage[margin=1in]{geometry}
\usepackage[utf8]{inputenc}
\usepackage{graphicx}
\usepackage{hyperref}
\usepackage{nicefrac}

\usepackage{xcolor}
\usepackage{dsfont}

\usepackage[]{appendix}

\def\RR{\mathbb R}
\def\NN{\mathbb N}
\def\eps{\varepsilon}
\def\ff{\mathcal{F}}
\def\G{\mathcal{G}}

\def\d{\,\mathrm{d}}
\newcommand{\ud}{\mathrm{d}}
\def\crl{\mathrm{curl}}
\def\de{\partial}

\newcommand{\vx}{\vec x}

\newcommand{\Xik}{X_{k}^i}
\newcommand{\Xjl}{X_{\ell}^j}
\newcommand{\xik}{x_{k}^i}
\newcommand{\xjl}{x_{\ell}^j}

\DeclareMathOperator{\supp}{supp}
\DeclareMathOperator{\dv}{div}

\newcommand{\nm}[1]{\left\| #1 \right\|}

\newtheorem{theor}{Theorem}
\newtheorem*{theor*}{Theorem}
\newtheorem{lem}{Lemma}[section]
\newtheorem{cor}[lem]{Corollary}
\newtheorem{claim}[lem]{Claim}
\newtheorem{prop}[lem]{Proposition}

\theoremstyle{definition}
\newtheorem{defi}[lem]{Definition}
\newtheorem{rem}[lem]{Remark}
\newtheorem{ex}{Example}

\begin{document}

\begin{frontmatter}

\title{Existence and regularity for a system of porous medium equations with small cross-diffusion and nonlocal drifts}

\author[GSSIaddress]{Luca Alasio}\corref{mycorrespondingauthor}
\ead{luca.alasio@gssi.it}

\author[Cambridgeaddress,FAUaddress]{Maria Bruna}
\ead{bruna@maths.cam.ac.uk}

\author[Aquilaaddress]{Simone Fagioli}
\ead{simone.fagioli@univaq.it}

\author[Cambridgeaddress]{Simon Schulz}
\ead{sms79@cam.ac.uk}

\cortext[mycorrespondingauthor]{Corresponding author}

\address[GSSIaddress]{Gran Sasso Science Institute, Viale F. Crispi 7, L’Aquila 67100, Italy}
\address[Cambridgeaddress]{Department of Applied Mathematics and Theoretical Physics, University of Cambridge, Cambridge CB3 0WA, UK}
\address[FAUaddress]{Department Mathematik, Friedrich-Alexander-Universit\"at Erlangen-N\"urnberg,
Cauerstrasse 11, 91058 Erlangen, Germany}
\address[Aquilaaddress]{Dipartimento di Ingegneria e Scienze dell’Informazione e Matematica, Universit\`a degli Studi dell’Aquila, Via Vetoio 1, 67100 Coppito, L’Aquila, Italy}

\begin{abstract}
We prove existence and Sobolev regularity of solutions of a nonlinear system of degenerate-parabolic PDEs with self- and cross-diffusion, transport/confinement and nonlocal interaction terms.
The macroscopic system of PDEs is formally derived from a large particle system and models the evolution of an arbitrary number of species with quadratic porous-medium interactions in a bounded domain $\Omega$ in any spatial dimension. 
The cross interactions between different species are scaled by a parameter $\delta<1$, with the $\delta= 0$ case corresponding to no interactions across species. 
A smallness condition on $\delta$ ensures existence of solutions up to an arbitrary time $T>0$ in a subspace of $L^2(0,T;H^1(\Omega))$. This is shown via a Schauder fixed point argument for a regularised system combined with a vanishing diffusivity approach.
The behaviour of solutions for extreme values of $\delta$ is studied numerically.
\end{abstract}


\end{frontmatter}

\begin{footnotesize}
\tableofcontents
\end{footnotesize}
\section{Introduction, motivation and set-up}
Partial differential equations involving nonlinear diffusion in combination with local or nonlocal transport are commonly used to model the macroscopic behaviour of a large number of agents or individuals in the natural, life, and social sciences. 
Systems involving many species or types of agents have been used to study multiple chemotactic populations in competition for nutrient \cite{conca_espejo_vilches_2011,EspSteVel}, tumour growth \cite{cristini2010multiscale,perthame2015some,preziosi}, pedestrian dynamics \cite{ApRDegMot,ColLec}, and  opinion formation \cite{DurMarPieWol}. Further applications can be found in population biology \cite{Bruna:2012ch,CheKol} and in semiconductor devices \cite{chen2007analysis,markowich1989system}. 

In this paper we study a class of drift-diffusion systems taking the following form:
    \begin{equation}\label{eq:main_introduction}
\de_t u_i = \dv\left[ u_i(\nabla u_i +\nabla L_i(t,x,u_i) + \delta F_i(t,x,u, \nabla u)) \right], \qquad i=1,\dots, M,
\end{equation}
where $u = (u_1, \dots, u_M)$ is a vector of non-negative functions defined on a bounded domain $Q_T = (0,T)\times\Omega$ describing the densities of $M$ subpopulations. The transport term $L_i$ models the presence of external forces and nonlocal self-interactions, $L_i(t,x,u_i)=V_i(t,x) + \big(W_i * u_i(t,\cdot)\big)(x)$,
while $F_i$ represents the interaction with the other species. The latter is of particular importance since it includes the contribution of \emph{cross-diffusion terms} through the dependence on $\nabla u$.

One of the main features of \eqref{eq:main_introduction} is that it can be used to describe cells sorting and the resulting pattern formation. This is a reorganisation process in which cells of different species---which  in principle react differently to external forces such as chemical signal or attraction/repulsion with the other species---have the propensity to group together in a delimited region \cite{KanPriVel,LemKin,Shigesada:1979tq,VasWei}.  
This biological phenomenon can be also interpreted as the inhibition or activation of growth whenever two populations occupy the same habitat, which can be  be attributed to volume or size constraints of the individual cells forming the different populations. 
In the seminal papers \cite{BerGurHil,BerGurHilPel,BerDalMim,GurPip} it was shown that segregation is induced by the presence of cross-diffusion terms. 
Nonlinear diffusion may also help in describing volume filling effects and in preventing blow-up in biological aggregation models, see \cite{bruna2017cross,CC06,HittmeierJungel,PaiHil,wilson2017reactions}.

The main goal of the present paper is to provide a well-posedness result for  \eqref{eq:main_introduction}, which can be seen as a $\delta$-order perturbation of a set of $M$ decoupled drift-diffusion equations with degenerate diffusion of porous-media type and nonlocal interactions.
One of the main difficulties for these systems is the lack of a suitable maximum principle, meaning that Sobolev estimates must be obtained in an alternative fashion. 
Under a \emph{smallness} assumption on the perturbation parameter $\delta$, we are able to estimate each density in $L^2(0,T;H^1(\Omega))$, as well as the relative time derivatives in a dual Sobolev space.
Such estimates allow us to show convergence of weak solutions for a proper regularised approximating sequence.  
The structure of the system and the numerical simulations indicate that the expected critical value for the $L^2(0,T;H^1(\Omega))$-framework is $\delta=1$. However, technical conditions impose a stricter bound on $\delta$. The $BV$ norm is the only norm expected to remain bounded as $\delta \to 1$.

The paper is organised as follows. 
In Section \ref{sec:derivation} we present a formal derivation of system  \eqref{eq:main_introduction} from a system of interacting particles.
In Section \ref{sec:general} we provide the general set of  assumptions and the statement of the main result,  Theorem \ref{theor:main}. 
We sketch the main steps in the strategy of the proof of the main result in Section \ref{sec:strategy}, and we give a brief overview of different approaches and existing results in the literature in Section \ref{sec:null}. 
Section \ref{sec:numerics} is devoted to the numerical investigation of a particular case of system \eqref{eq:main_introduction}, which helps to highlight the influence of the parameter $\delta$ in the time evolution of solutions and its norms. 
The remainder of the paper focuses on the proof of Theorem \ref{theor:main}. In Section \ref{sec:regularised} we introduce a regularised system, providing existence, uniqueness and uniform estimates for these regularised solutions. 
In Section \ref{sec:fixed_point} we present a fixed point argument and show compactness properties for the solution map of the regularised system.
Convergence of approximate solutions to weak solutions of system \eqref{eq:main_introduction} is established in Section \ref{sec:vanishing} through vanishing diffusivity. 
Appendix \ref{sec:appendix united} collects some technical results used in the paper.

\subsection{Model derivation}\label{sec:derivation}
In this subsection we sketch a formal derivation  of \eqref{eq:main_introduction} starting from the interacting particle system with $M$ species, each composed of $N_i$ identical particles, $i = 1, \dots, M$. 
To simplify the presentation, the derivation is shown for the case for $L_i(t,x,u_i) \equiv V_i(x)$, that is, we drop the nonlocal interactions $W_i$ and the time-dependence, and a simple cross-diffusion term $F_i(t,x,u,\nabla u) \equiv \sum_{j=1,j\ne i}^M \nabla u_j$. The addition of nonlocal interactions and time is straightforward.
We denote by $N$ the total number of particles, $N = \sum_i N_i$. We consider the following model:
\begin{subequations}
\label{amb}
\begin{align}
\label{amb1}
	\ud \Xik(t) &= -\nabla V_i(\Xik) \ud t - \sum_{j = 1}^M \sum_{\substack{\ell=1 \\ (\ell,j) \ne (k,i)}}^{N_i} \nabla K_{ij}(\Xik - \Xjl) \ud t,\\
	\label{amb2}
	\Xik(0)&= \xi_{k}^i, \qquad k = 1, \dots, N_i,
\end{align}
where $\Xik(t)$ is the position of the $k$-th particle in the $i$-th species at time $t$, evolving in a bounded domain $\Omega \subset \mathbb R^d$ such that $|\Omega| = 1$. Particles are initialised with $\xi_1^i, \dots, \xi_{N_i}^i$ independent and identically distributed random variables with the common probability density function $u_{i,0}$. 
\end{subequations} 
Here $K_{ii}$ denotes the self-interaction potential in species $i$, and $K_{ij}$ denotes the cross-interaction potential between species $i$ and $j$. Note that $K_{ij}$ and $K_{ji}$ may differ to represent an asymmetric interaction between the two species. The potentials are assumed to be obtained from some fixed function $K_0$ by the scaling
\begin{equation} \label{scaling}
	K_{ij}(x) = \chi_{ij} K_0 \left(\frac{|x|}{\varepsilon_{ij}} \right),
\end{equation}
where the parameters $\chi_{ij}$ and $\varepsilon_{ij}$ represent the strength and the range of the interactions respectively, and depend on $N_j$ in a way that will be made specific later on. The scale-free potential $K_0: \mathbb R ^d \to \mathbb R$ is a radial, nonnegative function whose gradient is locally Lipschitz outside the origin. Moreover, it is assumed that $\|K_0\|_{L^1} <\infty$. Without loss of generality, we set $\|K_0\|_{L^1}  = 1$.

Depending on $\chi$ and $\varepsilon$, one expects different limit equations \cite{Bodnar:2005kv}. For example, when the interactions are long range ($\varepsilon \sim 1$) and weak ($\chi \sim N^{-1}$), then one particle interacts on average with an order $N$ particles as $N \to \infty$ one recovers a mean-field limit for weakly interacting particles. In contrast, the case of moderately interacting particles corresponds to stronger but more localised interactions, so one particle interacts with fewer particles. As a result, one expects interactions to emerge as local terms in the limit equation. 

We define the total interaction potential of the $i$-th species as
\begin{equation}
K_i(\vec x) = \sum_{k = 1}^{N_i} \bigg[ \sum_{\ell> k}^{N_i} K_{ii}(\xik - x_\ell^i) + \sum_{\substack{j = 1 \\ j \ne i}}^M \sum_{\ell=1}^{N_j}	K_{ij}(\xik- \xjl) \bigg]
\end{equation}
where $\vx = (\xik)_{k = 1, \dots, N_i, i = 1, \dots, M}$. 
Then the joint probability density $P_N(\vx,t)=\text{Prob}(X(t) = \vx)$  of $N$ particles evolving according to \eqref{amb} satisfies the following equation
\begin{subequations}
\begin{align} \label{eqN}
\de_t P_N = \sum_{i  =1}^M\sum_{k=1}^{N_i} 
	\nabla_{\xik} \cdot  \big[ \nabla V_i(\xik) P_N +  \nabla_{\xik} K_i(\vx) P_N  \big], \qquad \vx\in \Omega^{N} , t>0,
\end{align}
together with boundary conditions
\begin{align} \label{bc1}
	\nu \cdot \big[ \nabla V_i(\xik) P_N +  \nabla_{\xik} K_i(\vx) P_N  \big] &= 0, \qquad \xik \in \partial \Omega,t > 0, 
\end{align}
for $k = 1, \dots, N_i$ and $i = 1, \dots, M$, where $\nu$ is the outward  normal on $\partial \Omega$ and the other coordinates are in $\Omega$.
\end{subequations}

We consider the one-particle densities for each species as
\begin{equation}
	u_i(x,t) = \int_{\Omega^N} P_N(\vx, t)  \delta(x_1^i - x) \ud \vx,
\end{equation}
where we note that the choice of $x_1^i$ is unimportant (since within a subpopulation, particles are indistinguishable).  To obtain the equation for $u_1(x)$, we integrate \eqref{eqN} over all particle positions except one particle in the first species and use the boundary conditions \eqref{bc1}:
\begin{align} \label{equ1}
\de_t u_1 = \dv \left[ \nabla V_1(x) u_1 +  G_1(x,t) \right],
\end{align}
where
\begin{align*}
G_1(x,t) &= \int_{\Omega^N} \nabla_{x_1^1} K_1(\vx) P_N \delta(x_1^i - x) \ud \vx \\
&= \int_{\Omega^N} \bigg[\sum_{\ell = 2}^{N_1} \nabla_{x} K_{11}(x-x_\ell^1) P_N  + \sum_{j = 2}^M \sum_{\ell = 1}^{N_j} \nabla_{x} K_{1j}(x-x_\ell^j) P_N \bigg] \ud \vx \\
&= (N_1-1) \int_\Omega \nabla_{x} K_{11} (x-y) P_2^{11}(x,y,t) \ud y  + \sum_{j = 2}^M N_j \int_\Omega \nabla_{x} K_{1j} (x-y) P_2^{1j}(x,y,t) \ud y.
\end{align*}
Here $P_2^{1j}, j = 1, \dots, M$ stands for the following two-particle density
\begin{align*}
P_2^{1j}(x,y,t) &= \int_{\Omega^N} P_N(\vx,t) \delta(x_1^1-x)\delta(x_2^j-y) \ud \vx.
\end{align*}

Oelschl\"ager \cite{oelschlager1990large} proved propagation of chaos (meaning that any fixed number of particles remains approximately independent in time despite the interaction) for the single-species cases similar to \eqref{amb} under quite restrictive initial conditions. These conditions were relaxed by Philipowski \cite{Philipowski:2007gd} by means of using regularising Brownian motions, that is, adding terms $\epsilon \ud B_k^i$ in \eqref{amb1} and taking $\epsilon \to 0$ at a suitable rate depending on $N, \varepsilon$. So including such a term would make sense when considering a rigorous derivation. An alternative approach taken by \cite{chen2020rigorous} in the multiple species case is to include the interactions between particles in the diffusion term (note that, in their case, the mean-field limit model still contains linear diffusion terms). For our purposes, here we simply assume that an analogous propagation of chaos for $P_N$ exists. 
In particular, this  means that the two-particle marginals may be approximated by
$$
P_2^{11}(x,y,t)  = u_1(x,t) u_1(y,t), \qquad P_2^{1j}(x,y,t) = u_1(x,t) u_j(y,t).
$$
Using these expressions in $G_1$, it reduces to 
$$
G_1(x,t) = (N_1 - 1) u_1 \nabla (K_{11} \ast u_1) + \sum_{j = 2}^M N_j u_1 \nabla (K_{1j} \ast u_j).
$$
Finally, if we consider the scaling \eqref{scaling} with $\varepsilon_{11}, \varepsilon_{1j} \ll 1$, we can localise the convolution terms and arrive at
\begin{equation}
	\label{interaction_term}
	G_1(x,t) = (N_1 - 1)\varepsilon_{11}^d \chi_{11} u_1 \nabla u_1 + \sum_{j = 2}^M N_j  \varepsilon_{1j}^d \chi_{1j}  u_1 \nabla u_j,
\end{equation}
using that $\|K_0\|_{L^1} = 1$. The analogous calculation can be done for any of the other species to obtain $G_i(x,t)$ for $i = 1, \dots, M$.
Now we can determine the suitable scaling for interactions that lead to the structure in \eqref{eq:main_introduction}. Namely, we set the strengths to be
$$
\chi_{ii} = \frac{1}{(N_i-1)\varepsilon_{ii}^{d}},\qquad \chi_{ij} = \frac{\delta}{N_j \varepsilon_{ij}^d},
$$
with $\delta \ll 1$, and we let $N_i \to \infty, \varepsilon_{ij} \to 0$ in such a way that $N_j \gg 1/\varepsilon_{ij}$. Using this and combining \eqref{equ1} and \eqref{interaction_term}, we arrive at 
\begin{equation}\label{eq:main_derivation}
\de_t u_i = \dv\big[
u_i\big(\nabla u_i +\nabla V_i(x) + \delta \sum_{j = 1 , j \ne i}^M \nabla u_j \big) \big]. 
\end{equation}

\subsection{Assumptions and notion of solution}\label{sec:general}

Let $\Omega\subset \RR^d$ be an open bounded domain of class $C^2$ with outward normal denoted by $\nu$, and let $T>0$. We denote the parabolic cylinder by $Q_T := (0,T)\times\Omega$, the lateral boundary by $\Sigma_T := (0,T) \times \partial\Omega$, and the closure by $\bar{Q}_T := [0,T]\times\bar{\Omega}$. 
Consider also $M$ functions 
$$F_i: \bar{Q}_T\times\RR^M\times\RR^{d\times M} \to \RR^{d}
$$ such that the dependence on the last argument is affine, namely
\begin{equation}\label{eq:F structure 0}
F_i(t,x,z,p) = G_i^0(t,x,z) + \sum_{j=1}^M G_{ij}^1(t,x,z)p_j, \qquad \text{for } i \in \{1,\dots,M\}, 
\end{equation}
where $\{G_i^0\}_{i=1}^M$ are vector functions with values in $\mathbb{R}^d$ and $\{G_{ij}^1\}_{i,j=1}^M$ take values in the space of $d\times d$ matrices; all of the functions above are assumed to be $C^2$-regular and uniformly bounded with respect their  arguments. 
The assumption of regularity is not optimal and can be relaxed, but it avoids unpleasant technicalities in Section \ref{subsec:regularised} (see also Remark \ref{rem:unpleasant technicalities} in Appendix \ref{sec:appendix proof russian}). 
We emphasise that, in \eqref{eq:F structure 0}, for each $j\in\{1,\dots,M\}$, the quantity $p_j$ is a column vector in $\mathbb{R}^d$, while $p$ is a $d \times M$ matrix whose $j$-th column is the column vector $p_j$. 
We assume that there exists a positive constant $C_F = C_F(T,\Omega)$ such that
\begin{equation}\label{eq:fbound0}
|F_i(t,x,z,p)|  \leq C_F(1+|p|), \quad \forall (t,x,z,p) \in \bar{Q}_T\times\RR^M\times\RR^{d\times M}.
\end{equation}
For $i \in \{1,\dots,M\}$ and $\delta \in \mathbb{R}$ some small constant to be made precise, consider the following system of equations:
\begin{equation}\label{eq:main0}
\left\lbrace\begin{aligned}
&\de_t u_i = \dv\left[
u_i(\nabla u_i +\nabla L_i(t,x,u_i) + \delta F_i(t,x,u, \nabla u))
\right]
\qquad && \text{in } Q_T,
\\
&0 = \nu \cdot\left[
u_i(\nabla u_i +\nabla L_i(t,x,u_i) + \delta F_i(t,x,u, \nabla u))
\right]
\qquad && \text{on } \Sigma_T,
\\
&u_i(0,\cdot) = u_{i,0}
\qquad && \text{on } \Omega,
\end{aligned}\right.
\end{equation}
where $u = (u_i)_{i=1}^{M}$ is the unknown vector-valued function, and each $u_{i,0}$ is a given non-negative function in $L^p(\Omega)$ for $p>1$. The transport terms $\{L_i\}_{i=1}^M$ are prescribed by 
\begin{equation}\label{eq:Li terms}
 L_i(t,x,u_i)=V_i(t,x) + \big(W_i * u_i(t,\cdot)\big)(x), 
\end{equation}
where we assume $W_i$ to be radially symmetric for each $i \in \{1,\dots,M\}$, and the convolution to be only with respect to the space variable, i.e., 
\begin{equation*}
    \big(W_i * u_i(t,\cdot)\big)(x) = \int_\Omega W_i(x-y) u_i(t,y) \d y \qquad \text{a.e.~} (t,x) \in Q_T. 
\end{equation*}
Again, in order to avoid unpleasant technicalities in Section \ref{subsec:regularised}, we assume $\{V_i,W_i\}_{i=1}^M$ to be $C^2$-regular and uniformly bounded in all of their respective arguments. 
In particular, there exists a positive constant $C_L$ such that 
\begin{equation}\label{eq:drift V bound}
    \max_{i \in \{1,\dots,M\}} \Vert V_i \Vert_{C^2(\mathbb{R}^{d+1})} + \max_{i \in \{1,\dots,M\}} \Vert W_i \Vert_{C^2(\mathbb{R}^d)}  \leq C_L. 
\end{equation}

\begin{rem}
We note that minor adaptations of our approach allow to treat additional cross-interaction terms of the form $\dv \left( u_i \nabla \sum_{j,k=1}^M W_j*u_k(t,\cdot)(x)\right)$ in \eqref{eq:main0}. However, in order to do this, cross-interaction terms must be included as part of the term $F_i$, i.e., they must be premultiplied by the small parameter $\delta$. 
\end{rem}

The definition of the function space that we use depends on the initial data as follows.

\begin{defi}[Function space]\label{defi:functional space}
We define the following Banach space:
\begin{equation*}
\begin{aligned}
    \Xi := \bigg\lbrace u : Q_T \to \mathbb{R}^M \big| u \in ( L^2(0,T&;H^1(\Omega)) )^M \text{, } \partial_t u \in (X')^M \text{, and, for } i \in \{1,\dots,M\}, \\ 
& u_i \geq 0 \text{ a.e.~in } Q_T, \quad \int_\Omega u_i(t,x) \d x = \int_\Omega u_{i,0}(x) \d x \text{ a.e.~}t \in [0,T] \bigg\rbrace,
\end{aligned}
\end{equation*}
where 
\begin{equation*}
X := L^{r}(0,T;W^{1,r}(\Omega)), \qquad X' = L^{r'}(0,T;(W^{1,r}(\Omega))')
\end{equation*}
with $r := 2(d+1)$ and $r' = (2d+2)/(2d+1)$.
\end{defi}

\begin{defi}[Weak solution]\label{defi:weak sol 0}
Fix an arbitrary $T>0$. Given the non-negative functions  $(u_{i,0})_{i=1}^M$ belonging to $L^p(\Omega)$ for $p>1$, we say that the vector-valued function $ u = (u_i)_{i=1}^M \in \Xi$ , is a weak solution of \eqref{eq:main0} if:
\begin{enumerate}
    \item for any test function $\phi \in C^1(\bar{Q}_T)$ and, for each $i \in \{1,\dots,M\}$, there holds 
\begin{equation}\label{eq:weak solution}
    \begin{aligned}
        \langle \partial_t u_i , \phi \rangle_{X'\times X} +\int_{Q_T} u_i(\nabla u_i + \nabla L_i(t,x,u_i) + \delta F_i(t,x,u,\nabla u))\cdot\nabla\phi \d x \d t = 0;
    \end{aligned}
\end{equation}
\item for each $i\in\{1,\dots,M\}$, the function $u_i$ is non-negative a.e.~in $Q_T$ and conserves its initial mass, i.e., 
\begin{equation}\label{eq:conserve initial mass requirement}
    \int_\Omega u_i(t,x) \d x = \int_\Omega u_{i,0}(x) \d x \qquad \text{a.e.~} t \in (0,T); 
\end{equation}
\item for each $i \in \{1,\dots,M\}$, we have $u_i \in C([0,T];(W^{1,r}(\Omega))')$ (\textit{cf.}~Remark \ref{rem:embed into cts negative sobolev space}) and the initial datum is satisfied in the $(W^{1,r}(\Omega))'$ sense, i.e.,
$
    \lim_{t \to 0^+}\Vert u_i(t,\cdot) - u_{i,0} \Vert_{(W^{1,r}(\Omega))'} = 0.
$
\end{enumerate}
\end{defi}

\begin{rem}\label{rem:comments about weak formulation}
Observe that the second condition in Definition \ref{defi:weak sol 0} implies that any weak solution $u=(u_i)_{i=1}^M$, for $i \in \{1,\dots,M\}$, belongs to $L^\infty(0,T;L^1(\Omega))$ and
$
\Vert u_i(t,\cdot) \Vert_{L^1(\Omega)} = \int_\Omega u_i(t,x) \d x = \int_\Omega u_{i,0}(x) \d x$, $\forall t \in [0,T]. 
$
\end{rem}

\begin{rem}\label{rem:embed into cts negative sobolev space} \label{rem:clarify duality prod under weak sol def}

Observe that $\Xi \subset (C([0,T];(W^{1,r}(\Omega))'))^M$. Indeed, let $u \in \Xi$ and $\varphi \in X$ with $\Vert \varphi \Vert_X \leq 1$ be arbitrary. Then, using the H\"{o}lder inequality, for any $i \in \{1,\dots,M\}$, 
\begin{equation*}
    \left|\int_{Q_T} u_i \varphi \d x \d t \right| \leq \Vert u_i \Vert_{L^2(Q_T)} (|\Omega|T)^{\frac{d}{r}} \Vert \varphi \Vert_X. 
\end{equation*}
Thus, $\Vert u_i \Vert_{X'} \leq \Vert u_i \Vert_{L^2(Q_T)} (|\Omega|T)^{\frac{d}{r}}$. Meanwhile, we also have $\partial_t u_i \in X'$ by the definition of $\Xi$, and it therefore follows that $u_i$ belongs to $W^{1,r'}(0,T;(W^{1,r}(\Omega))')$. By \cite[Theorem 2 of Section 5.9.2]{evans1998partial}, it follows that $u_i \in C([0,T];(W^{1,r}(\Omega))')$ for every $i\in\{1,\dots,M\}$. 

We identify 
\begin{equation}\label{eq:what the dual means}
    \langle \partial_t u_i , \phi \rangle_{X'\times X} = \int_0^T \langle \partial_t u_i(t,\cdot) , \phi(t,\cdot) \rangle_\Omega \d t \qquad \forall \phi \in C^1(\bar{Q}_T), 
\end{equation}
where $\langle \cdot ,\cdot \rangle_\Omega$ is the duality product of $W^{1,r}(\Omega)$. Moreover, given $u=(u_i)_{i=1}^M \in \Xi$, the weak formulation \eqref{eq:weak solution} is equivalent to 
\begin{equation*}
    \begin{aligned}
        \int_{Q_T} \big( -u_i \partial_t\phi + u_i(\nabla u_i + \nabla L_i(t,x,u_i) + \delta F_i(t,x,u,\nabla u))\cdot\nabla\phi \big) \d x \d t = 0, 
    \end{aligned}
\end{equation*}
for $i \in \{1,\dots,M\}$, for any $\phi \in C^1(\bar{Q}_T)$ with $\phi(0,\cdot)=\phi(T,\cdot)=0$ in $\bar{\Omega}$; see Lemma \ref{lem:reg sol mirrors weak sol def}. 
\end{rem}

Our main result is the following. 
\begin{theor}\label{theor:main}
Let $(u_{i,0})_{i=1}^M$ be non-negative functions belonging to $L^p(\Omega)$ for $p>1$, and $\delta \in \mathbb{R}$ be such that
\begin{equation} \label{delta_bound}
\alpha \delta^2 C_F^2 C_\Omega < 1,
\end{equation}
where $C_F$ is specified in \eqref{eq:fbound0}, $\alpha$ depends only on $\Omega$ and the smoothing operator \eqref{smoothing operator def}, and $C_\Omega$ is given in \eqref{eq:big constant for L2H1 bound 0}. 
Then there exists a weak solution $u = (u_i)_{i=1}^M$ of \eqref{eq:main0}, in the sense of Definition \ref{defi:weak sol 0}. Moreover, there exists a positive constant $C=C(\Omega,T,d,\delta)$, prescribed by 
\begin{equation}
    C = C_\Omega \big( 1 - \delta^2 C_F^2 C_\Omega \big)^{-1},
\end{equation}
 such that, for $i\in\{1,\dots,M\}$, 
\begin{equation}\label{eq:estimates for main theorem i}
\Vert u_i \Vert_{L^2(0,T;H^1(\Omega))}^2 \leq
C \bigg( 1 + \Vert u_{i,0} \Vert_{L^1(\Omega)}^2 + \int_\Omega u_{i,0}\log u_{i,0} \d  x \bigg), 
\end{equation}
and there exists another positive constant $C'=C'(\Omega,T,d,\delta)$, such that, for $i\in \{1,\dots,M\}$, 
\begin{equation}\label{eq:estimates for main theorem ii}
\nm{\de_t u_i}_{X'}
\leq
C' \bigg( 1 + \Vert u_{i,0} \Vert_{L^1(\Omega)}^2 + \int_\Omega u_{i,0}\log u_{i,0} \d  x \bigg). 
\end{equation}
\end{theor}

\subsection{Strategy}\label{sec:strategy}
We summarise the strategy for the proof of Theorem \ref{theor:main} as follows:
\begin{itemize}
\item {\em Weak solution for regularised frozen system} [Subsection \ref{subsec:regularised}]: We consider a decoupled, regularised system with unknown $z$ instead of $u$ to distinguish it from the solution of the original coupled system. The decoupled system, namely \eqref{eq:pdefrozen}, is obtained by ``freezing'' the cross-diffusion terms. In particular, we replace the unknown vector-valued function $z$ with a given function $\bar{z}$ and, eventually, we shall identify $z$ and $\bar{z}$ via a fixed point argument. We study solutions $z \in (C^{2,1}(\bar{Q}_T))^M$ according to Definition \ref{defi:weak sol regularised}. Existence, uniqueness, non-negativity and mass preservation of the solutions are shown in Lemma \ref{lem:existence of viscous approximates 0}.
\item {\em Uniform estimates and uniqueness for regularised frozen system} [Subsections \ref{subsec:uniform} \& \ref{subsec:unique}]: In Lemmas \ref{lem:energy estimates 0} and \ref{lem:time compactness bound 0}, we derive uniform estimates with respect to the regularisation parameters for solutions of the regularised system.
We obtain $H^1$-type bounds for $z$ and bounds in a dual Sobolev space for $\de_t z$. Uniqueness in $\Xi$ is obtained introducing a suitable dual problem.
\item {\em Weak compactness of the solution map} [Subsection \ref{section:weak compactness}]: We construct a solution operator $S$ for the regularised system \eqref{eq:pdefrozen} associating $z$ to $\bar{z}$ as in \eqref{eq:solmap}.
In Lemma \ref{lem:weak compactness}, we show that the map $S$, composed with a suitable regularising operator \eqref{eq:Rnu introduce}, is sequentially weakly compact in a suitable Sobolev space.
\item {\em Strong compactness of the solution map} [Subsection \ref{subsec:strong}]: In Lemma \ref{lem:strong compactness} we improve the compactness result and show that the solution map is strongly compact in $\Xi$. To do this, we exploit the lower semicontinuity of the Fisher information and apply the div-curl Lemma. 
\item {\em Vanishing diffusivity} [Section \ref{sec:vanishing}]: Thanks to a variant of Schauder's Fixed Point Theorem, in Proposition \ref{cor:existence coupled sys}, we obtain existence of solutions of the coupled system \eqref{pdecoupled 0}, which corresponds to original system \eqref{eq:main0} with artificial diffusivity.
Finally, we let the diffusivity vanish and prove Theorem \ref{theor:main}.
\end{itemize}

\subsection{Null results: what we tried and did not work}\label{sec:null}

The one-species counterpart of system \eqref{eq:main_introduction}
has been largely studied in the literature, see for example \cite{ambrosio2008gradient,bertozzi2009existence} and references therein. Nevertheless, a complete well-posedness theory for cross-diffusion systems in presence of transport term is not currently available.

Indeed, the separation process between species described at the beginning of this section may lead to discontinuities of the densities and in their derivatives at the interface between different species.
These issues are also accentuated by the presence of degenerate diffusion terms that, on the one hand, determine finite speed of propagation in the supports of the solutions and, on the other hand, cause a possible loss of regularity at their boundary.
In order to better explain the difficulties that cross-diffusion terms bring to the analysis, let us consider the following special case of \eqref{eq:main_introduction} for $M=2$:
 \begin{equation}\label{eq:easy_intro}
  \partial_{t}u_i= \partial_x\left[u_i\partial_x(u_i+\delta u_j)-u_i V_i'(x)\right], \quad i,j=1,2,\,j\neq i.
 \end{equation}
Consider an abstract splitting scheme built as follows: given  $\bar{u}_i$, $i=1,2$, we solve the 
 \emph{decoupled} equations 
\[
\begin{cases}\partial_t \rho_i = \partial_x\left[\rho_i \partial_x\rho_i-\rho_i V_i'(x)\right],\\
 \rho_i(t_0,x)=\bar{u}_i(x),
\end{cases}\quad i=1,2,
\quad t\in [t_0,t_1]
\]
which can be done by means of several results in the literature for nonlinear-diffusion and transport equations, see \cite{ambrosio2008gradient,KimZhangII}. In the next step we use the densities $\rho_i$ obtained above and evaluated at time $t_1$ as initial data for the the following hyperbolic equations 
\[
\begin{cases}
\partial_t w_i = \delta\partial_x\left[ w_i \partial_x \rho_j \right],\\
 w_i(t_1,x)=\rho_i,
\end{cases} \quad i,j=1,2,\,j\neq i,
\quad t\in [t_1,t_2],
\]
where $\partial_x \rho_j$ is given by the previous iteration and is frozen at time $t_1$. Each of these two steps is then iterated for a sequence of times $t_0<t_1<t_2<t_3<\dots<T$.
The regularity of the first (diffusive) step, induced by the quadratic porous medium term (see \cite{KimZhangII,vazquez2007porous}), is insufficient to ensure well-posedness in the second step, see  \cite[Chap.~1, Sect.~2]{AmbrosioNotes} and \cite{AmbCri,BouJam}. 
This simple argument highlights how cross-diffusion models have a strongly hyperbolic nature. 

According to the classical theory in \cite{ladyzhenskaia1988linear}, the well-posedness of \eqref{eq:easy_intro} is related to positive definiteness of the diffusion matrix
\begin{equation*}
    D(u_1,u_2)= \begin{pmatrix} u_1 & \delta u_1 \\ \delta u_2 & u_2 \end{pmatrix}.
\end{equation*}
Since the matrix above is not symmetric we must consider its symmetric part $ \frac12(D+D^T)$, which has determinant 
\[
\det\left( \frac12(D+D^T)\right)=\left(1-\frac{\delta^2}{2}\right)u_1u_2-\frac{\delta^2}{4}(u_1^2+u_2^2).
\]
From the above it is clear that the presence of cross-diffusion may induce a negative quadratic form. 
The lack of uniform parabolicity (namely the failure of $D(u_1,u_2) \geq c\mathbb{I}$, for some $c > 0$)  is present in several applications such as \cite{blanchet2006two,stara1995some}, and has attracted a lot of interest in recent years.

As pointed out in the splitting scheme sketched above, the main issue lies in the difficulty of providing \emph{a priori} estimates for the single components $u_1$ and $u_2$ and on their space derivatives.
Several attempts were made in this direction, even trying to extend the concepts of parabolicity, see the classical references \cite{amann1989dynamic,le2006everywhere,PieSch}.

In presence of reaction terms instead of the transport terms, an existence theory was obtained in a one-dimensional $BV$-setting for the case $\delta=1$ in \cite{BerGurHil,BerGurHilPel,carrillo2018splitting}; see also \cite{Pert2019} for a multi-dimensional result.  
The $BV$-setting is somehow natural since the emergence of ``segregated solutions'' is highlighted in several contexts \cite{Bur2020,Bur2018,Carrillo:2017uq}, see also \cite[Chap.~1, Sect.~4]{AmbrosioNotes}. 
A general existence theory for systems with arbitrary cross-diffusion terms and local/nonlocal transport is far from being completed.

In order to achieve a satisfactory theory, many results in the literature have been inspired by the gradient flow structure that can be associated to systems in the form of \eqref{eq:easy_intro}. These can be split into two categories: formal gradient flow structure and a Wasserstein gradient flow theory.
In the first group we mention the works \cite{Burger:2010gb,Jun2006,desvillettes2015entropic,jungel2015boundedness,jungel2012entropy}, where a formal gradient flow formulation provides the estimates needed to prove global existence.
The second approach concerns the many-species version of Wasserstein gradient flow theory of \cite{ambrosio2008gradient,San15a}. 
Such an approach has been already successfully used in \cite{DifFag1}  for a system of nonlocal interaction equations with two species and non symmetric cross-interactions, and first used in a system with cross-diffusion terms but no transport terms in \cite{LauMat}. 
Other results, only apply to diagonal diffusion and in some cases only in bounded domains \cite{CarLab1,CarLab2}, or with \emph{dominant diagonal parts} \cite{di2018nonlinear}, see also \cite{alasio2018stability,desvillettes2015new}.
In one space dimension, \cite{KimMes} provides an existence result for cross-diffusion systems with \emph{ordered} external potentials.

\section{Numerical investigation} \label{sec:numerics}

In this section we present numerical simulations of \eqref{eq:main_introduction} with two species ($M=2$) in one dimension ($d=1$). We consider the case with $L_i(t,x, u_i) =V_i(x)$ and $F_i(t,x,u, \nabla u) \equiv u_j$, leading to the following system of equations:
\begin{subequations}
	\label{eq:numerics}
\begin{align}
\de_t u_1 &= \dv\left[
u_1\left(\nabla u_1 +\nabla V_1(x) + \delta \nabla u_2 \right) \right],\\
\de_t u_2 &= \dv\left[ u_2\left(\nabla u_2 +\nabla V_2(x) + \delta \nabla u_1 \right) \right].
\end{align}	
\end{subequations}
Throughout this section we consider the domain $\Omega = [-1,1]$ with no-flux boundary conditions, and initial conditions $u_1(0,x) = u_{1,0} (x)$ and $u_2(0,x) = u_{2,0}(x)$ with unit mass. 
We solve \eqref{eq:numerics} using the positivity-preserving finite-volume scheme presented in \cite{Carrillo:2017uq}, which is first order in space and time. We use $J = 64$ grid points in space and a fixed timestep $\Delta t = 10^{-6}$.

We consider what the bound \eqref{delta_bound} on $\delta$ is for our particular examples. For our choice of cross-term $F_i$, we have that $C_F = 1$ (see \eqref{eq:fbound0}) and $C_\Omega$ in \eqref{eq:big constant for L2H1 bound 0} simplifies to 
\begin{equation*}
    C_\Omega(T) = 2  \max\left\lbrace (1+C_P) , \frac{T}{2} + 4(1+C_P)\big( e^{-1} + 2 T C_L^2 \big) \right\rbrace, 
\end{equation*}
with Poincar\'e constant $C_P = (2/\pi)^2$ and $C_L = \max_i \| V_i\|_{L^2}$. For the purposes of the numerical simulation, it is convenient to consider $C_\Omega(\Delta t)$. In the limit of $\Delta t \to 0$, we have
$$
C_\Omega(0) = 8 e^{-1} (1 + C_P) \approx 4.136,
$$
independent of the external potentials $L_i$. Therefore, the upper bound on $\delta$ is given by
$$
\delta < \delta_{\max} = 1/\sqrt{C_\Omega(0)} \approx 0.492.
$$
Below we numerically investigate the behaviour beyond such value, and close to the critical value $\delta = 1$.

\begin{ex}[Left and right initial conditions] \label{ex:1} In the first example we consider the initial conditions
$$
u_{1,0}(x) = C_1[(x+0.5)(-0.9-x)]_+, \qquad u_{2,0}(x) = C_2[(x-0.5)(0.9-x)]_+,
$$
where $C_1,C_2$ are such that the initial densities are normalised to unit mass.  We consider the external potentials  $V_1 = 0$ and $V_2(x) = 2 x^2$, and four different values of $\delta = 0.4, 0.6, 0.8, 0.99$. For this choice of potentials, we have $C_L = 6$ and 
$$
\delta_{\max}(T) = 0.0203, \quad \delta_{\max}(\Delta t) = 0.492,
$$
where the upper bound on $\delta$ is $\delta_{\max}(t) = 1/\sqrt{C_\Omega(t)}$.

In the left column of Fig.~\ref{fig:timeevol1} we show the time evolution at ten equally spaced times  $t_k$ between 0 and $T = 3$ of $u_1$ and $u_2$ (solid blue and red lines, respectively) as well as the corresponding steady states $u_{1,\infty}(x)$ and $u_{2,\infty}(x)$ (dashed green and purple lines respectively), which are computed as the minimisers of the energy
\begin{equation}
	\label{energy}
	E[u_1,u_2](t) = \int_\Omega \left[ \frac{1}{2} u_1^2 + \frac{1}{2} u_2^2 + \delta u_1 u_2 + V_1 u_1 + V_2 u_2\right] \ud x.
\end{equation}
As we increase $\delta$ closer to one (the value at which $E$ stops being strictly convex), we observe the formation of sharper interface between the two components. For the smallest value of $\delta$, $\delta = 0.5$, there is no ``vacuum region''  for $u_1$ due to $u_2$ (that is, $\supp u_1 = \Omega$) and by $t = 3$ the solution is very close to the steady state. Increasing $\delta$ changes this: for larger $\delta$, the stationary solution $u_1$ has a vacuum region in the middle of the domain in which $u_1 = 0$ to numerical precision
(which grows closer to $\supp u_2$ as $\delta$ approaches one). This vacuum region implies that it takes much longer for half of the mass of $u_1$ to transfer from the left to the right on the domain, implying that the equilibration to the stationary state is much slower (this can be clearly seen in the bottom row, where the final time solution $u_1(T, x)$ is still very far from the steady state minimiser $u_{1,\infty}$).
\begin{figure}
\begin{center}
	\includegraphics[width = 0.75\textwidth]{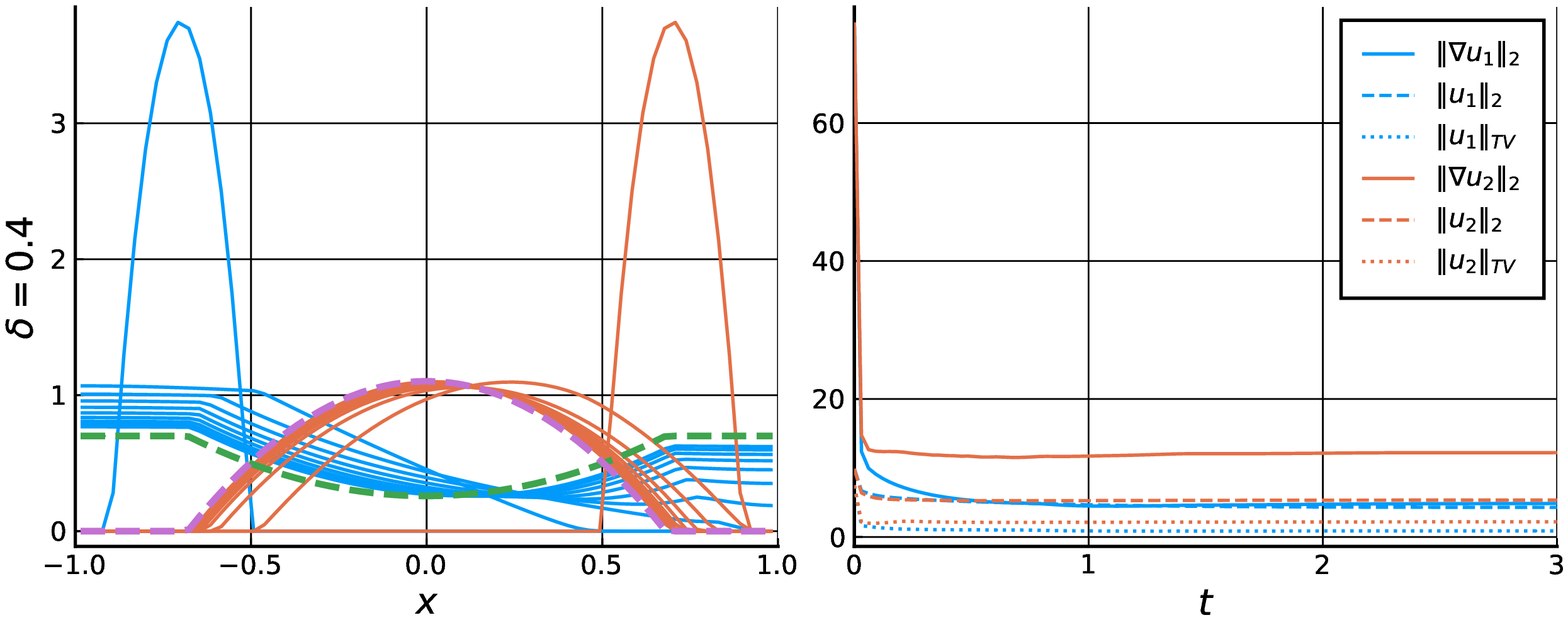}\\
	\includegraphics[width = 0.75\textwidth]{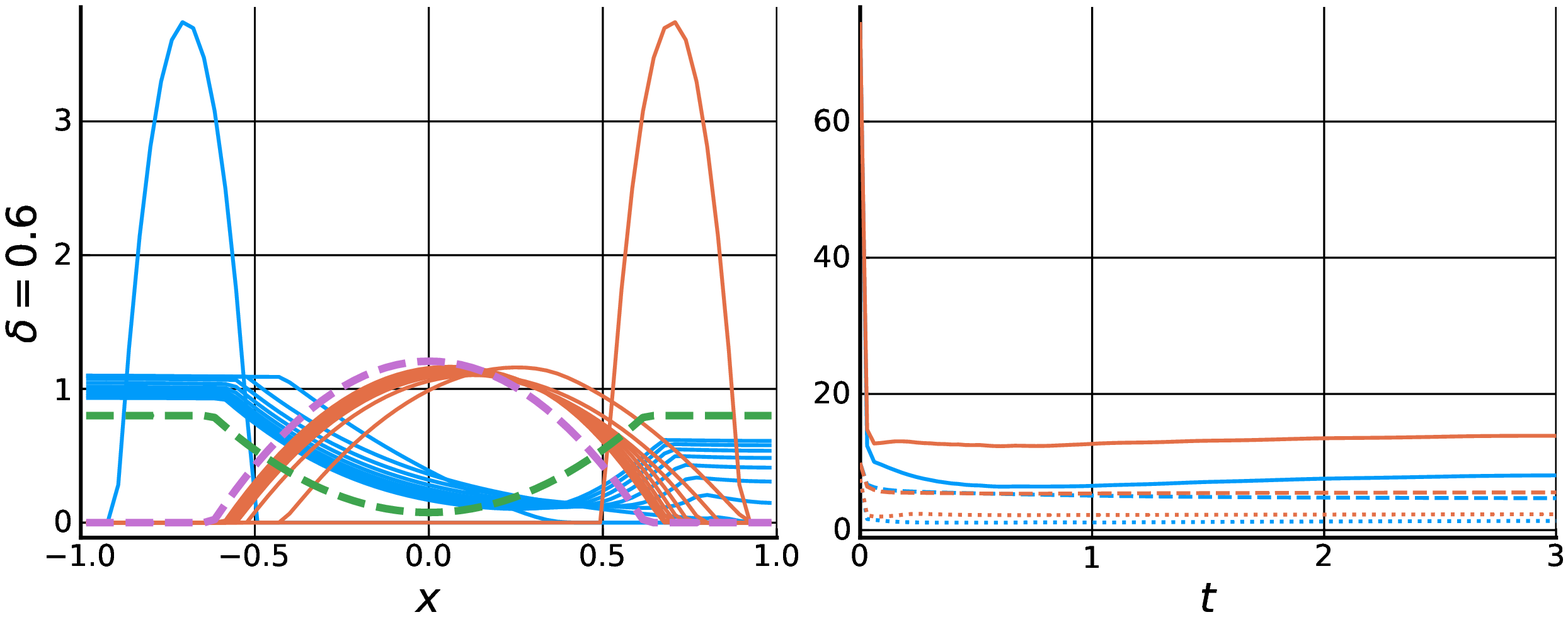}\\
	\includegraphics[width = 0.75\textwidth]{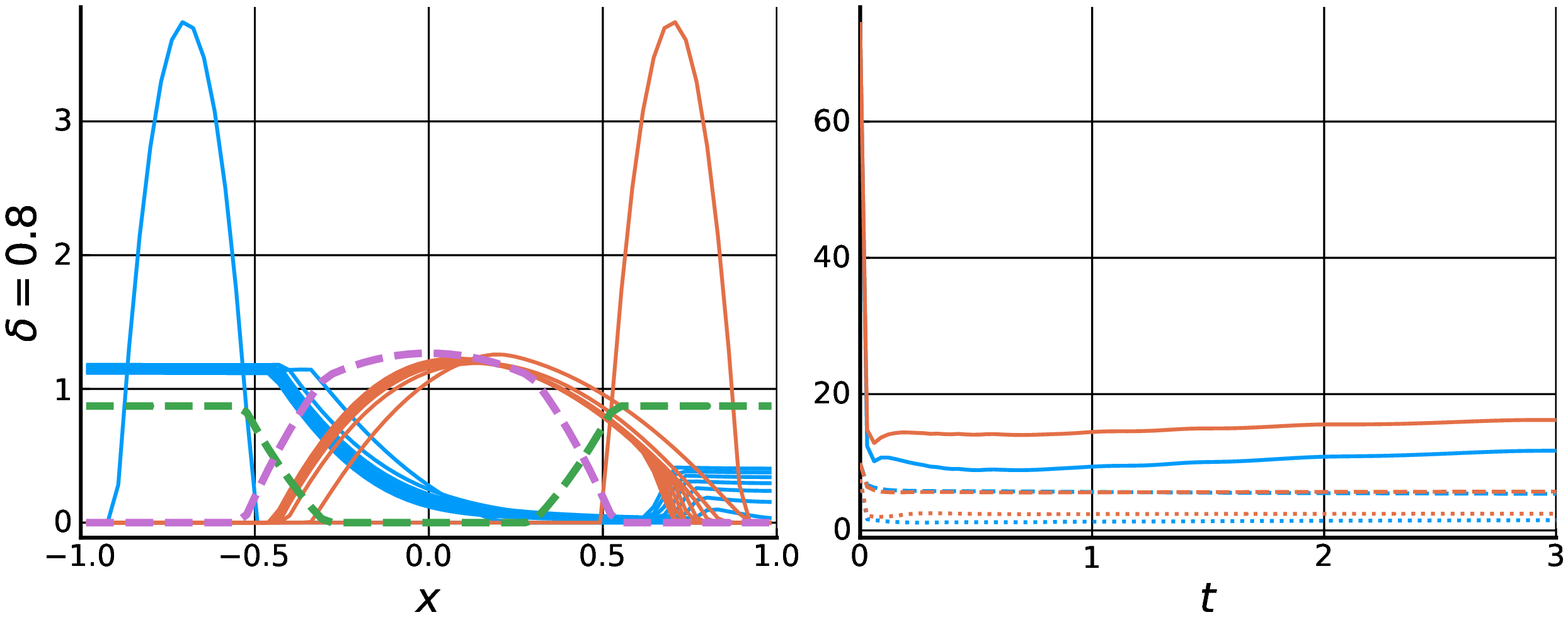}\\
	\includegraphics[width = 0.75\textwidth]{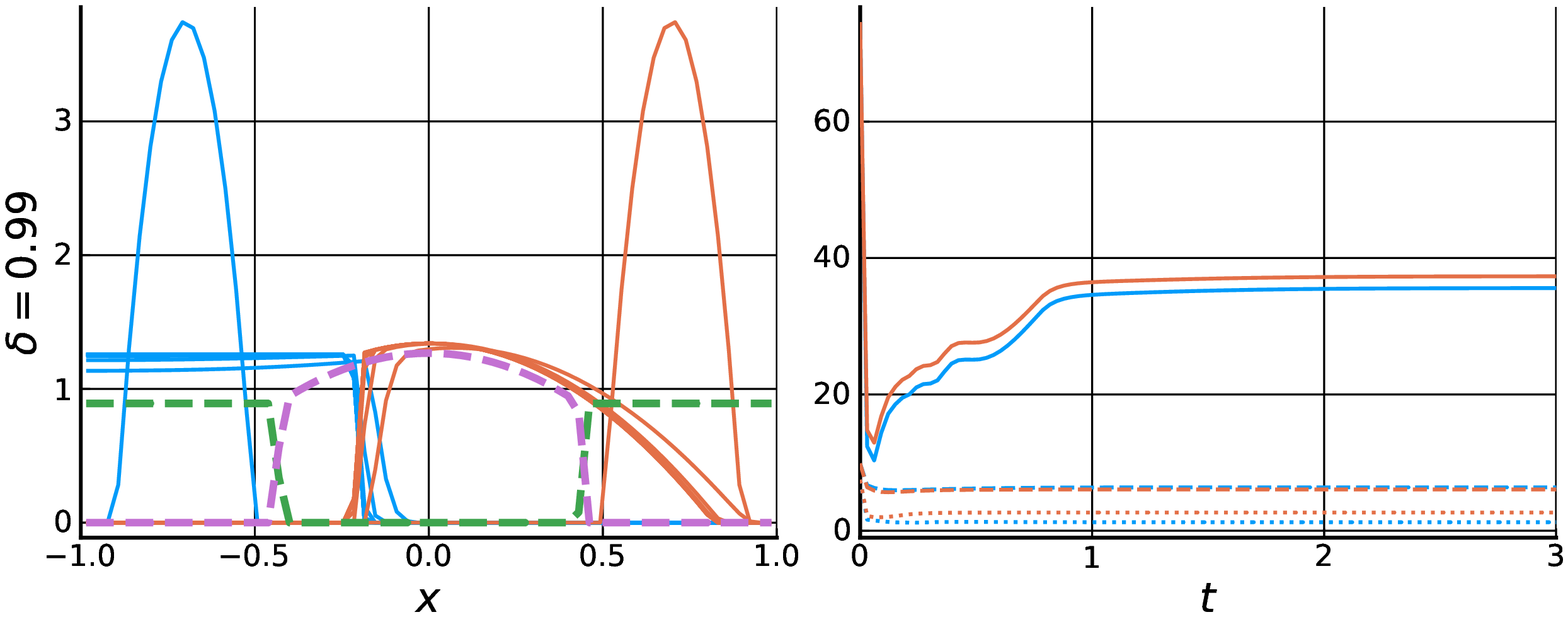}
\end{center}
\caption{Time evolution with left and right initial conditions, $V_1 = 0$ and $V_2 = 2x^2$ and final time $T = 3$, for various values of $\delta$ (Example \ref{ex:1}). The left column shows ten equally spaced timepoints in $[0,T]$ in solid lines, and the stationary stated in dashed lines. The right column shows the temporal evolution of three norms.}
 \label{fig:timeevol1}	
\end{figure}

In the right column of Fig.~\ref{fig:timeevol1} we plot the time evolution of the spatial $L^2$ norms of $u_i$ and $\nabla u_i$, as well as the Total Variation (TV), all computed using the partition given by the spatial grid used in the finite-volume scheme. The key point to note is that the effect of increasing $\delta$ is noticed markedly by the semi-norm $\|\nabla u_i\|_{L^2}$, whereas the other two norms, $\|u_i\|_{L^2}$ and $\|u_i\|_{TV}$, remain mostly unchanged by $\delta$. 
\end{ex}

\begin{ex}[Uniform initial conditions] \label{ex:2} Here we consider exactly the same set-up as in the previous example, except that now both components start with uniform initial conditions
$$
u_{1,0} (x) = u_{2,0}(x) = 1/|\Omega|.
$$
Therefore, we expect the same stationary states (since for $\delta <1$, $E$ is strictly convex). We show the results of this example in Fig.~\ref{fig:timeevol2}.
Because of the symmetry in the initial conditions, in this case the convergence to the steady state is much faster, as there is no mass that has to ``cross'' through the vacuum region as the latter is formed. 
\begin{figure}
\begin{center}
	\includegraphics[width = 0.75\textwidth]{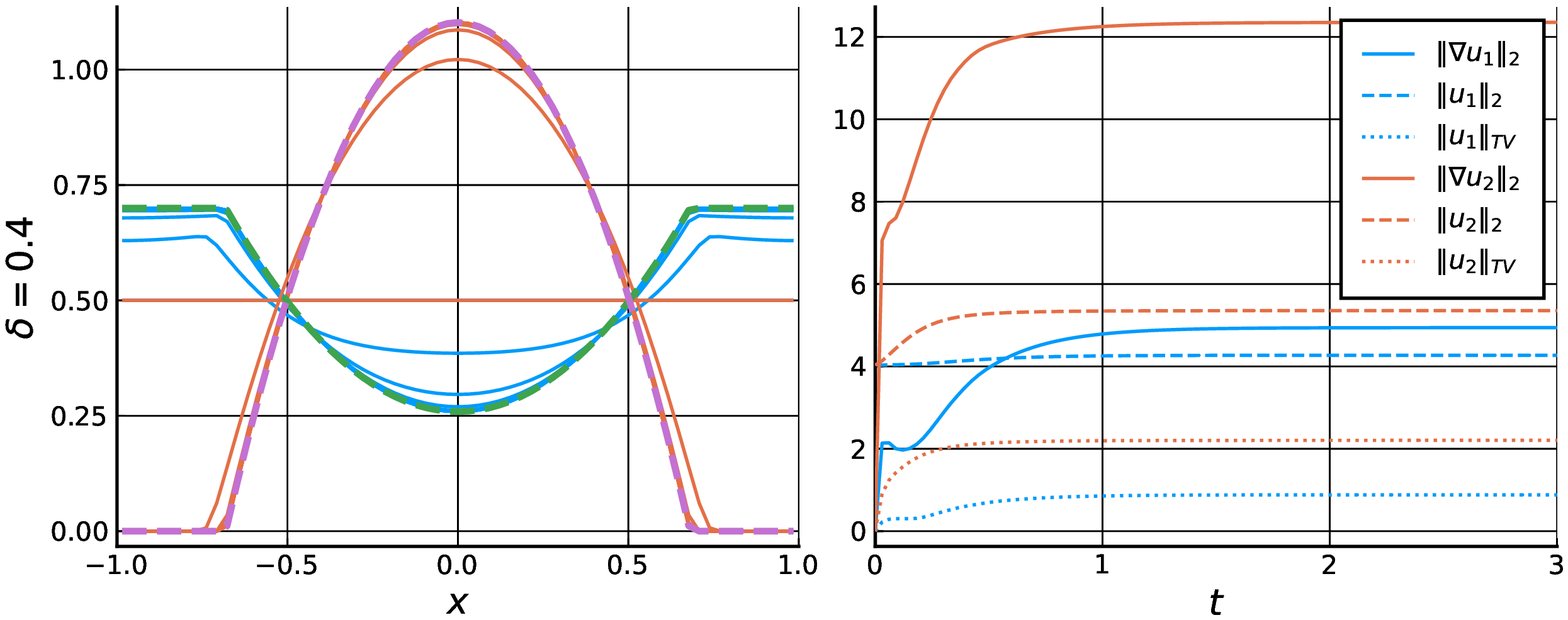}\\
	\includegraphics[width = 0.75\textwidth]{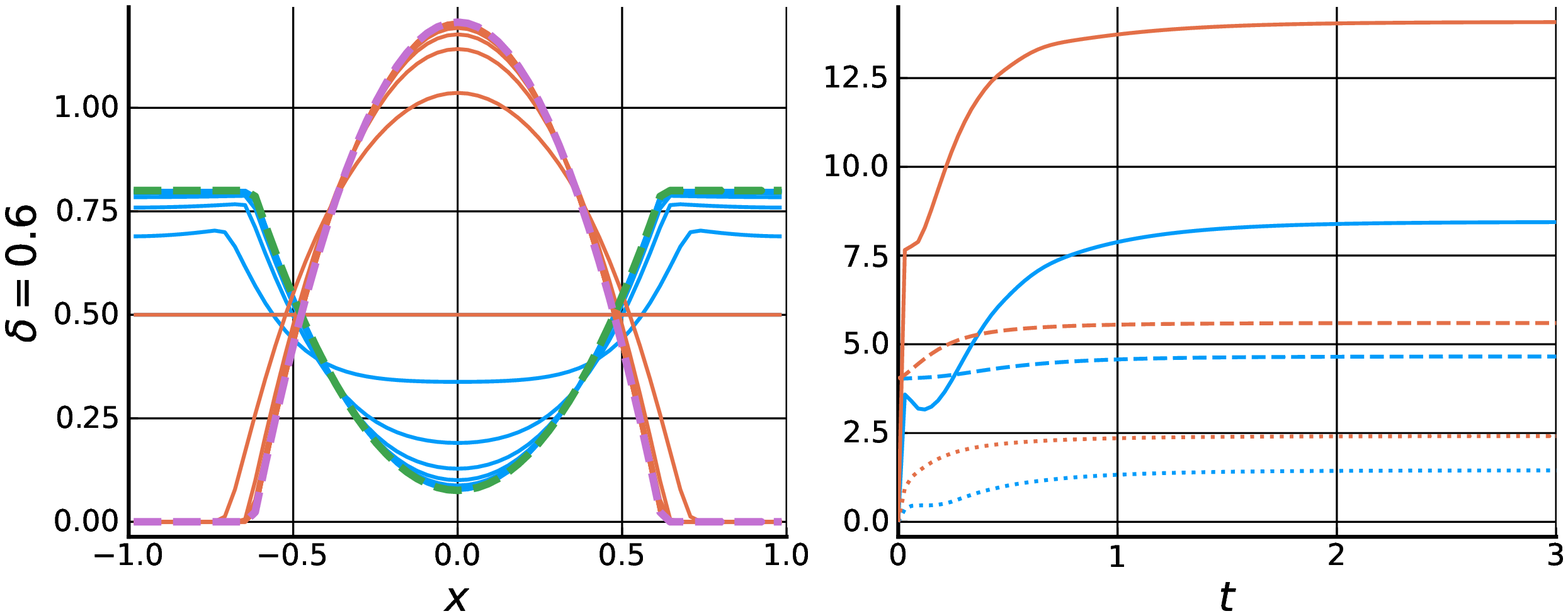}\\
	\includegraphics[width = 0.75\textwidth]{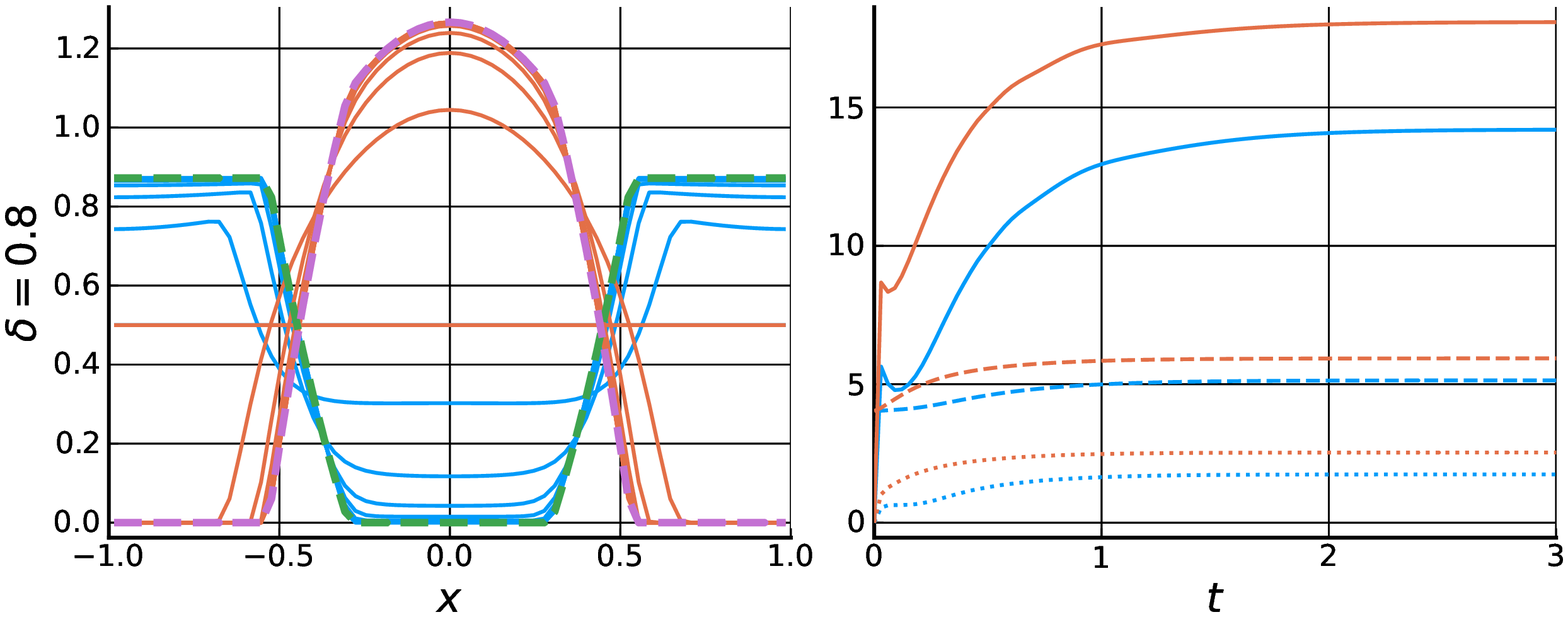}\\
	\includegraphics[width = 0.75\textwidth]{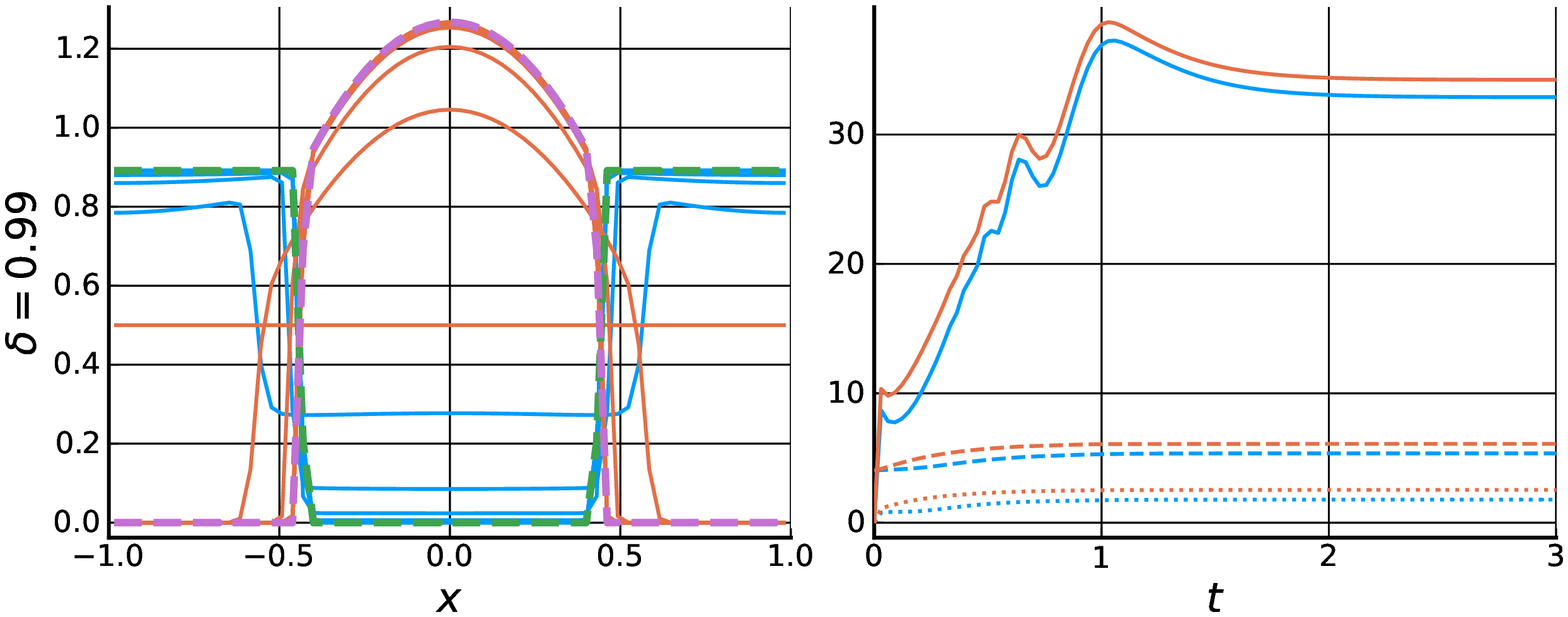}
\end{center}
\caption{Time evolution with uniform initial conditions, $V_1 = 0$ and $V_2 = 2x^2$ and final time $T = 3$, for various values of $\delta$ (Example \ref{ex:2}). The left column shows ten equally spaced timepoints in $[0,T]$ in solid lines, and the stationary stated in dashed lines. The right column shows the temporal evolution of three norms.
}
 \label{fig:timeevol2}	
\end{figure}	
\end{ex}

\begin{ex}[Stronger external potentials] \label{ex:3} We now consider the left and right initial conditions as in Example \ref{ex:1} while changing the external potentials to $V_1(x) = x^2/2$ and $V_2(x) = 50 x^2$. For this choice, we have that
$$
\delta_{\max}(T) = 0.00126, \quad \delta_{\max}(\Delta t) = 0.480,
$$
that is, a ten-fold reduction in $\delta_{\max}(T)$ with respect to Example \ref{ex:1} but a barely noticeable change in $\delta_{\max}(\Delta t)$ (as expected given that $\delta_{\max}$ is independent of $L_i$ in the small time limit).
The stronger confinement potential in the second species leads to a vacuum region in the first species for smaller values of $\delta$ than in the previous examples, and the associated slower convergence (see Figure \ref{fig:timeevol3}).
\begin{figure}[ht]
\begin{center}
\includegraphics[width = .95\textwidth]{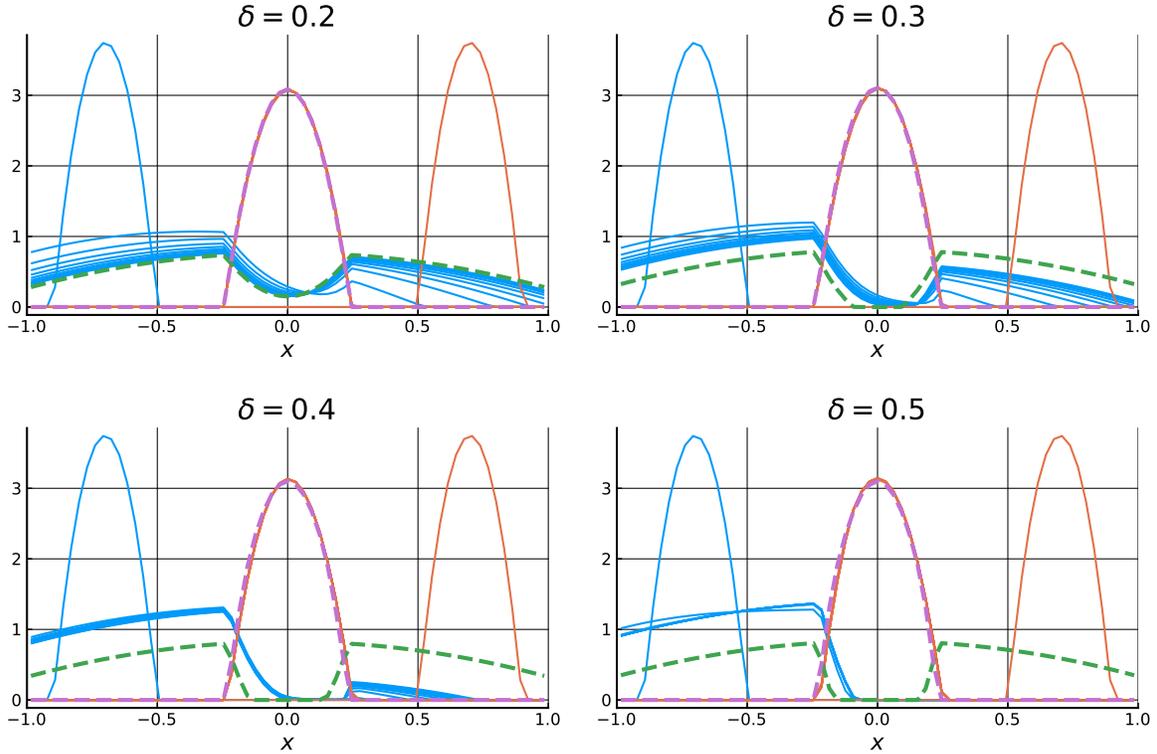}
\end{center}
\caption{Time evolution with uniform initial conditions, $V_1(x) = x^2/2$ and $V_2(x) = 50 x^2$ and final time $T = 3$, for various values of $\delta$ (Example \ref{ex:3}). There are ten equally spaced timepoints in $[0,T]$ in solid lines, and stationary stated in dashed lines.}
\label{fig:timeevol3}
\end{figure}	
\end{ex}

\begin{ex}[Evolution of norms in time and space as a function of $\delta$] \label{ex:4} In this final example, we look at the evolution of the norms in $\delta$ instead of time. To this end, we consider the following integrated-in-time norms
\begin{align*}
	\|u\|_{2,T} &= \left\{  \sum_{t_k = 0}^T \left[ \|u_1(t_k,\cdot)\|_{L^2}^2 + \|u_2(t_k,\cdot)\|_{L^2}^2  \right] \right\}^{1/2},\\
	\|\nabla u\|_{2,T} &= \left\{ \sum_{t_k = 0}^T \left[ \|\nabla u_1(t_k,\cdot)\|_{L^2}^2 + \|\nabla u_2(t_k,\cdot)\|_{L^2}^2   \right] \right\}^{1/2},\\
	\|u\|_{TV,T} &= \left\{ \sum_{t_k = 0}^T \left[ \|u_1(t_k,\cdot)\|_{TV} + \|u_2(t_k,\cdot)\|_{TV}   \right] \right\}^{1/2}.
\end{align*}

We use uniform initial conditions, a final time $T = 5$, and values for $\delta = 0, 0.1, \dots, 0.9, 0.95, 0.99$. We show the evolution of the three norms in Fig. \ref{fig:sweep_delta} for two cases: first, for the potentials used in Examples 1 and 2, namely $V_1 = 0$ and $V_2 = 2x^2$; and second, for $V_1 = x^2/2$ and $V_2 = 50x^2$. In the latter, the combination of external potentials makes the interface between the two components sharper (since the second component has a very strong confining potential, but also the first component now wants to concentrate around the origin). This fact is clearly visible in the trend of $\|\nabla u\|_{2,T}$ for increasing $\delta$. In contrast, the TV norm remains unchanged. 
As mentioned in the introduction, this observation indicates that a smallness assumption on $\delta$ is necessary in order to keep to the functional framework of $L^2(0,T;H^1(\Omega))$ for the analysis. The plot also illustrates that there is hope for a more general existence theory for solutions belonging to the space $BV$ when $\delta$ is close to $1$.

\begin{figure}
\begin{center}
	\includegraphics[width = .95\textwidth]{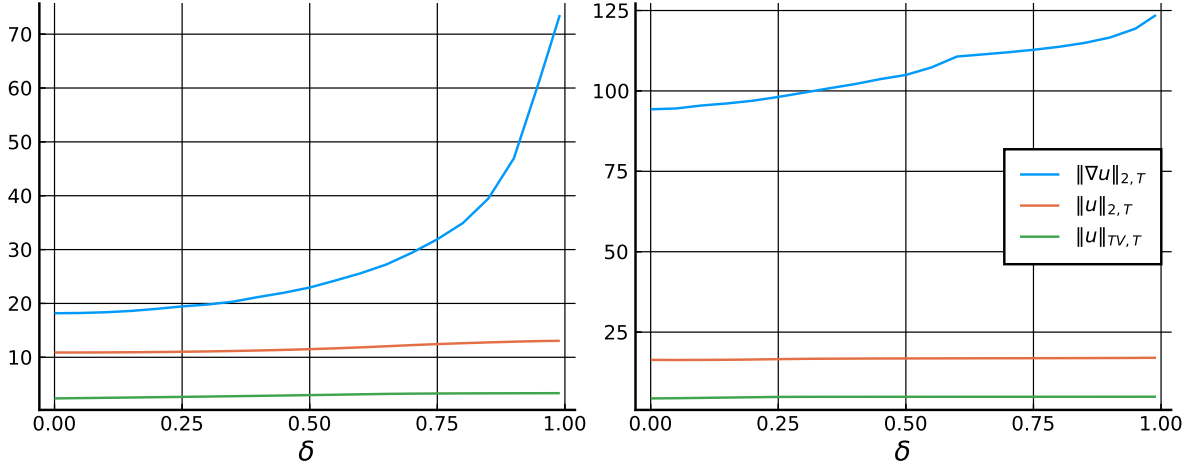}
	\end{center}
\caption{Norms (rescaled by $T |\Omega|$) in time and space as a function of $\delta$ with uniform initial conditions and final time $T = 3$. (Left) $V_1 = 0$ and $V_2 = 2x^2$. (Right) $V_1 = x^2/2$ and $V_2 = 50x^2$.}
 \label{fig:sweep_delta}	
\end{figure}
\end{ex}

\section{Regularised frozen system}\label{sec:regularised}
We introduce below the regularised system with frozen cross-diffusion. Let $\bar{z} = (\bar{z}_i)_{i=1}^M \in (C^\infty(\bar{Q}_T))^M$ be a given vector function. Throughout this section, we denote by $z=(z_i)_{i=1}^M$ the solution of the \emph{regularised frozen system} 
\begin{equation}\label{eq:pdefrozen}
\left\lbrace\begin{aligned}
& \partial_t z_i = \dv [z_i
(\nabla z_i + \nabla L_i(t,x,z_i) + \delta F_i(t,x,\bar{z},\nabla\bar{z})) + \eps \nabla z_i] \qquad &&\text{in } Q_T, \\ 
& 0 = \nu \cdot [z_i
(\nabla z_i + \nabla L_i(t,x,z_i) + \delta F_i(t,x,\bar{z},\nabla\bar{z})) + \eps \nabla z_i] \qquad && \text{on } \Sigma_T, \\ 
& z_i(0,\cdot) = z_{i,0} \qquad &&\text{on } \Omega. 
\end{aligned}\right.
\end{equation}

\begin{rem}\label{rem:fixing initial data for pdefrozen}
The constant $\varepsilon>0$ and the vector of non-negative functions $z_0 = (z_{i,0})_{i=1}^M \in (C^\infty_c(\Omega))^M$ do not change throughout the present section and Section \ref{sec:fixed_point}. The initial functions $z_{i,0}$ are chosen such that $\int_\Omega z_{i,0} \d x = \int_\Omega u_{i,0} \d x$ for $i\in\{1,\dots,M\}$. 
\end{rem}

In the next subsection we introduce the definition of weak solution to the above regularised frozen system, and show the mass conservation and non-negativity for such weak solutions. We then prove the existence of these solutions, and deduce from their regularity that they satisfy the system of equations in the classical sense. Then we prove some Sobolev estimates independent of $\varepsilon$, and conclude with a uniqueness result. 

In what follows, we will sometimes use the shorthand $\bar{F}_i$ to refer to the function 
\begin{equation}\label{eq:F bar notation}
    \bar{F}_i(t,x) := F_i(t,x,\bar{z}(t,x),\nabla\bar{z}(t,x)) \qquad \forall (t,x) \in \bar{Q}_T. 
\end{equation}

\begin{rem}\label{rem:Mbarz bound}
In view of the condition that $\{G^0_i\}_{i=1}^M$ and $\{G^1_{ij}\}_{i,j=1}^M$ in \eqref{eq:F structure 0} be $C^2$-regular and uniformly bounded with respect to all arguments (see Subsection \ref{sec:general}), for each $\bar{z}= (\bar{z}_i)_{i=1}^M \in (C^\infty(\bar{Q}_T))^M$ fixed, there exists a positive constant $\Lambda_{\bar{z}}$, where 
\begin{equation*}
\Lambda_{\bar{z}}=\Lambda_{\bar{z}}\left( \max_{1\leq i \leq M}\Vert \bar{z}_i \Vert_{C^2(\bar{Q}_T)} , \max_{1 \leq i \leq M}\Vert F_i \Vert_{C^2(\bar{Q}_T \times \mathbb{R}^M \times \mathbb{R}^{Md})}\right), 
\end{equation*}
such that, for all $(t,x) \in \bar{Q}_T$ and every $i \in \{1,\dots,M\}$, 
\begin{equation}\label{eq:boundedness of bar F z}
    \begin{aligned}
    |\bar{F}_i(t,x)| + \sum_{j=1}^d \left| \frac{\partial\bar{F}_i}{\partial x_j} (t,x) \right| + \left| \frac{\partial\bar{F}_i}{\partial t} (t,x) \right| + \sum_{j=1}^d\left| \frac{\partial^2 \bar{F}_i}{\partial t\partial x_j}(t,x) \right| + \sum_{k=1}^d \sum_{j=1}^d\left| \frac{\partial^2 \bar{F}_i}{\partial x_k \partial x_j}(t,x) \right| \leq \Lambda_{\bar{z}}. 
    \end{aligned}
\end{equation}

Additionally, there exists a positive constant $\Lambda_0$ depending only on $(z_{i,0})_{i=1}^M$ such that 
\begin{equation}\label{eq:M zero}
    \Vert z_{i,0} \Vert_{C^2(\bar{\Omega})} \leq \Lambda_0. 
\end{equation}
\end{rem}

\subsection{Definition and existence of regularised solutions}\label{subsec:regularised}

In accordance with the definition of weak solution given in \cite[Section 5.7]{ladyzhenskaia1988linear}, we provide the following notion of solution to the regularised frozen problem. 

\begin{defi}[Weak solution for regularised frozen system]\label{defi:weak sol regularised}
We say that $z \in (C^{2,1}(\bar{Q}_T))^M$  \emph{solves the weak form} of \eqref{eq:pdefrozen} if, for any test function $\phi \in C^1(\bar{Q}_T)$, for $i \in \{1,\dots,M\}$, for $t \in [0,T]$, 
\begin{equation}\label{eq:weak sol regularised}
    \begin{aligned}
        \int_\Omega z_i(t,x) \phi(t,x) \d x &- \int_\Omega z_{i,0}(x) \phi(0,x) \d x - \int_0^t \int_{\Omega} z_i \partial_t\phi \d x \d t \\ 
        &+\int_0^t \int_{\Omega} \big[ z_i(\nabla z_i + \nabla L_i(t,x,z_i) + \delta F_i(t,x,\bar{z},\nabla\bar{z})) + \varepsilon \nabla z_i \big]\cdot\nabla\phi \d x \d t = 0. 
    \end{aligned}
\end{equation}
Correspondingly, we define $S_\varepsilon : \bar{z} \to z$ to be the \emph{solution operator}, whose image is the weak solution of \eqref{eq:pdefrozen}.
\end{defi}

\begin{rem}
Note that the following compatibility condition has been implicitly imposed in the previous weak formulation, 
\begin{equation}\label{eq:compatibility 0}
    0 = \nu \cdot \big[ z_{i,0}(\nabla z_{i,0} + \nabla L_i(0,x,z_{i,0}) + \delta F_i(0,x,\bar{z}(0,x),\nabla\bar{z}(0,x))) + \varepsilon \nabla z_{i,0} \big] \qquad \text{on } \partial\Omega, 
\end{equation}
which is manifestly satisfied for all choices of $\bar{z} \in (C^\infty(\bar{Q}_T))^M$, as the fixed initial data $z_{i,0}$ is identically zero on the boundary $\partial\Omega$ due to its compact support (see Remark \ref{rem:fixing initial data for pdefrozen}). 
\end{rem}

\begin{lem}[Mass conservation for regularised frozen system]\label{lem:mass cons regularised}
Suppose there exists a weak solution $z$ of the problem \eqref{eq:pdefrozen} in the sense of Definition \ref{defi:weak sol regularised}. Then, for $i\in\{1,\dots,M\}$, 
\begin{equation*}
    \int_\Omega z_i(t,x) \d x = \int_\Omega z_{i,0}(x) \d x \qquad \forall t \in [0,T]. 
\end{equation*}
\end{lem}
\begin{proof}
The assertion is immediate from using the test function $\phi = 1$ in Definition \ref{defi:weak sol regularised}. 
\end{proof}

\begin{lem}[Sign preservation for regularised frozen system]\label{lem:sign preserved regularised}
Suppose there exists a weak solution $z$ of the problem \eqref{eq:pdefrozen} in the sense of Definition \ref{defi:weak sol regularised}. Then, for $i \in \{1,\dots,M\}$, 
\begin{equation*}
    z_i(t,x) \geq 0 \qquad \text{for a.e. } (t,x) \in \bar{Q}_T. 
\end{equation*}
\end{lem}
\begin{proof}
Define the function $\theta(t,x) := [z_i(t,x)]_-$ to be the negative part of the weak solution in question. 
Noting that $\theta = -z_i \mathds{1}_{z_i \leq 0}$, we observe that this function is non-negative and supported in the set 
$\{(t,x) \in \bar{Q}_T : z_i(t,x) \leq 0\}$. Moreover, we find that 
\begin{equation*}
    \nabla \theta = -\nabla z_i \mathds{1}_{z_i \leq 0}, \qquad \partial_t \theta = -\partial_t z_i \mathds{1}_{z_i \leq 0}
\end{equation*}
in the sense of distributions. It follows that $\theta \in L^2(0,T;H^1(\Omega)) \cap L^\infty(0,T;L^1(\Omega))$ and $\partial_t \theta \in X' \cap L^2(0,T;(H^1(\Omega))')$. Using standard density arguments in Sobolev spaces, we may test against $\theta$ in the weak formulation of Definition \ref{defi:weak sol regularised}. In turn, we obtain, for a.e.~$t\in (0,T)$, 
\begin{equation}\label{eq:test against negative part prior to time integral}
    \frac{d}{dt}\int_\Omega \frac{1}{2}\theta^2(t,x) \d x + \int_\Omega \big( \theta |\nabla\theta|^2 + \theta \nabla \theta \cdot \nabla L_i(t,x,z_i) + \delta \theta \nabla\theta\cdot F_i(t,x,\bar{z},\nabla\bar{z}) + \varepsilon |\nabla\theta|^2 \big) \d x = 0. 
\end{equation}
Given the form of the terms $\{L_i\}_{i=1}^M$ from \eqref{eq:Li terms}, we have 
\begin{equation*}
    \bigg| \int_\Omega \theta \nabla \theta \cdot \nabla L_i(t,x,z_i) \d x \bigg| \leq \big( \Vert \nabla V_i(t,\cdot) \Vert_{L^\infty(\Omega)} + \Vert \nabla W_i*z_i(t,\cdot) \Vert_{L^\infty(\Omega)} \big) \int_\Omega \theta |\nabla \theta| \d x, 
\end{equation*}
and 
\begin{equation*}
    \Vert \nabla W_i*z_i(t,\cdot) \Vert_{L^\infty(\Omega)} \leq \Vert \nabla W_i \Vert_{L^\infty(\mathbb{R}^d)} \Vert z_i(t,\cdot) \Vert_{L^1(\Omega)} \qquad \text{a.e.~} t \in (0,T). 
\end{equation*}
It therefore follows, using the fact that $z_i \in L^\infty(0,T;L^1(\Omega))$ since $z_i \in C^{2,1}(\bar{Q}_T)$ as per Definition \ref{defi:weak sol regularised}, that 
\begin{equation*}
    \bigg| \int_\Omega \theta \nabla \theta \cdot \nabla L_i(t,x,z_i) \d x \bigg| \leq C_L \big( 1 + \Vert z_i \Vert_{L^\infty(0,T;L^1(\Omega))} \big) \int_\Omega \theta |\nabla \theta| \d x \qquad \text{a.e.~} t \in (0,T). 
\end{equation*}
Meanwhile, using the boundedness of $\bar{z}$ and that of $F$ in $C^2$ to control $F_i(t,x,\bar{z},\nabla\bar{z})$ from \eqref{eq:boundedness of bar F z} (see Remark \ref{rem:Mbarz bound}), we obtain 
\begin{equation*}
    \bigg| \int_\Omega \delta \theta \nabla \theta \cdot F_i(t,x,\bar{z},\nabla\bar{z}) \d x \bigg| \leq |\delta|\Lambda_{\bar{z}} \int_\Omega \theta |\nabla \theta| \d x \qquad \text{a.e.~} t \in (0,T). 
\end{equation*}

Integrating \eqref{eq:test against negative part prior to time integral} with respect to the time variable, and using the previous estimates, we find
\begin{equation*}
    \begin{aligned}
        \int_\Omega \frac{1}{2}\theta^2(t,x) \d x + \int_0^t \int_\Omega \theta |\nabla \theta|^2 \d x \d \tau + \varepsilon \int_0^t \int_\Omega |\nabla \theta|^2 \d x \d \tau \leq & \Lambda \int_0^t \int_\Omega \theta |\nabla \theta| \d x, 
    \end{aligned}
\end{equation*}
where the positive constant 
\begin{equation*}
    \Lambda := |\delta| \Lambda_{\bar{z}} + C_L \big( 1 + \Vert z_i \Vert_{L^\infty(0,T;L^1(\Omega))} \big) 
\end{equation*} 
is independent of time. An application of the Cauchy--Young inequality gives 
\begin{equation*}
    \begin{aligned}
        \int_\Omega \frac{1}{2}\theta^2(t,x) \d x + \int_0^t \int_\Omega \theta |\nabla \theta|^2 \d x \d \tau + \frac{\varepsilon}{2}\int_0^t & \int_\Omega |\nabla \theta|^2 \d x \d \tau \leq \frac{\Lambda^2 }{\varepsilon} \int_0^t \int_\Omega \frac{1}{2}\theta^2 \d x \d \tau.
    \end{aligned}
\end{equation*}
Dropping the last two terms in the left-hand side of the inequality above, Gr\"{o}nwall's Lemma yields 
\begin{equation*}
    \int_\Omega \theta^2(t,x) \d x \leq \bigg( \int_\Omega \theta^2(0,x) \d x \bigg) e^{\frac{\Lambda^2}{\varepsilon}t}
    = 0 \qquad \text{for a.e. } t \in (0,T), 
\end{equation*}
where the final equality follows from the non-negativity of the initial data $z_{i,0}$ (Remark \ref{rem:fixing initial data for pdefrozen}). The result follows. 
\end{proof}

\begin{rem}
It is a priori not clear how to prove such a sign preservation result for the original system \eqref{eq:main0} directly from Definition \ref{defi:weak sol 0}, due to the presence of cross-terms of the form $\int_\Omega \nabla u_j(t,x) [u_i(t,x)]_- \d x$ with $i \neq j$. The non-negativity of the solution of the original system \eqref{eq:main0} will therefore be deduced via a limiting procedure from the non-negativity of the regularised solutions of \eqref{eq:pdefrozen}. 
\end{rem}

In what follows, we apply the classical theory of Ladyzhenskaya, Solonnikov, and Uraltseva \cite[Chap.~5, Sec.~7, Thm.~7.4]{ladyzhenskaia1988linear} to deduce the existence and uniqueness of classical solutions to the regularised system \eqref{eq:pdefrozen}. The proof is given in Appendix \ref{sec:appendix proof russian}. 
\begin{lem}[Existence and uniqueness of regularised solutions]\label{lem:existence of viscous approximates 0}
There exists a unique  $z = (z_i)_{i=1}^M \in C^{2,1}(\bar{Q}_T)$ solving \eqref{eq:pdefrozen} as a pointwise equality between continuous functions. Moreover, for $i \in \{1,\dots,M\}$, 
\begin{equation*}
    \begin{aligned}
        z_i(t,x) \geq 0 \quad \text{for } (t,x) \in \bar{Q}_T, \quad z_i(0,x) = z_{i,0}(x) \quad \text{for } x \in \bar{\Omega}, \quad \int_\Omega z_i(t,x) \d x = \int_\Omega z_{i,0}(x) \d x  \quad \text{for } t \in [0,T]. 
    \end{aligned}
\end{equation*}
\end{lem}

\subsection{Uniform estimates}\label{subsec:uniform}

In this subsection, we derive the uniform estimates on the solutions of the regularised frozen system by testing the equation against the logarithm of the solution.

\begin{lem}[Energy estimates]\label{lem:energy estimates 0}
Let $z = (z_i)_{i=1}^M \in (C^{2,1}(\bar{Q}_T))^M$ be the solution of \eqref{eq:pdefrozen} provided by Lemma \ref{lem:existence of viscous approximates 0}. Then, for any $t \in [0,T]$, there holds the estimate 
\begin{multline}\label{eq:first energy est with u log u 0}
\int_\Omega z_i(t)  (\log z_i(t))_+ \d x + \frac{1}{2} \int_0^t \int_\Omega |\nabla z_i|^2 \d x \d \tau + \varepsilon\int_0^t \int_\Omega  \frac{|\nabla z_i|^2}{z_i} \d x \d \tau 
\\
\leq (e^{-1}|\Omega|+2|\Omega| T C_L^2) + \int_\Omega z_{i,0} \log z_{i,0} \d x + 2|\Omega| T C_L^2 \left( \int_\Omega z_{i,0} \d x \right)^2 
+  \int_0^T \int_\Omega \delta^2 |\bar{F}_i|^2 \d x \d \tau, 
\end{multline}
for $i \in \{1,\dots,M\}$, along with the Sobolev estimate 
\begin{equation}\label{eq:L2H1 bound 0}
\begin{aligned}
    \sup_{t \in [0,T]} \int_\Omega z_i(t) (\log z_i(t))_+ \d x + \Vert z_i &\Vert_{L^2(0,T;H^1(\Omega))}^2 \\ 
    &\leq
C_\Omega \bigg( 1 + \Vert z_{i,0} \Vert_{L^1(\Omega)}^2 + \int_\Omega z_{i,0}\log z_{i,0} \d  x
+ \int_0^T\int_\Omega
\delta^2 |\bar{F}_i|^2
\d x \d t \bigg), 
\end{aligned}
\end{equation}
for $i\in\{1,\dots,M\}$. The positive constant $C_\Omega = C_\Omega(\Omega,T,d)$, which is independent of $\varepsilon$ and $z_0 = (z_{i,0})_{i=1}^M$, is given by 
\begin{equation}\label{eq:big constant for L2H1 bound 0}
    C_\Omega = 2 \cdot \max\left\lbrace (1+C_P) , \frac{T}{|\Omega|} + 2 |\Omega| (1+C_P)\big( e^{-1} + 2 T C_L^2 \big) \right\rbrace, 
\end{equation}
where $C_P = C_P(\Omega,d)$ is the Poincar\'{e} constant. 
\end{lem}
\begin{proof}
Begin by assuming that $z_i$ is strictly positive in $Q_T$ and multiply \eqref{eq:pdefrozen} by $\log z_i$. By rewriting the time derivative, we get 
\begin{equation*}
    \partial_t (z_i \log z_i) - \partial_t z_i = (\log z_i) \dv [z_i(\nabla z_i + \nabla L_i(t,x,z_i) + \delta \bar{F}_i) + \varepsilon\nabla z_i] \qquad \forall (t,x) \in Q_T. 
\end{equation*}
Integrating in space and time then yields, for any $t \in [0,T]$, 
\begin{equation}\label{eq:precede a priori}
    \begin{aligned}
    \int_0^t \frac{d}{d\tau}\bigg(\int_\Omega z_i(\tau) \log z_i(\tau&) \d x \bigg) \d \tau - \int_0^t \frac{d}{d\tau}\bigg( \int_\Omega z_i(\tau) \d x  \bigg) \d \tau = \\ 
    &- \int_0^t \int_\Omega  \left[ \nabla z_i \cdot (\nabla z_i + \nabla L_i(t,x,z_i) + \delta \bar{F}_i) + \varepsilon \frac{|\nabla z_i|^2}{z_i} \right] \d x \d \tau, 
    \end{aligned}
\end{equation}
where the no-flux boundary condition makes the boundary term vanish.

When $z_i$ is merely non-negative, we multiply by $\log(\beta+z_i)$ for some $\beta>0$ and get 
\begin{equation*}
    \partial_t (z_i \log (\beta+ z_i)) - (\partial_t z_i)\frac{z_i}{\beta+z_i} = (\log (\beta+z_i)) \dv [z_i(\nabla z_i + \nabla L_i(t,x,z_i) + \delta \bar{F}_i) + \varepsilon\nabla z_i], \qquad \forall (t,x) \in Q_T.
\end{equation*}
Integrating in time and space, and using that $|\partial_t z_i|$ is integrable on $Q_T$ since $z_i$ is $C^{2,1}(\bar{Q}_T)$, the Dominated Convergence Theorem implies 
$$ \lim_{\beta \to 0} \int_0^t \int_\Omega (\partial_t z_i)\frac{z_i}{\beta+z_i} \d x \d \tau = \int_0^t \int_\Omega (\partial_t z_i) \d x \d \tau .$$
All other terms may be treated similarly with the exception of the Fisher information, where one obtains 
\begin{equation*}
    \lim_{\beta \to 0} \int_0^t \int_\Omega \frac{|\nabla z_i|^2}{\beta +  z_i} \d x \d \tau = \int_0^t \int_\Omega \frac{|\nabla z_i|^2}{z_i} \d x \d \tau,
\end{equation*}
by the Monotone Convergence Theorem, since the sequence of integrands $\left \{|\nabla z_i|^2/(\beta+z_i) \right \}_{\beta>0}$ is pointwise increasing as $\beta \to 0$ and non-negative by virtue of Lemma \ref{lem:existence of viscous approximates 0}. In turn, we recover \eqref{eq:precede a priori}. 

Recall from Lemma \ref{lem:existence of viscous approximates 0} that  $\int_\Omega z_i(\tau) \d x$ is constant. Equation \eqref{eq:precede a priori} therefore simplifies to 
\begin{equation}\label{eq:first energy est prelim 0}
    \begin{aligned}
    \int_\Omega z_i(t) \log z_i(t) \d x  = \int_\Omega z_{i,0} \log z_{i,0} \d x  - \delta \int_0^t \int_\Omega \nabla z_i \cdot \bar{F}_i \d x \d \tau &- \int_0^t \int_\Omega  |\nabla z_i|^2 \d x \d \tau \\ 
    - \int_0^t \int_\Omega \nabla z_i \cdot \nabla L_i(t,x,z_i) \d x \d \tau &- \varepsilon\int_0^t \int_\Omega  \frac{|\nabla z_i|^2}{z_i} \d x \d \tau, 
    \end{aligned}
\end{equation}
for any $t \in [0,T]$. 
Observe also that 
$ z \log z = z \log z \mathds{1}_{\{z \geq 1\}} + z \log z \mathds{1}_{\{0 \leq z < 1\}}$,
$ \forall z \in [0,\infty),
$
and, since $x \mapsto x \log x$ is non-positive over the interval $[0,1]$ and achieves its minimum (with value $-e^{-1}$) at the point $x = e^{-1}$, it follows that 
\begin{equation}\label{eq:positive part log lower bound}
    z \log z \geq z \log z \mathds{1}_{\{z \geq 1\}} - e^{-1} \qquad \forall z \in [0,\infty). 
\end{equation}
Also, using \eqref{eq:Li terms}, \eqref{eq:drift V bound}, and an application of the triangle inequality followed by H\"{o}lder's inequality yields 
\begin{align} \label{eq:holder grad z grad L}
\begin{aligned}
\bigg| \int_0^t \int_\Omega \nabla z_i \cdot \nabla L_i(t,x,z_i) \d x \d \tau \bigg| &\leq \Vert \nabla V_i \Vert_{L^2(Q_T)} \Vert \nabla z_i \Vert_{L^2(Q_T)} \\ 
        &\qquad + \Vert \nabla z_i \Vert_{L^2(Q_T)} \bigg( \int_{0}^t \int_\Omega |\nabla W_i*z_i(\tau,\cdot)(x)|^2 \d x \d \tau  \bigg)^{\frac{1}{2}}\\
&  \leq (|\Omega|T)^{\frac{1}{2}} C_L \left( 1 + \int_\Omega  z_{i,0}(y) \d y \right) \Vert \nabla z_i \Vert_{L^2(Q_T)}, 
\end{aligned}
\end{align}
where we bounded the convolution term as follows
\begin{equation}\label{eq:NO YOUNG}
\begin{aligned}
 \int_{0}^t \int_\Omega |\nabla W_i*z_i(\tau,\cdot)(x)|^2 \d x \d \tau  
 &= \int_{0}^t \int_\Omega \bigg|\int_\Omega \nabla W_i(x-y)z_i(\tau,y) \d y\bigg|^2 \d x \d \tau \\ 
&\leq C_L^2 \int_0^t \int_\Omega \bigg( \int_\Omega z_i(\tau , y) \d y \bigg)^2 \d x \d \tau = C_L^2 |\Omega| T \left( \int_\Omega  z_{i,0}(y) \d y \right)^2,  
    \end{aligned}
\end{equation}
where we used the boundedness of $\{\nabla W_i\}_{i=1}^M$ inherited from \eqref{eq:drift V bound}, the non-negativity of $z_i$ due to Lemma \ref{lem:sign preserved regularised}, and the fact that $x \mapsto x^2$ is increasing on $[0,\infty)$ to obtain the inequality, and the mass conservation from Lemma \ref{lem:mass cons regularised} to obtain the final equality. 
Using estimates \eqref{eq:positive part log lower bound} and \eqref{eq:holder grad z grad L}, along with an application of the weighted Cauchy--Young inequality to the terms on the right-hand side of \eqref{eq:holder grad z grad L} and to $\int_0^t \int_\Omega \nabla z_i \cdot \bar{F}_i \d x \d t$, we deduce the estimate \eqref{eq:first energy est with u log u 0} from \eqref{eq:first energy est prelim 0}. 

The Poincar\'{e}--Wirtinger inequality 
\begin{equation*}
\int_\Omega \left(z_i - \frac{1}{|\Omega|}\int_\Omega z_i \d x \right)^2 \d x
\leq C_P \int_\Omega |\nabla z_i|^2 \d x,
\end{equation*}
where $C_P = C_P(\Omega,d)$ is the Poincar\'{e} constant, implies
\begin{equation*}
\begin{aligned}
C_P \int_\Omega |\nabla z_i|^2 \d x  &\geq \int_\Omega \left(z_i - \frac{1}{|\Omega|}\int_\Omega z_i \d x \right)^2 \d x \\ 
&=
\int_\Omega \left[
z_i^2 - \frac{2}{|\Omega|}z_i \int_\Omega z_i \d x + \frac{1}{|\Omega|^2} \left( \int_\Omega z_i \d x \right)^2 \right] \d x = \int_\Omega
z_i^2 \d x
- \frac{1}{|\Omega|} \left( \int_\Omega z_{i,0} \d x \right)^2,
\end{aligned}
\end{equation*}
where we used the conservation of mass in the final equality.
Substituting back into \eqref{eq:first energy est with u log u 0}, we get 
\begin{multline*}
\int_0^T \int_\Omega \frac{1}{2} z_i^2 \d x \d t \leq \frac{T}{2|\Omega|}\left( \int_\Omega z_{i,0} \d x \right)^2  + C_P(e^{-1}|\Omega|+2|\Omega| T C_L^2) + C_P\int_\Omega z_{i,0} \log z_{i,0} \d x 
\\ + 2C_P|\Omega| T C_L^2 \left( \int_\Omega z_{i,0} \d x \right)^2 + C_P\int_0^T \int_\Omega \delta^2 |\bar{F}_i|^2 \d x \d \tau. 
\end{multline*}
Combining the above with \eqref{eq:first energy est with u log u 0}, we obtain 
\begin{equation*}
\begin{aligned} \int_0^T \int_\Omega \frac{1}{2} (z_i^2 + |\nabla z_i|^2 ) \d x \d t \leq \  &\left(\frac{T}{2|\Omega|} + 2(1+C_P)|\Omega| T C_L^2 \right)\left( \int_\Omega z_{i,0} \d x \right)^2 \\ 
&+ (1+C_P)(e^{-1}|\Omega|+2|\Omega| T C_L^2) + (1+C_P) \int_\Omega z_{i,0} \log z_{i,0} \d x \\ 
        &+ (1+C_P) \int_0^T \int_\Omega \delta^2 |\bar{F}_i|^2 \d x \d t, 
    \end{aligned}
\end{equation*}
from which we recover \eqref{eq:L2H1 bound 0} with the appropriate constant $C_\Omega$ given by \eqref{eq:big constant for L2H1 bound 0}. 
\end{proof}

Before proceeding to the next lemma, which covers the uniform estimate for the time derivative of the regularised solutions, we recall the interpolation result of Di Benedetto. 
\begin{prop}[Proposition 3.2 of \cite{dibenedetto2012degenerate}]\label{prop:DiBenedetto interpolation est 0}
Let $m,p\geq 1$ and assume that $\de\Omega$ is piecewise smooth. There exists a constant $\gamma$ depending only on $d, m, p$ and the structure of $\de\Omega$ such that, for any $v\in V^{m.p}$, where
\begin{equation}\label{eq:V12}
V^{m.p} := L^\infty(0,T;L^m(\Omega)) \cap L^p(0,T;W^{1,p}(\Omega)),
\end{equation}
we have
\begin{equation}
\nm{v}_{L^q(Q_T)}
\leq \gamma\left(1 + \frac{T}{|\Omega|^{\frac{d(p-m) + pm }{dm} } }\right)
\nm{v}_{V^{m.p}},
\qquad \text{where } q = p\frac{d+m}{d}. 
\end{equation}
\end{prop}

\begin{lem}[Time derivative estimate of the regularised solutions]\label{lem:time compactness bound 0}
Recall the space $X$ introduced in Definition \ref{defi:functional space}. There holds the uniform estimate in the dual space 
\begin{equation}\label{eq:unif est on dt u 0}
    \begin{aligned}
        \nm{\de_t z_i}_{X'} \leq C_{X'} \bigg( 1 + \Vert z_{i,0} \Vert_{L^1(\Omega)}^2 + \int_\Omega z_{i,0} \log z_{i,0} \, \d x + \delta^2 \nm{\bar{F}_i}_{L^2(Q_T)}^2 \bigg), 
    \end{aligned}
\end{equation}
where the positive constant $C_{X'}=C_{X'}(\Omega,T,d)$, which is independent of $\varepsilon$ and $z_0 = (z_{i,0})_{i=1}^M$, is given by 
\begin{equation}\label{eq:CXprime expression}
    C_{X'} = 2 \gamma \left[\left(1 + \frac{T}{|\Omega|^{\frac{d+2}{d} } }\right)  + (|\Omega|T)^{\frac{d}{r}} \right] \left( 2+3C_\Omega + 2|\Omega|TC_L^2 \right). 
\end{equation}
\end{lem}
\begin{proof}
Applying Proposition \ref{prop:DiBenedetto interpolation est 0} with $v=z_i$, $m=1$, and $p=2$ yields 
\begin{equation}\label{eq:apply dibenedetto 0}
\nm{z_i}_{L^q(Q_T)}
\leq \gamma\left(1 + \frac{T}{|\Omega|^{\frac{d+2}{d} } }\right)
\nm{z_i}_{V^{1.2}},
\quad q = 2\frac{d+1}{d}, 
\end{equation}
where the space $V^{1,2}$ is defined in \eqref{eq:V12}. Notice that $q>2$ for all choices of dimension $d$. 

Fix $\eta \in C^\infty(\bar{Q}_T)$. Going back to \eqref{eq:pdefrozen}, writing $\langle \de_t z_i, \eta \rangle = \int_0^T \int_\Omega \partial_t z_i \eta \d x \d t$, and using the divergence theorem in conjunction with the no-flux boundary condition, we find 
\begin{equation}\label{eq:once this is established}
    \begin{aligned}
    | \langle \de_t z_i, \eta \rangle| = \bigg| \int_0^T \int_\Omega \nabla \eta \cdot [z_i (\nabla z_i + \nabla L_i(t,x,z_i) + \delta \bar{F}_i) + \varepsilon \nabla z_i] \d x \d t \bigg|, 
    \end{aligned}
\end{equation}
from which we obtain, using the triangle inequality, 
\begin{equation*}
    | \langle \de_t z_i, \eta \rangle| \leq \int_0^T \int_\Omega |\nabla \eta| z_i|\nabla z_i + \nabla L_i(t,x,z_i) + \delta \bar{F}_i| \d x \d t + \varepsilon \int_0^T\int_\Omega |\nabla \eta| |\nabla z_i| \d x \d t.  
\end{equation*}
Then, using H\"older's inequality, we get 
\begin{equation}\label{eq:weakform2 0}
\begin{aligned}
| \langle \de_t z_i, \eta \rangle|
\leq 
\nm{z_i}_{L^q(Q_T)} & \nm{\nabla z_i + \nabla L_i(\cdot,\cdot,z_i) + \delta \bar{F}_i}_{L^2(Q_T)}
\nm{\nabla \eta}_{L^r(Q_T)} + \varepsilon \nm{\nabla \eta}_{L^2(Q_T)} \nm{\nabla z_i}_{L^2(Q_T)}, 
\end{aligned}
\end{equation}
where $r$ satisfies
$
\frac{1}{q} + \frac{1}{r} = \frac{1}{2}, 
$
and $q$ is as given in \eqref{eq:apply dibenedetto 0}. Hence $r=2(d+1)>2$ and, since $Q_T = (0,T) \times \Omega$ is a bounded domain, an application of the H\"{o}lder inequality shows that 
\begin{equation*}
    \begin{aligned}
    \nm{\nabla \eta}_{L^2(Q_T)} &= \bigg( \int_0^T \int_\Omega |\nabla \eta|^2 \d x \d t  \bigg)^{\frac{1}{2}} \\ 
    &\leq  \bigg(\int_0^T \int_\Omega |\nabla \eta|^{2(d+1)} \d x \d t \bigg)^\frac{1}{2(d+1)}  \bigg( \int_0^T \int_\Omega 1 \d x \d t \bigg)^{\frac{d}{2(d+1)}} = (|\Omega|T)^{\frac{d}{r}} \nm{\nabla \eta}_{L^r(Q_T)}. 
    \end{aligned}
\end{equation*}
Combining the above with \eqref{eq:apply dibenedetto 0} and \eqref{eq:weakform2 0} shows that, with $C_\gamma=C_\gamma(\Omega,T,d)$ a positive constant given  by
\begin{equation*}
C_\gamma =  \gamma\left(1 + \frac{T}{|\Omega|^{\frac{d+2}{d} } }\right)  + (|\Omega|T)^{\frac{d}{r}}, 
\end{equation*}
which is independent of $\varepsilon \in (0,1)$ and $z_0 = (z_{i,0})_{i=1}^M$, there holds 
\begin{equation*}
\begin{aligned}
| \langle \de_t z_i, \eta \rangle| \leq C_\gamma 
\nm{\nabla \eta}_{L^r(Q_T)}
\bigg(  &
\nm{z_i}_{L^2(0,T;H^1(\Omega))}
\\&+  
\nm{z_i}_{V^{1,2}} \big(\nm{z_i}_{L^2(0,T;H^1(\Omega))} + \Vert \nabla L_i(\cdot,\cdot,z_i) \Vert_{L^2(Q_T)} 
+ \delta \nm{\bar{F}_i}_{L^2(Q_T)} \big)
\bigg) , 
\end{aligned} 
\end{equation*}
for any $\eta \in C^\infty(\bar{Q}_T)$, where we also used Minkowski's inequality. Using the Cauchy--Young inequality on both of the terms inside the large brackets gives
\begin{equation}\label{eq:bound dt u mid way 0}
    \begin{aligned}
    | \langle \de_t z_i, \eta \rangle| \leq 2C_\gamma \bigg( 1 + \nm{z_i}_{V^{1,2}}^2 &+\nm{z_i}_{L^2(0,T;H^1(\Omega))}^2 + \Vert \nabla L_i(\cdot,\cdot,z_i) \Vert_{L^2(Q_T)}^2 + \delta^2 \nm{\bar{F}_i}_{L^2(Q_T)}^2 \bigg) \nm{\nabla \eta}_{L^r(Q_T)}, 
    \end{aligned} 
\end{equation}
for any $\eta \in C^\infty(\bar{Q}_T)$. Observe then that 
\begin{equation*}
    \Vert \nabla L_i(\cdot,\cdot,z_i) \Vert_{L^2(Q_T)} \leq \Vert \nabla V_i \Vert_{L^2(Q_T)} + \bigg(\int_0^T \int_\Omega |\nabla W_i*z_i(t,\cdot)(x)|^2 \d x \d t\bigg)^{\frac{1}{2}}, 
\end{equation*}
where we estimate the second term on the right-hand side in the same way as in \eqref{eq:NO YOUNG}. It follows that 
\begin{equation}\label{eq:bounding Li in minkowski for time deriv bound}
\Vert \nabla L_i(\cdot,\cdot,z_i) \Vert_{L^2(Q_T)} \leq (|\Omega|T)^{\frac{1}{2}}C_L \left( 1 + \int_\Omega z_{i,0} \d x \right). 
\end{equation}
Note also that the mass conservation of Lemma \ref{lem:mass cons regularised} yields 
\begin{equation*}
    \nm{z_i}_{V^{1,2}} \leq \nm{z_{i,0}}_{L^1(\Omega)} + \nm{z_i}_{L^2(0,T;H^1(\Omega))}. 
\end{equation*}
By combining the above with \eqref{eq:bound dt u mid way 0} along with \eqref{eq:bounding Li in minkowski for time deriv bound}, and with the estimate \eqref{eq:L2H1 bound 0} of Lemma \ref{lem:energy estimates 0}, we obtain 
\begin{equation*}
    \begin{aligned}
        | \langle \de_t z_i, \eta \rangle| \leq C_{X'} \bigg( 1 + \nm{z_{i,0}}_{L^1(\Omega)}^2 + \int_\Omega z_{i,0} \log z_{i,0} \d x + \delta^2 \nm{\bar{F}_i}_{L^2(Q_T)}^2 \bigg) \nm{\nabla \eta}_{L^r(Q_T)}, 
    \end{aligned}
\end{equation*}
for any $\eta \in C^\infty(\bar{Q}_T)$, with 
$
C_{X'} = 2 C_\gamma \big( 2+3C_\Omega + 2|\Omega|TC_L^2 \big),
$
i.e., as given in \eqref{eq:CXprime expression}, which is independent of $\varepsilon$ and $z_0=(z_{i,0})_{i=1}^M$. Using the density of the smooth functions in the space $L^r(0,T;W^{1,r}(\Omega))$, we take the supremum over all test functions $\eta\in L^r(0,T;W^{1,r}(\Omega))$ in the previous estimate, and deduce the uniform estimate \eqref{eq:unif est on dt u 0}.
\end{proof}

We also make note of the following quantitative estimate on the second derivatives of the regularised solutions. The proof is contained in Appendix \ref{appendix:H2bounds}. 

\begin{lem}[Second derivative estimate]\label{lem:quantitativeH2}
For the regularised frozen system \eqref{eq:pdefrozen}, there holds, for $i \in \{1,\dots,M\}$, the estimate 
\begin{equation}\label{eq:H2bound quantitative}
		\begin{aligned}
		& \Vert \Delta z_i \Vert_{L^2(Q_T)} \leq C(\varepsilon,\delta,T,\Omega,C_L,\Vert z_{i,0} \Vert_{C^1(\bar{\Omega})},\Vert \bar{F}_i\Vert_{L^\infty(0,T;W^{1,\infty}(\Omega))},\Vert \partial_t \bar{F}_i\Vert_{L^1(0,T;L^2(\Omega))}), 
		\end{aligned}
\end{equation}
where the right-hand side is a positive quantity depending only on the parameters in its parentheses. 
\end{lem}

\subsection{Uniqueness of weak solutions to the regularised frozen problem}\label{subsec:unique}

The following result provides a correspondence between equivalent weak formulations of the regularised frozen problem \eqref{eq:pdefrozen}. Note that the proof also holds for the weak formulation of the original problem, as mentioned in Remark \ref{rem:clarify duality prod under weak sol def}. 

\begin{lem}[Equivalence of weak formulations of \eqref{eq:pdefrozen}]\label{lem:reg sol mirrors weak sol def} 
Let $z=(z_i)_{i=1}^M \in \Xi$ and denote the flux by $\ff_i := z_i(\nabla z_i + \nabla L_i(t,x,z_i) + \delta F_i(t,x,\bar{z},\nabla\bar{z})) + \varepsilon \nabla z_i$. 
The following formulations are equivalent:
\begin{enumerate}
\item for each $i \in \{1,\dots,M\}$, for every $\eta \in C^1(\bar{\Omega})$,
for a.e.~$t \in [0,T]$, 
\begin{equation}\label{eq:gregory weak sol regularised}
    \begin{aligned}
        \langle \partial_t z_i(t,\cdot), \eta \rangle_\Omega +\int_{\Omega} \ff_i \cdot\nabla\eta \d x = 0, 
    \end{aligned}
\end{equation}
where $\langle \cdot , \cdot \rangle_\Omega$ is the duality product of $W^{1,r}(\Omega)$; 
\item
for each $i \in \{1,\dots,M\}$, 
for every $\phi \in C^1(\bar{Q}_T)$, 
\begin{equation}\label{eq:weak sol regularised alternative form}
    \langle \partial_t z_i , \phi \rangle_{X'\times X} + \int_{Q_T} \ff_i\cdot\nabla\phi \d x \d t = 0; 
\end{equation}
\item for each $i \in \{1,\dots,M\}$, 
for every $\phi \in C^1(\bar{Q}_T)$ with $\phi(0,\cdot)=\phi(T,\cdot)=0$ on $\bar{\Omega}$, 
\begin{equation}\label{eq:weak sol regularised original}
\int_{Q_T} \big( -z_i \partial_t \phi + \ff_i\cdot\nabla\phi \big) \d x \d t = 0. 
\end{equation}
\end{enumerate}
\end{lem}
\begin{proof}
Before proving the lemma, we recall some useful facts. Firstly, as already noted in Remark \ref{rem:embed into cts negative sobolev space}, the definition of the space $\Xi$ implies $z_i \in W^{1,r'}(0,T;(W^{1,r}(\Omega))')$. This implies, by \cite[Theorem 2 of Section 5.9.2]{evans1998partial}, that $z_i \in C([0,T];(W^{1,r}(\Omega))')$ and 
\begin{equation*}
    z_i(t_2,\cdot) - z_i(t_1,\cdot) = \int_{t_1}^{t_2} \partial_t z_i(t,\cdot) \d t \quad \text{in } (W^{1,r}(\Omega))' \quad \forall 0 \leq t_1 < t_2 \leq T, 
\end{equation*}
and hence, given any $\eta \in C^1(\bar{\Omega})$, there holds 
\begin{equation}\label{eq:duality with dt t1 t2}
    \int_\Omega z_i(t_2,x)\eta(x) \d x - \int_\Omega z_i(t_1,x)\eta(x) \d x = \int_{t_1}^{t_2} \langle \partial_t z_i(t,\cdot) , \eta \rangle \d t \quad \forall 0 \leq t_1 < t_2 \leq T, 
\end{equation}
where the order of Bochner integration and the duality product were interchanged, which is justified by \cite[Appendix E.5, Theorem 8]{evans1998partial} and the summability
\begin{equation}\label{eq:summability}
    \int_0^T \Vert \partial_t z_i (t,\cdot) \Vert_{(W^{1,r}(\Omega))'} \d t \leq T^{\frac{1}{r}} \Vert \partial_t z_i \Vert_{X'}. 
\end{equation}
Recall also the definition of weak time derivative in terms of Bochner integration (\textit{cf.}~\cite[Section 5.9.2]{evans1998partial}). That is, for $z_i \in L^{r'}(0,T;(W^{1,r}(\Omega))')$, the element $\partial_t z_i \in L^{r'}(0,T;(W^{1,r'}(\Omega)')$ is such that 
\begin{equation}\label{eq:weak deriv recall}
    \int_0^T \theta'(t) z_i(t,\cdot) \d t = - \int_0^T \theta(t) \partial_t z_i(t,\cdot) \d t \qquad \forall \theta \in C^1_c((0,T)), 
\end{equation}
where the equality holds in the sense of $(W^{1,r}(\Omega))'$. We are now in a position to prove the lemma. 

\textbf{Step I [2.$\iff$1.]: } 
Begin with the forward implication. Fix any $\eta \in C^1(\bar{\Omega})$ and subsequently choose $\phi(t,x) := \theta(t)\eta(x)$ where $\theta \in C^1([0,T])$ is arbitrarily chosen. Then, \eqref{eq:weak sol regularised alternative form} and the identification \eqref{eq:what the dual means} made in Remark \ref{rem:clarify duality prod under weak sol def} imply $\int_0^t \big( \langle \partial_t z_i(\tau,\cdot), \eta \rangle_\Omega + \int_\Omega \ff_i \cdot \nabla \eta \d x \big) \theta(\tau) \d \tau = 0$. We deduce \eqref{eq:gregory weak sol regularised}, using the arbitrariness of $\theta$ and the Fundamental Lemma of Calculus of Variations. 

The converse follows from the density in $C^1(\bar{Q}_T)$ of the subalgebra $\{ \theta(t)\eta(x) : \theta \in C^1([0,T]), \eta \in C^1(\bar{\Omega})\}$, due to Nachbin's version of the Weierstra{\ss} Approximation Theorem \cite{nachbin}, using the continuity of $\partial_t z_i$ as an element of $X'$ for the duality term (indeed, convergence in $C^1(\bar{Q}_T)$ implies convergence in $X$) and the Dominated Convergence Theorem for the integral term.

\textbf{Step II [2.$\iff$3.]: } Testing \eqref{eq:weak deriv recall} with $\eta \in C^1(\bar{\Omega})$, using the summability \eqref{eq:summability} to exchange the order of Bochner integration and the duality product in $W^{1,r}(\Omega)$, as well as the Tonelli--Fubini Theorem, we obtain 
\begin{equation}\label{eq:about to conclude duality product time}
    \int_{Q_T} z_i(t,x)\theta'(t)\eta(x) \d x \d t = - \int_0^T  \langle \partial_t z_i(t,\cdot) , \theta(t)\eta(\cdot) \rangle_\Omega  \d t \qquad \forall \theta \in C^1_c((0,T)), \forall \eta \in C^1(\bar{\Omega}). 
\end{equation}
Using the identification \eqref{eq:what the dual means} from Remark \ref{rem:clarify duality prod under weak sol def}, we recognise the right-hand side of the above as $\langle \partial_t z_i , \theta \eta \rangle_{X'\times X}$. Thus, from the above, we obtain the equivalence of \eqref{eq:weak sol regularised alternative form} and \eqref{eq:weak sol regularised original} in the particular instance of test functions of the form $\phi(t,x) = \theta(t)\eta(x)$ with $\theta \in C^1_c((0,T))$ and $\eta \in C^1(\bar{\Omega})$. Note that this subset of functions is dense in $\{ \phi \in C^1(\bar{Q}_T) : \phi(0,\cdot)=\phi(T,\cdot)=0 \text{ on } \bar{\Omega} \}$ endowed with the subspace topology of $C^1(\bar{Q}_T)$. The result then follows from the continuity of $\partial_t z_i$ as an element of $X'$ for the duality term, and the Dominated Convergence Theorem for all remaining terms. 
\end{proof}

The next result is the focal point of this subsection, and is concerned with the uniqueness of solutions to the regularised frozen problem in the larger space $\Xi$. 

\begin{lem}[Uniqueness for regularised frozen problem in $\Xi$]\label{lem:uniqueness in weak space}
There exists a unique $z=(z_i)_{i=1}^M\in\Xi$ such that,  for every $i \in \{1,\dots,M\}$, $\int_\Omega z_i(t,x)\d x = \int_\Omega z_{i,0}(x)\d x$, $z_i(0,\cdot) = z_{i,0}$ in $(W^{1,r}(\Omega))'$  and, for every $\eta \in C^1(\bar{Q}_T)$,
\begin{equation}\label{eq:weak form}
\begin{aligned}
&\langle \de_t z_i , \eta \rangle_{X'\times X} + 
\int_{Q_T} \left[ z_i(
\nabla z_i + \nabla L_i(t,x,z_i) + \delta F_i(t,x,\bar{z},\nabla \bar{z})
)+ \eps \nabla z_i \right]
\cdot \nabla \eta \d x \d t = 0.
\end{aligned}
\end{equation}
\end{lem}

\begin{proof}
Let $z$ and $z^*$ satisfy the hypotheses and let $w:=z_i - z_i^*$. 
Problem \eqref{eq:weak form} is nonlinear, however we can formally derive a dual equation through integration by parts and properties of the convolution. Specifically, supposing that the convolution kernel $H$ is radially symmetric, and that the functions below are sufficiently integrable, we note the following property:
$$
\int_\Omega f(x)(H*g)(x) \d x
=
\int_{\Omega\times\Omega} f(x)H(x-y)g(y) \d y \d x 
=
\int_\Omega (H*f)(y)g(y) \d y, 
$$
which follows directly from the Tonelli--Fubini Theorem. 
In the sequel we use the notation $\bar{F}_i$ introduced in \eqref{eq:F bar notation}, and we drop the arguments of $\bar{F}_i$, $V_i$ and $W_i$. Define $\phi\in L^2(0,T;H^2(\Omega)) \cap H^1(0,T;L^2(\Omega))$ as a strong solution of the following linear and strongly parabolic \emph{dual equation}:
\begin{equation}
\begin{aligned}
\de_t\phi
-(\eps+a_\kappa)\Delta \phi 
+ 
( \nabla V_i + \delta \bar{F}_i
+\nabla W_i*z_i
)\cdot \nabla \phi 
+
\sum\nolimits_{k=1}^d
(\de_{x_k} W_i)
*(z_i^*\de_{x_k} \phi)
&= -\xi && \text{ in } Q_T,
\\
\nabla\phi\cdot\nu &= 0
&& \text{ on } \Sigma_T,
\\
\phi(0,\cdot)&= 0 && \text{ on } \Omega,
\end{aligned}
\end{equation}
where $\xi\in C_c^\infty(Q_T)$ is arbitrary and $\{a_\kappa\}_{\kappa \in \mathbb{N}}$ is a monotone sequence of bounded functions that approximate $\frac{z_i+z_i^*}{2}$; in particular, we choose \begin{equation*}
a_\kappa := \min
\left\{ \left(\frac{z_i+z_i^*}{2}\right), \kappa
\right\}.
\end{equation*}
Note that $a_\kappa \geq 0$ belongs to $L^\infty(Q_T)$ for every $\kappa \in \mathbb{N}$. Consider now the equation for the time-shifted function $\psi(t,x) = \phi(T-t,x)$:
\begin{equation}\label{eq: dual}
\begin{aligned}
\de_t\psi
+(\eps+a_\kappa)\Delta \psi
- 
( \nabla V_i + \delta \bar{F}_i
+\nabla W_i*z_i
)\cdot \nabla \psi 
-
\sum\nolimits_{k=1}^d
(\de_{x_k} W_i)
*(z_i^*\de_{x_k} \psi) &= \xi && \text{ in } Q_T,
\\
\nabla\psi\cdot\nu &= 0
&& \text{ on } \Sigma_T,
\\
\psi(T,\cdot)&= 0 && \text{ on } \Omega,
\end{aligned}
\end{equation}
and test against a bounded $C^1$ function $\theta: [0,T]\to [1,\infty)$ such that $\de_t \theta \geq 1$. We have, for every $t \in [0,T]$, 
\begin{equation*}
\begin{aligned}
\int_t^T \int_\Omega 
\de_t \psi \, \theta \Delta\psi
\d x\d \tau  &= 
-\int_t^T \int_\Omega 
\frac{\theta}{2}
\de_t \left(|\nabla\psi|^2 \right)
\d x\d \tau
\\&= 
-\int_\Omega
\frac{1}{2}\left(
\theta(T)|\nabla\psi(T,x)|^2
-\theta(t)|\nabla\psi(t,x)|^2
\right) \d x
+\int_t^T \int_\Omega 
|\nabla\psi|^2
\de_t \left(\frac{\theta}{2} \right)
\d x\d \tau
\\&\geq \int_\Omega \frac{1}{2}|\nabla \psi(t,x)|^2 \d x +
\int_t^T \int_\Omega
\frac{1}{2}
|\nabla\psi|^2
\d x\d \tau,
\end{aligned}
\end{equation*}
where we used $\phi(0,\cdot)=\psi(T,\cdot)=0$ and the lower bound on $\theta$. On the other hand, we also have
\begin{equation}\label{eq:on the other hand before firstly we estimate}
\begin{aligned}
\int_t^T \int_\Omega
\de_t \psi \, \theta \Delta\psi
\d x\d \tau  &= 
\int_t^T \int_\Omega \theta
\bigg[
- (\eps+a_\kappa)|\Delta\psi|^2
+ 
( \nabla V_i + \delta \bar{F}_i
+\nabla W_i*z_i
)\cdot \nabla \psi
\Delta\psi
\\ & \qquad\qquad +
\sum\nolimits_{k=1}^d
(\de_{x_k} W_i)
*(z_i^*\de_{x_k} \psi)
\Delta\psi
+\xi \Delta\psi
\bigg]
\d x\d \tau.
\end{aligned}
\end{equation}
Firstly, we estimate the drift terms: letting $\G := \nabla V_i + \delta \bar{F}_i
+\nabla W_i*z_i$,
\begin{equation*}
\begin{aligned}
\int_t^T \int_\Omega \theta
 \G\cdot \nabla \psi
\,\Delta\psi
\d x \d \tau
&=
\int_t^T \int_\Omega
\frac{\theta}{2}
|\nabla \psi|^2
\dv\G
\d x \d \tau
-
\int_t^T \int_\Omega
\theta
\sum\nolimits_{k,l} \de_{x_k} \psi \, \de_{x_l}\G_k \,
\de_{x_l} \psi
\d x \d \tau
\\&\leq
\left(\frac{1}{2}+d^2\right)\Vert \theta \Vert_{L^\infty([0,T])}\Vert \nabla \G \Vert_{L^\infty(Q_T)}
\int_{Q_T} 
|\nabla \psi|^2
\d x \d t,
\end{aligned}
\end{equation*}
where we used that $\mathcal{G}$ and its space derivative are bounded, since the nonlocal drift is bounded by $\Vert \nabla W_i * z_i \Vert_{L^\infty(Q_T)} \leq \Vert W_i \Vert_{C^1(\mathbb{R}^d)} \Vert z_{i,0} \Vert_{L^1(\Omega)}$. Secondly, we estimate the nonlocal term:
\begin{equation*}
\begin{aligned}
\bigg|\int_t^T \int_\Omega \theta
\sum\nolimits_{k=1}^d
(\de_{x_k} W_i)
*(z_i^* & \de_{x_k} \psi)
\Delta\psi
\d x \d \tau \bigg|
\\ 
&= \bigg|
\int_t^T \int_\Omega \theta
\sum\nolimits_{k=1}^d
(\nabla\de_{x_k} W_i)
*(z_i^*\de_{x_k} \psi)
\cdot\nabla \psi
\d x \d \tau \bigg|
\\&\leq
d\Vert \theta \Vert_{L^\infty([0,T])}\Vert  W_i \Vert_{C^2(\mathbb{R}^d)}
\int_t^T \int_\Omega
|\nabla \psi(\tau,x)| \bigg( \int_\Omega |z_i^*(\tau,y) \nabla \psi(\tau,y)| \d y \bigg)
\d x \d \tau
\\&\leq
d|\Omega|^{\frac{1}{2}}\Vert \theta \Vert_{L^\infty([0,T])}\Vert  W_i \Vert_{C^2(\mathbb{R}^d)} \int_t^T \Vert z^*_i(\tau,\cdot) \Vert_{L^2(\Omega)} \Vert \nabla \psi(\tau,\cdot) \Vert^2_{L^2(\Omega)} \d \tau, 
\end{aligned}
\end{equation*}
where we integrated by parts to obtain the first equality and used the H\"{o}lder inequality to obtain both the second and the final lines. By integrating the term in \eqref{eq:on the other hand before firstly we estimate} involving $\xi \Delta \psi$ by parts, and then using the Cauchy--Schwarz integral inequality followed by the Young inequality, we have obtained 
\begin{equation}\label{eq:pre gronwall estimates for uniqueness}
\begin{aligned}
    \frac{1}{2} \Vert & \nabla \psi (t,\cdot) \Vert_{L^2(\Omega)}^2 + \int_t^T \int_\Omega \big(
\frac{1}{2} 
|\nabla\psi|^2
+ (\eps+a_\kappa)|\Delta\psi|^2
\big)
\d x\d \tau
\\ 
\leq &
\frac{1}{2}\Vert \nabla \xi \Vert_{L^2(Q_T)}^2 + C(d,|\Omega|,\Vert\theta\Vert_{L^\infty},\Vert \nabla\G \Vert_{L^\infty},\Vert W_i \Vert_{C^2})
\int_t^T (1+ \Vert z^*_i(\tau,\cdot) \Vert_{L^2(\Omega)})
\frac{1}{2}\Vert \nabla \psi(\tau,\cdot) \Vert^2_{L^2(\Omega)}
\d x \d \tau. 
\end{aligned}
\end{equation}
From the above we deduce, by first estimating the integrand in the right-hand side by its supremum over $[\tau,T]$ and then bounding the left-hand side similarly, that for every $t \in [0,T]$ there holds 
\begin{equation*}
\begin{aligned}
    \frac{1}{2} \sup_{\tau \in [t,T]}\Vert \nabla \psi (\tau,\cdot) \Vert_{L^2(\Omega)}^2 \leq &
\frac{1}{2}\Vert \nabla \xi \Vert_{L^2(Q_T)}^2 + C
\int_t^T (1+ \Vert z^*_i(\tau,\cdot) \Vert_{L^2(\Omega)})
\frac{1}{2}\sup_{y \in [\tau,T]}\Vert \nabla \psi(y,\cdot) \Vert^2_{L^2(\Omega)}
\d x \d \tau. 
\end{aligned}
\end{equation*}
where we used the shorthand $C$ to denote the quantity with the same name appearing on the right-hand side of \eqref{eq:pre gronwall estimates for uniqueness}. Applying the Gr\"{o}nwall Lemma (starting from the initial point $t=T$, where $\psi(T,\cdot)=\phi(0,\cdot)=0$) to this latter inequality yields 
\begin{equation*}
\begin{aligned}
    \sup_{t \in [0,T]}\Vert \nabla \psi (t,\cdot) \Vert_{L^2(\Omega)}^2 \leq 
\Vert \nabla \xi \Vert_{L^2(Q_T)}^2 \exp\left(C\int_0^T (1+\Vert z^*_i(\tau,\cdot) \Vert_{L^2(\Omega)}) \d \tau \right), 
\end{aligned}
\end{equation*}
Using the H\"{o}lder inequality, the integral inside the exponential is bounded by $T + T^{\frac{1}{2}}\Vert z^*_i \Vert_{L^2(Q_T)}^2$. Hence, returning to \eqref{eq:pre gronwall estimates for uniqueness} and using the previous estimate, we obtain
\begin{equation}\label{eq:L2 of Laplace bound}
\int_{Q_T}  (\eps+a_\kappa)|\Delta\psi|^2
\d x\d t
\leq
 C_\xi, 
\end{equation}
where $C_\xi$ depends on
$d$, $|\Omega|$, $T$, $\Vert \nabla \xi \Vert_{L^2(Q_T)}$,
$\Vert\theta\Vert_{L^\infty}$,
$\Vert \nabla\G \Vert_{L^\infty}$,
$\Vert W_i \Vert_{C^2}$,
$\Vert z^*_i \Vert_{L^2(Q_T)} $,
and not on $\kappa$. 

The regularity of $\phi$ implies that $\psi\in X$. Recalling that $w(0,\cdot) = \psi(T,\cdot) = 0$, by proceeding as in the proof of Lemma \ref{lem:reg sol mirrors weak sol def}, we obtain
\begin{equation*}
\begin{aligned}
\langle \partial_t w , \psi \rangle_{X'\times X}
=
-\int_{Q_T} w\de_t\psi \d x \d t &=-\int_{Q_T} 
\left[ \eps \nabla w +
z_i
\nabla z_i 
-
z_i^*
\nabla z_i^* 
\right]
\cdot \nabla \psi \d x \d t
\\&\quad
- \int_{Q_T} \left[
w( \nabla V_i + \delta \bar{F}_i
)\cdot \nabla \psi 
+
\left(
z_i
\nabla W_i*z_i
-
z_i^*
\nabla  W_i*z_i^*
\right)\cdot\nabla \psi
\right]\d x \d t.
\end{aligned}
\end{equation*}
Testing against the solution of the dual problem \eqref{eq: dual}, we have
\begin{equation*}
\begin{aligned}
\int_{Q_T} 
w \bigg[
\de_t \psi 
+ \left(\eps+\frac{z_i+z_i^*}{2}\right)\Delta\psi
&-
( \nabla V_i + \delta \bar{F}_i
+\nabla W_i*z_i
)\cdot \nabla \psi
\\  & +
\sum\nolimits_{k=1}^d
(\de_{x_k} W_i)
*(z_i^*\de_{x_k} \psi)
\bigg]
\d x\d t
=
\int_{Q_T} 
w\,\xi  
\d x\d t.
\end{aligned}
\end{equation*}
It follows that
\begin{equation*}
\begin{aligned}
\int_{Q_T} 
w\,\xi  
\d x\d t
&=
\int_{Q_T} 
w\left(\frac{z_i+z_i^*}{2}-a_\kappa\right)\Delta\psi
\d x\d t
\\&\leq
\left(\int_{Q_T} 
(\eps+a_\kappa) |\Delta\psi|^2
\d x\d t
\right)^\frac{1}{2}
\left(\int_{Q_T} 
\frac{w^2}{\eps + a_\kappa}
\left(\frac{z_i+z_i^*}{2}-a_\kappa
\right)^2 
\d x\d t
\right)^\frac{1}{2}.
\end{aligned}
\end{equation*}
Thus, using the non-negativity of $a_\kappa$ and the estimate \eqref{eq:L2 of Laplace bound}, we have obtained
\begin{equation*}
\int_{Q_T} 
w\,\xi  
\d x\d t
\leq
\eps^{-1}C_\xi^{\frac{1}{2}}
\left(\int_{Q_T} 
w^2
\left(\frac{z_i+z_i^*}{2}-a_\kappa
\right)^2 
\d x\d t
\right)^\frac{1}{2}, 
\end{equation*}
and we recall that $C_\xi$ is independent of $\kappa$. Even if the integral on the right-hand side above may be infinite, since $a_\kappa \nearrow \frac{z_i+z_i^*}{2}$ a.e.~monotonically as $\kappa\to\infty$, a direct application of the Monotone Convergence Theorem leads to
$$\int_{Q_T} 
w\,\xi  
\d x\d t
\leq 0$$
for any $\xi\in C_c^\infty(Q_T)$.
In particular, we can choose $\xi$ to be strictly positive or strictly negative, so we must have $w=z_i-z^*_i=0$ a.e.~in $Q_T$.
\end{proof}

\section{Fixed point with diffusivity}\label{sec:fixed_point}

In this section, we use a fixed point argument in order to go from the regularised frozen system \eqref{eq:pdefrozen} to the \emph{regularised coupled system}
\begin{equation}\label{pdecoupled 0}
\left\lbrace\begin{aligned}
& \partial_t z_i = \dv [z_i
(\nabla z_i + \nabla L_i(t,x,z_i) + \delta F_i(t,x,z,\nabla z)) + \eps \nabla z_i] \qquad &&\text{in } Q_T, \\ 
& 0 = \nu \cdot [z_i
(\nabla z_i + \nabla L_i(t,x,z_i) + \delta F_i(t,x,z,\nabla z)) + \eps \nabla z_i] \qquad && \text{on } \Sigma_T, \\ 
& z_i(0,\cdot) = z_{i,0} \qquad &&\text{on } \Omega. 
\end{aligned}\right.
\end{equation}

We begin by recalling the Leray--Schauder--Schaefer Fixed Point Theorem and its simple corollary. 

\begin{theor*}[Leray--Schauder--Schaefer]
Let $S$ be a compact map from a Banach space $B$ into itself.
Suppose that the set $\{\xi \in B : \xi = \lambda S(\xi) \text{ for some } \lambda \in [0,1] \}$ is bounded. Then $S$ has a fixed point.
\end{theor*}

\begin{cor}\label{corfp}
Let $S$ be a compact map from a Banach space $B$ into itself.
Suppose that there exist two constants $a\in [0,1)$ and $b>0$ such that $\nm{S(\xi)}_B \leq a\nm{\xi}_B + b$ for all $\xi\in B$. Then $S$ has a fixed point.
\end{cor}

We emphasise that, throughout this entire section, $\varepsilon>0$ and the initial data $z_{0}=(z_{i,0})_{i=1}^M \in (C^\infty_c(\Omega))^M$ prescribed in Section \ref{sec:regularised} (\textit{cf.}~Remark \ref{rem:fixing initial data for pdefrozen}) are fixed.

\subsection{Weak compactness of the solution map}\label{section:weak compactness}

Recall the regularised frozen system \eqref{eq:pdefrozen} of Section \ref{sec:regularised}, i.e., 
\begin{equation*}
\left\lbrace\begin{aligned}
& \partial_t z_i = \dv [z_i
(\nabla z_i + \nabla L_i(t,x,z_i) + \delta F_i(t,x,\bar{z},\nabla\bar{z})) + \eps \nabla z_i] \qquad &&\text{in } Q_T, \\ 
& 0 = \nu \cdot [z_i
(\nabla z_i + \nabla L_i(t,x,z_i) + \delta F_i(t,x,\bar{z},\nabla\bar{z})) + \eps \nabla z_i] \qquad && \text{on } \Sigma_T, \\ 
& z_i(0,\cdot) = z_{i,0} \qquad &&\text{on } \Omega. 
\end{aligned}\right.
\end{equation*}
From Section \ref{sec:regularised} we know that the above admits a unique classical solution for each smooth vector function $\bar{z}$. We also recall the space $\Xi$ introduced in Definition \ref{defi:functional space}, and note that, since it is a closed subspace of a Banach space, it is itself Banach. 

Consider the solution operator $S_\eps$ of Definition \ref{defi:weak sol regularised}: 
\begin{equation}\label{eq:solmap}
\begin{aligned}
S_\eps : (C^\infty(\bar{Q}_T))^M &\to (C^{2,1}(\bar{Q}_T))^M \subset \Xi
\\
\bar{z} &\mapsto z,
\end{aligned}
\end{equation}
where $z$ solves \eqref{eq:pdefrozen} and $(C^\infty(\bar{Q}_T))^M$ is equipped with the subspace topology of $\Xi$. Let us introduce, for every $\mu>0$, the following linear smoothing operator $R_\mu$ (which is defined by extension and mollification in Appendix \ref{appendix:smoothing op}): 
\begin{equation}\label{eq:Rnu introduce}
\begin{aligned}
R_\mu : \Xi &\to (C^\infty(\bar{Q}_T))^M
\\
w &\mapsto R_\mu(w),
\end{aligned}
\end{equation}
with the property that $R_\mu(w)$ converges to $w$ strongly in $(L^2(0,T;H^1(\Omega)))^M$  as $\mu \to 0$, along with 
\begin{equation}\label{eq:L2H1 bound on Rnu}
    \nm{R_\mu(w)}_{(L^2(0,T;H^1(\Omega)))^M} \leq C_{reg} \nm{w}_{(L^2(0,T;H^1(\Omega)))^M}, 
\end{equation}
for some fixed constant $C_{reg}$ independent of $\mu$, and 
\begin{equation}\label{eq:not unif L2H2 bound on Rnu}
    \nm{R_\mu(w)}_{(L^\infty(0,T;W^{2,\infty}(\Omega)))^M} + \nm{\partial_t R_\mu(w)}_{(L^\infty(0,T;W^{1,\infty}(\Omega)))^M} \leq C_\mu \nm{w}_{(L^2(0,T;H^1(\Omega)))^M}, 
\end{equation}
for some positive constant $C_\mu$ depending on $\mu$ (\textit{cf.}~Lemma \ref{lem:smoothing operator} and Remark \ref{rem:time deriv L2H1 depends mu reg}). This constant explodes in the limit as $\mu \to 0$.

We will also repeatedly make use of the following two technical lemmas in later arguments. Their proofs are contained in Appendix \ref{sec:two technical seq lemmas}. 

\begin{lem}\label{lem:preserve sign and mass}
Let $\{\zeta^n\}_{n\in\mathbb{N}}$ be a sequence in $L^2(Q_T)$ such that:
\begin{enumerate}
    \item $\{\zeta^n\}_{n\in\mathbb{N}}$ converges weakly to $\zeta$ in $L^2(Q_T)$;
    \item $\zeta^n$ is non-negative a.e.~in $Q_T$ for every $n\in\mathbb{N}$; 
    \item $\int_{Q_T} \zeta^n(t,x) \d x = \Lambda$ a.e.~$t \in (0,T)$ for some non-negative constant $\Lambda$. 
\end{enumerate}
Then $\zeta \in L^\infty(0,T;L^1(\Omega))$, 
$\displaystyle
    \zeta \geq 0
$ a.e. in $Q_T$, and 
$\displaystyle
    \int_\Omega \zeta(t,x) \d x = \Lambda. 
$
\end{lem}

\begin{lem}\label{lem:preserve initial data}
Let $\{\zeta^n\}_{n\in\mathbb{N}}$ be a sequence in $L^2(Q_T)$ such that: 
\begin{enumerate}
\item $\{\zeta^n\}_{n\in\mathbb{N}}$ converges weakly to $\zeta$ in $L^2(Q_T)$; 
\item For every $n\in\mathbb{N}$, $\zeta^n(0,\cdot)=\zeta_0$ in $(W^{1,r}(\Omega))'$ for some fixed $\zeta_0 \in L^p(\Omega)$; 
\item $\{\partial_t\zeta^n\}_{n\in\mathbb{N}}$ converges weakly-* to $\partial_t \zeta$ in $X'$. 
\end{enumerate}
Then, $\zeta \in C([0,T];(W^{1,r}(\Omega))')$, and there exists a positive constant $\Lambda$, independent of $n$ such that, given any $\phi \in W^{1,r}(\Omega)$, there holds 
\begin{equation}\label{eq:to show cty for zeta}
    \bigg|\int_\Omega \zeta(t,x) \phi(x) \d x - \int_\Omega \zeta(s,x) \phi(x) \d x \bigg| \leq (t-s)^{\frac{1}{r}} \Lambda \Vert \phi \Vert_{W^{1,r}(\Omega)} \quad \forall 0 < s \leq t \leq T.
\end{equation}
Moreover, 
\begin{equation}\label{eq:to show cty for zeta with initial}
    \Vert \zeta(t,\cdot) - \zeta_0 \Vert_{(W^{1,r}(\Omega))'} \leq t^{\frac{1}{r}} \Lambda  \quad \forall 0 \leq t \leq T, 
\end{equation}
so that $\zeta(0,\cdot) = \zeta_0$ in $(W^{1,r}(\Omega))'$. 
\end{lem}

\begin{lem}[Weak compactness]\label{lem:weak compactness}
The map
$
S_\eps \circ R_\mu
$
from $\Xi$ into itself is weakly sequentially compact with respect to the subspace topology of $(L^2(0,T;H^1(\Omega))^M$.
\end{lem}
\begin{proof}
For $i \in \{1,\dots,M\}$, let $\{\bar{z}^n_i\}_{n\in\NN}$ be a uniformly bounded sequence in $\Xi$, i.e., there exists $C_i>0$, independent of $n,\mu,\eps$, such that 
\begin{equation}\label{eq:unibound}
\nm{\bar{z}^n_i }_{L^2(0,T;H^1(\Omega))} \leq C_i, \qquad \nm{\partial_t \bar{z}^n_i }_{X'} \leq C_i \qquad \forall n\in\NN.
\end{equation}
Let us introduce 
\begin{equation*}
\hat{z}^n := R_\mu(\bar{z}^n),
\qquad
z^n := S_\eps(\hat{z}^n),
\end{equation*}
and notice that by the estimate \eqref{eq:not unif L2H2 bound on Rnu}, up to the multiplicative constant $C_{reg}$, we have that $\bar{z}=(\bar{z}_i)_{i=1}^M$  and $\hat{z}=(\hat{z}_i)_{i=1}^M$ satisfy the same uniform bound \eqref{eq:unibound}. We also define 
\begin{equation*}
\hat{F}^n_i := F_i(t,x,\hat{z}^n,\nabla\hat{z}^n).
\end{equation*}
and notice that, due to \eqref{eq:fbound0} and \eqref{eq:unibound}, we have 
\begin{equation*}
    \Vert \hat{F}^n_i \Vert_{L^2(Q_T)} \leq C_F(1+C_i) \qquad \forall n\in\NN, 
\end{equation*}
where we omitted the $C_{reg}$ factor, for clarity of presentation. 
By Lemmas \ref{lem:existence of viscous approximates 0}, \ref{lem:energy estimates 0}, and \ref{lem:time compactness bound 0}, there exist a positive constants $C_i'$ independent of $n, \mu,\eps$ such that 
\begin{equation*}
\nm{z^n_i }_{L^2(0,T;H^1(\Omega))} \leq C_i', \qquad 
\nm{\partial_t z^n_i }_{X'} \leq C_i',
\qquad
\forall n\in\NN.
\end{equation*}

It follows that all of $\{\bar{z}_i^n\}_{n\in\NN}$, $\{\hat{z}_i^n\}_{n\in\NN}$, and $\{z_i^n\}_{n\in\NN}$ are bounded sequences in $\Xi$. An application of the theorems of Banach--Alaoglu and Aubin--Lions \cite[Theorem II.5.16]{boyer2012mathematical} implies that there exists a common subsequence (still indexed by $n$) such that
\begin{equation*}
\begin{aligned}
\bar{z}_i^n \rightharpoonup \bar{z}_i, \quad \hat{z}^n_i \rightharpoonup \hat{z}_i, \quad 
& {z}_i^n \rightharpoonup z_i \qquad &&\text{weakly in } L^2(0,T;H^1(\Omega)), \\ 
&\hat{F}^n_i \rightharpoonup \hat{F}^*_i \qquad &&\text{weakly in } L^2(Q_T), \\
\bar{z}^n_i \to \bar{z}_i, \quad &z_i^n \to z_i \qquad &&\text{strongly in } L^2(Q_T), \\ 
&\partial_t z^n_i \overset{*}{\rightharpoonup} v_i \qquad &&\text{weakly-* in } X', 
\end{aligned}
\end{equation*}
where $\bar{z}_i, \hat{z}_i , z_i \in L^2(0,T;H^1(\Omega))$, $\hat{F}^*_i \in L^2(Q_T)$, and $v_i \in X'$. From the linearity and the continuity property of the regularisation operator (\textit{cf.}~Corollary \ref{cor:smoothing operator L2 cylinder conv}), 
\begin{equation*}
    \Vert \hat{z}^n_i - R_\mu\bar{z}_i \Vert_{L^2(Q_T)} = \Vert R_\mu \bar{z}^n_i - R_\mu\bar{z}_i \Vert_{L^2(Q_T)} \leq C \Vert \bar{z}^n_i - \bar{z}_i \Vert_{L^2(Q_T)} \to 0 \qquad \text{as } n \to \infty, 
\end{equation*}
for some positive constant $C$, which, incidentally, does not depend on $\mu$. Hence it follows that $\hat{z}_i = R_\mu\bar{z}_i$ as elements of $L^2(Q_T)$, and that, additionally, $\hat{z}^n_i \to \hat{z}_i$ strongly in $L^2(Q_T)$. Moreover, given any $\Theta \in (C^1_c(Q_T))^d$, there holds, from the definition of weak derivative, 
\begin{equation*}
	\begin{aligned}
		\bigg|\int_{Q_T} \Theta \cdot \nabla (R_\mu \bar{z}_i - \hat{z}_i ) \d x \d t\bigg| = \bigg|\int_{Q_T} \dv \Theta (R_\mu \bar{z}_i - \hat{z}_i ) \d x \d t\bigg| \leq \Vert \dv \Theta \Vert_{L^2(Q_T)} \Vert R_\mu \bar{z}_i - \hat{z}_i \Vert_{L^2(Q_T)}, 
		\end{aligned}
\end{equation*}
from which we deduce that $\Vert \nabla (R_\mu \bar{z}_i  - \hat{z}_i ) \Vert_{L^2(Q_T)} = 0$, due to the density of $C^1_c(Q_T)$ in $L^2(Q_T)$. Thus, $\hat{z}_i = R_\mu \bar{z}_i$ as elements of $L^2(0,T;H^1(\Omega))$.

Additionally, we note that, in view of Lemma \ref{lem:preserve sign and mass}, we have $z_i \geq 0$ a.e.~in $Q_T$ and 
\begin{equation*}
	\int_\Omega |z_i(t,x)| \d x = \int_\Omega z_i(t,x) \d x = \int_\Omega z_{i,0} \d x \qquad \text{a.e.~}t \in (0,T), 
\end{equation*}
which fulfils the requirement for $z_i \in L^\infty(0,T;L^1(\Omega))$.

Furthermore, using the structure of $\hat{F}_i^n$ provided by \eqref{eq:F structure 0}, 
we know that 
\begin{equation*}
\hat{F}^*_i = F_i(t,x,\hat{z},\nabla\hat{z}) = G^0_i(t,x,\hat{z}) + \sum_{j=1}^M G^1_{ij}(t,x,\hat{z}) \nabla \hat{z}_j. 
\end{equation*}
Indeed, letting $\Theta \in (C^1_c(Q_T))^d$ be arbitrary, the structure \eqref{eq:F structure 0} implies that 
\begin{equation*}
	\int_{Q_T} \Theta \cdot \hat{F}^n_i \d x \d t = \int_{Q_T} \Theta \cdot G^0_i(t,x,\hat{z}^n) \d x \d t + \sum_{j=1}^M \int_{Q_T}  \Theta \cdot \nabla G^1_{ij}(t,x,\hat{z}^n) \nabla \hat{z}_j \d x \d t. 
\end{equation*}
Then, the fundamental theorem of calculus and the strong convergence in $L^2(Q_T)$ yield 
\begin{equation}\label{eq:in weak cpctness strong L2 conv for G0}
	\Vert G^0_i(\cdot,\cdot,\hat{z}^n) - G^0_i(\cdot,\cdot,\hat{z}) \Vert_{L^2(Q_T)} \leq \Vert \nabla_z G^0_i \Vert_{L^\infty(Q_T \times \mathbb{R}^M)} \max_{i \in \{1,\dots,M\}}\Vert \hat{z}^n_i - \hat{z}_i \Vert_{L^2(Q_T)} \to 0, 
\end{equation}
and similarly 
\begin{equation}\label{eq:in weak cpctness strong L2 conv for G1}
	\Vert G^1_{ij}(\cdot,\cdot,\hat{z}^n) - G^1_{ij}(\cdot,\cdot,\hat{z}) \Vert_{L^2(Q_T)} \leq \Vert \nabla_z G^1_{ij} \Vert_{L^\infty(Q_T \times \mathbb{R}^M)} \max_{i \in \{1,\dots,M\}}\Vert \hat{z}^n_i - \hat{z}_i \Vert_{L^2(Q_T)} \to 0, 
\end{equation}
so that $G^0_i(\cdot,\cdot,\hat{z}^n) \to G^0_i(\cdot,\cdot,\hat{z})$ and $G^1_{ij}(\cdot,\cdot,\hat{z}^n) \to G^1_{ij}(\cdot,\cdot,\hat{z})$ strongly in $L^2(Q_T)$. The weak convergence also implies $\nabla\hat{z}_j^n \rightharpoonup \nabla\hat{z}_j$ weakly in $L^2(Q_T)$, so that, using the fact that the product of a strongly converging sequence with a weakly converging sequence converges itself in the weak sense, 
\begin{equation*}
	\lim_{n\to\infty} \int_{Q_T} \Theta \cdot \hat{F}^n_i \d x \d t = \int_{Q_T} \Theta \cdot \hat{F}_i \d x \d t. 
\end{equation*}

Similarly, with the term $L_i$, we have 
\begin{equation*}
	\begin{aligned}
		\Vert L_i(t,x,z^n_i) - L_i(t,x,z_i) \Vert_{L^2(Q_T)}^2 &= \int_0^T \int_\Omega \big| (W_i*(z^n_i(t,\cdot)-z_i(t,\cdot)))(x) \big|^2 \d x \d t \\ 
		&= \int_0^T \int_\Omega \bigg| \int_\Omega W_i(x-y)(z^n_i(t,y)-z_i(t,y)) \d y \bigg|^2 \d x \d t \\ 
		&\leq |\Omega|C_L^2 \int_0^T \int_\Omega \bigg( \int_\Omega |z^n_i(t,y)-z_i(t,y)|^2 \d y \bigg) \d x \d t \\ 
		&= |\Omega|^2 C_L^2 \Vert z^n_i - z_i \Vert_{L^2(Q_T)}^2 \to 0, 
		\end{aligned}
\end{equation*}
where we applied Jensen's inequality and used the boundedness \eqref{eq:drift V bound} to obtain the third line, and the Tonelli--Fubini theorem to obtain the final line. An identical strategy yields 
\begin{equation*}
	\begin{aligned}
		\Vert \nabla L_i(t,x,z^n_i) - \nabla L_i(t,x,z_i) \Vert_{L^2(Q_T)}^2 &= \int_0^T \int_\Omega \big| (\nabla W_i*(z^n_i(t,\cdot)-z_i(t,\cdot)))(x) \big|^2 \d x \d t \\ 
		&\leq |\Omega|^2 C_L^2 \Vert z^n_i - z_i \Vert_{L^2(Q_T)}^2 \to 0, 
		\end{aligned}
\end{equation*}
so that we obtain the convergence $L_i(\cdot,\cdot,z^n_i) \to L(\cdot,\cdot,z_i)$ strongly in $L^2(0,T;H^1(\Omega))$. 

Furthermore, given any $\theta \in C^1_c(Q_T)$, an integration by parts with respect to the time variable yields 
\begin{equation*}
\int_{Q_T} \theta \de_t z^n_i \d x \d t = -\int_{Q_T}   z^n_i \de_t\theta \d x \d t. 
\end{equation*}
By taking the weak limits on both sides, we get $\int_{Q_T} \theta v_i \d x \d t = -\int_{Q_T}   z_i \de_t\theta \d x \d t$ and deduce $v_i = \partial_t z_i$. Hence, since any weakly-* convergent sequence is bounded, we have that $\Vert \partial_t z_i \Vert_{X'} < +\infty$, and the final requirement for $z_i$ belonging to $\Xi$ is fulfilled. We therefore have the following convergences for the subsequence: 
\begin{equation}\label{eq:convergences table weak compactness}
\begin{aligned}
\bar{z}^n_i \rightharpoonup \bar{z}_i, \quad \hat{z}^n_i \rightharpoonup R_\mu \bar{z}_i , \quad 
&z^n_i \rightharpoonup z_i \qquad &&\text{weakly in } L^2(0,T;H^1(\Omega)),
\\ 
\bar{z}^n_i \to \bar{z}_i, \quad \hat{z}^n_i \to R_\mu\bar{z}_i, \quad &z^n_i \to z_i \qquad &&\text{strongly in } L^2(Q_T), \\ 
&L_i(\cdot,\cdot,z^n_i) \to L_i(\cdot,\cdot,z_i) \qquad && \text{strongly in } L^2(0,T;H^1(\Omega)), \\ 
&\hat{F}^n_i \rightharpoonup F_i(\cdot,\cdot,R_\mu\bar{z},\nabla R_\mu\bar{z}) \qquad &&\text{weakly in } L^2(Q_T), 
\\ 
&\partial_t z^n_i \overset{*}{\rightharpoonup} \partial_t z_i \qquad &&\text{weakly-* in } X', 
\end{aligned}
\end{equation}
with $z_i \in \Xi$, and the lemma is proved, since we have shown the weak convergence in $\Xi$ of a subsequence of $\{S_\varepsilon \circ R_\mu(\bar{z}^n)\}_{n\in\mathbb{N}}$ towards $z = (z_i)_{i=1}^M$ for any bounded sequence $\{\bar{z}_n\}_{n\in\mathbb{N}}$ in $\Xi$. 
\end{proof}

\begin{rem}
We emphasise that, in order to obtain \eqref{eq:convergences table weak compactness}, the Aubin--Lions Lemma (\textit{cf.}~\cite[Theorem II.5.16]{boyer2012mathematical}) was used to ensure the strong convergences $\bar{z}^n_i \to \bar{z}_i$ and $z^n_i \to z_i$ in $L^2(Q_T)$. The application of \cite[Theorem II.5.16]{boyer2012mathematical} is justified due to the uniform bound in $L^2(0,T;H^1(\Omega))$ for the sequence of functions $\{\bar{z}^n_i,z^n_i\}_{n\in\mathbb{N}}$ and the uniform bound in $X' = L^{r'}(0,T;(W^{1,r}(\Omega))')$ for the corresponding sequence of derivatives $\{\partial_t\bar{z}^n_i,\partial_t z^n_i \}_{n\in\mathbb{N}}$. The strong convergence $\hat{z}^n_i \to R_\mu \bar{z}_i$ in $L^2(Q_T)$ was not deduced directly from the Aubin--Lions Lemma (though, alternatively, this can also be done) and was later obtained from properties of the regularisation operator $R_\mu$ and the strong convergence $\bar{z}^n_i \to \bar{z}_i$ in $L^2(Q_T)$. 
\end{rem}

\begin{rem}\label{rem:about Rnu in convergence fixed point}
As a consequence of  \eqref{eq:convergences table weak compactness}, for any $\phi \in C^1(\bar{Q}_T)$, 
\begin{equation*}
\begin{aligned}
	&\int_{Q_T} 
z^n_i \nabla z^n_i  \cdot\nabla\phi  \d x \d t
\to
\int_{Q_T} 
z_i \nabla z_i  \cdot\nabla\phi  \d x \d t, \\ 
&\int_{Q_T} z^n_i \nabla L_i(t,x,z^n_i) \cdot \nabla \phi \d x \d t \to \int_{Q_T} z_i \nabla L_i(t,x,z_i) \cdot \nabla \phi \d x \d t, \\ 
&\int_{Q_T} 
z^n_i F_i(t,x,R_\mu \bar{z}^n,\nabla R_\mu \bar{z}^n)\cdot\nabla\phi  \d x \d t
\to
\int_{Q_T} 
z_i F_i(t,x,R_\mu\bar{z},\nabla R_\mu\bar{z})\cdot\nabla\phi \d x \d t. 
\end{aligned}
\end{equation*}
It is then clear that the limit function $z_i$ satisfies the following weak formulation: 
\begin{equation*}
\begin{aligned}
&\langle \de_t z_i , \phi \rangle_{X'\times X} + 
\int_{Q_T} \nabla \phi \cdot \big( z_i(\nabla z_i + \nabla L_i(t,x,z_i) + \delta F_i(t,x,R_\mu\bar{z},\nabla R_\mu\bar{z}) + \eps \nabla z_i \big) \d x \d t = 0, 
\end{aligned}
\end{equation*}
for every $\phi \in C^1(\bar{Q}_T)$. Note the similarity between the no-flux weak formulation above and the formulation \eqref{eq:weak sol regularised alternative form} in Lemma \ref{lem:reg sol mirrors weak sol def}. 
\end{rem}

\subsection{Strong compactness of the solution map}\label{subsec:strong}

In this subsection, we improve the compactness result in Lemma \ref{lem:weak compactness} and show that the solution map is actually strongly compact from $\Xi$ to itself. To begin with, we recall the following result concerning semicontinuity properties of the Fisher information. 

\begin{lem}[Properties of Fisher information, \cite{ambrosio2014calculus} Lemma 4.10]\label{lem:lowersemi}
Let $K$ be a closed subset of $\RR^d$. The Fisher information, i.e., the functional 
\begin{equation*}
\ff[w] := \int_{\{x\in K | w(x) > 0\}} \frac{|\nabla w|^2}{w} \d x, 
\end{equation*}
is convex and sequentially lower semicontinuous with respect to the weak topology of $L^1(K)$.
\end{lem}

As a consequence, we have the following lemma, the proof of which is contained in Appendix \ref{sec:two technical seq lemmas}. 
\begin{lem}\label{lem:lowersemi the one we use}
Let $\{\zeta^n\}_{n\in\mathbb{N}}$ be a sequence of non-negative functions in $L^1(Q_T)$, and suppose that $\zeta^n \to \zeta$ strongly in $L^1(Q_T)$. Then there exists a subsequence of $\{\zeta^n\}_{n\in\mathbb{N}}$, still indexed by $n$, for which there holds 
\begin{equation*}
    \int_{Q_T} \frac{|\nabla \zeta|^2}{\zeta} \d x \d t \leq \liminf_{n\to\infty} \int_{Q_T} \frac{|\nabla \zeta^n|^2}{\zeta^n} \d x \d t. 
\end{equation*}
\end{lem}

\begin{lem}[Strong compactness]\label{lem:strong compactness}
The map
$
S_\eps \circ R_\mu
$
from $\Xi$ into itself is strongly sequentially compact.
\end{lem}
\begin{proof}
We emphasise that $\varepsilon$ and $\mu$ are fixed throughout this proof. For $i \in \{1,\dots,M\}$, let $\{\bar{z}^n_i\}_{n\in\NN}$ be a uniformly bounded sequence in $\Xi$.
Arguing as in Lemma \ref{lem:weak compactness}, we consider $\hat{z}^n=R_\mu \bar{z}^n$ and $z^n = S_\eps \hat{z}^n = S_\eps \circ R_\mu (\bar{z}^n)$. Recalling \eqref{eq:convergences table weak compactness}, we have that, for a suitable subsequence (still indexed by $n$), 
\begin{equation*}
\begin{aligned}
\bar{z}^n_i \rightharpoonup \bar{z}_i, \quad \hat{z}^n_i \rightharpoonup R_\mu \bar{z}_i , \quad 
&z^n_i \rightharpoonup z_i \qquad &&\text{weakly in } L^2(0,T;H^1(\Omega)),
\\ 
\bar{z}^n_i \to \bar{z}_i, \quad \hat{z}^n_i \to R_\mu\bar{z}_i, \quad &z^n_i \to z_i \qquad &&\text{strongly in } L^2(Q_T), \\ 
&L_i(\cdot,\cdot,z^n_i) \to L_i(\cdot,\cdot,z_i) \qquad && \text{strongly in } L^2(0,T;H^1(\Omega)), \\ 
&\hat{F}^n_i \rightharpoonup F_i(\cdot,\cdot,R_\mu\bar{z},\nabla R_\mu\bar{z}) \qquad &&\text{weakly in } L^2(Q_T), 
\\ 
&\partial_t z^n_i \overset{*}{\rightharpoonup} \partial_t z_i \qquad &&\text{weakly-* in } X', 
\end{aligned}
\end{equation*}
for some non-negative $z_i \in \Xi$, and we define $\hat{z} := R_\mu \bar{z}$. By integrating the equality \eqref{eq:precede a priori} from Lemma \ref{lem:energy estimates 0} over the time interval $[t_0,t] \subsetneq [0,T]$, we obtain 
\begin{equation*}
\begin{aligned}
\int_\Omega \big[ z_i^n(t) & \log z_i^n(t) - z_i^n(t_0) \log z_i^n(t_0) \big] \d x 
\\ &=
-  \int_{t_0}^t \int_\Omega \left[
\delta\nabla z_i^n \cdot F_i(t,x,\hat{z}^n,\nabla\hat{z}^n) + \nabla z_i^n \cdot \nabla L_i(t,x,z^n_i)
+ |\nabla z_i^n|^2 
+ \eps\frac{|\nabla z_i^n|^2}{z_i^n} \right] \d x \d\tau, 
\end{aligned}
\end{equation*}
for a.e. $t>t_0>0$. By Remark \ref{rem:about Rnu in convergence fixed point}, the uniqueness in the space $\Xi$ due to Lemma \ref{lem:uniqueness in weak space}, and the added regularity inherited from Lemma \ref{lem:existence of viscous approximates 0}, we deduce that $z_i \in C^{2,1}(\bar{Q}_T)$ satisfies \eqref{eq:pdefrozen} classically (with $\hat{z}=R_\mu \bar{z}$ featuring in the arguments of $F_i$). Then, similarly to what we had in Section \ref{subsec:uniform}, we also obtain 
\begin{equation*}
\begin{aligned}
\int_\Omega \big[ z_i(t) & \log z_i(t) - z_i(t_0) \log z_i(t_0) \big] \d x 
\\ &=
-  \int_{t_0}^t \int_\Omega \left[
\nabla z_i \cdot \delta F_i(t,x,\hat{z},\nabla\hat{z}) + \nabla z_i \cdot \nabla L_i(t,x,z_i)
+ |\nabla z_i|^2 
+ \eps\frac{|\nabla z_i|^2}{z_i} \right] \d x \d\tau, 
\end{aligned}
\end{equation*}
for a.e. $t>t_0>0$.
Taking the difference of the two relations above we obtain
\begin{equation}\label{eq:difference}
\begin{aligned}
\int_{t_0}^t \int_\Omega 
&\left[\left(1+\frac{\eps}{z_i^n}\right)
|\nabla z_i^n|^2 
-
\left(1+\frac{\eps}{z_i}\right)
|\nabla z_i|^2 
\right]
\d x \d \tau
\\ &=
\int_\Omega  \big( z_i(t_0) \log z_i(t_0) - z_i^n(t) \log z_i^n(t) -
z_i(t_0) \log z_i(t_0) + z_i(t) \log z_i(t) \big) \d x
\\ 
&\quad -\int_0^t \int_\Omega \big( \nabla z^n_i \cdot \nabla L_i(t,x,z^n_i) - \nabla z_i \cdot \nabla L_i(t,x,z_i) \big) \d x \d \tau \\ 
&\quad 
- \delta \int_0^t \int_\Omega 
\big( \nabla z_i^n \cdot F_i(t,x,\hat{z}^n,\nabla\hat{z}^n) 
- \nabla z_i \cdot F_i(t,x,\hat{z},\nabla\hat{z}) 
\big) \d x \d \tau.
\end{aligned}
\end{equation}
One can show, by following Step I of the proof of Lemma \ref{lem:lowersemi the one we use} (\textit{cf.}~Appendix \ref{sec:two technical seq lemmas}), that the strong convergence $z^n_i \to z_i$ in $L^2(Q_T)$ implies that, for a subsequence, for a.e.~$t\in(0,T)$, we have $\Vert z^n_i(t,\cdot) - z_i(t,\cdot) \Vert_{L^2(\Omega)} \to 0$. Then, defining $f:x \mapsto x (\log x)\mathds{1}_{[0,\infty)}(x)$, we note that $|f(x)| \leq C(1+x^2)$ globally for some universal constant $C$. Using the Generalised Dominated Convergence Theorem, we deduce that $f$ maps $L^2(\Omega)$ continuously into $L^1(\Omega)$, whence the entire first term on the right-hand side of \eqref{eq:difference} vanishes for this subsequence. The second term on the right-hand side of \eqref{eq:difference} also vanishes, due to the strong convergence in $L^2(0,T;H^1(\Omega))$ of the terms involving $\nabla L_i$, \textit{cf.}~\eqref{eq:convergences table weak compactness}. On the other hand, the final term on the right-hand side of \eqref{eq:difference} requires additional work.

Recall that, due to the structure provided by \eqref{eq:F structure 0}, 
\begin{equation*}
	F_i(t,x,\hat{z},\nabla\hat{z}) = G^0_i(t,x,\hat{z}) + \sum_{j=1}^M G^1_{ij}(t,x,\hat{z}) \nabla \hat{z}_j, 
\end{equation*}
and we consider, in particular, the term 
\begin{equation}\label{eq:term}
\int_{t_0}^t \int_\Omega 
\big(
G^1_{ij}(\tau,x,\hat{z}^n) \nabla \hat{z}_j^n
\cdot\nabla z_i^n  
- G^1_{ij}(\tau,x,\hat{z}) \nabla \hat{z}_j
\cdot\nabla z_i 
\big) \d x \d \tau. 
\end{equation}
For what follows we define the $(d+1)$-dimensional vector fields 
\begin{equation*}
    \hat{v}^n_j(\tau,x) := (0,G^1_{ij}(\tau,x,\hat{z}^n(\tau,x)) \nabla \hat{z}_j^n(\tau,x)), \qquad v^n_i(\tau,x) := (0,\nabla z^n_i(\tau,x)). 
\end{equation*}
Note that the strong convergence in $L^2(Q_T)$ of the terms involving $G^1_{ij}$, \textit{cf.}~\eqref{eq:in weak cpctness strong L2 conv for G1}, and the weak convergences $\nabla \hat{z}^n_j \rightharpoonup \nabla \hat{z}_j , \nabla z^n_i \rightharpoonup \nabla z_i$ in $L^2(Q_T)$ implies that both sequences $\{\hat{v}^n_j\}_{n\in\mathbb{N}} , \{v^n_i\}_{n\in\mathbb{N}}$ are weakly convergent in $(L^2(Q_T))^{d+1}$. In the next paragraph, we pass to the limit in \eqref{eq:term} using the div-curl Lemma.

By Lemma \ref{lem:quantitativeH2}, $\Vert \Delta z^n \Vert_{L^2(Q_T)}$ is bounded by $\Vert R_\mu \bar{z}^n \Vert_{L^\infty(0,T;W^{2,\infty}(\Omega))} + \Vert \partial_t R_\mu \bar{z}^n \Vert_{L^\infty(0,T;W^{1,\infty}(\Omega))}$, which, from the estimate \eqref{eq:not unif L2H2 bound on Rnu}, is bounded by $\Vert \bar{z}^n \Vert_{L^2(0,T;H^1(\Omega))}$, and this is bounded independently of $n$. Therefore, for every $i\in\{1,\dots,M\}$, 
\begin{equation*}
    \dv_{t,x} v^n_i = \Delta z^n_i \qquad \text{is bounded in } L^2(Q_T) \text{ independently of } n, 
\end{equation*}
and thus, by the Rellich Theorem, the sequence $\{\dv_{t,x} v^n_i \}_{n\in\mathbb{N}}$ is confined to a compact subset of $H^{-1}(Q_T)$. Additionally, an explicit computation using the chain rule shows that, for every $j\in\{1,\dots,M\}$, 
\begin{equation*}
    \begin{aligned}
        |\crl_{t,x} \hat{v}^n_j| \leq \max_{i,j\in\{1,\dots,M\}}\Vert G^1_{ij} \Vert_{C^1(\bar{Q}_T\times\mathbb{R}^M)} \big( |\partial_t \nabla \hat{z}^n| &+ |\nabla^2 \hat{z}^n| + |\nabla \hat{z}^n| + |\partial_t \hat{z}^n||\nabla \hat{z}^n| + |\nabla \hat{z}^n|^2 \big). 
    \end{aligned}
\end{equation*}
Recall that $\mu>0$ is fixed. Hence, estimate \eqref{eq:not unif L2H2 bound on Rnu} shows that $\Vert \hat{z}^n_i \Vert_{L^\infty(0,T;W^{2,\infty}(\Omega))} = \Vert R_\mu \bar{z}^n_i \Vert_{L^\infty(0,T;W^{2,\infty}(\Omega))}$ is bounded independently of $n$ and likewise for  $\Vert \partial_t \hat{z}^n_i \Vert_{L^\infty(0,T;W^{1,\infty}(\Omega))} = \Vert \partial_t R_\mu \bar{z}^n_i \Vert_{L^\infty(0,T;W^{1,\infty}(\Omega))}$. In turn, $\{\crl_{t,x} \hat{v}^n_j\}_{n\in\mathbb{N}}$ is bounded in $L^\infty(Q_T)$ independently of $n$, and therefore also uniformly bounded in $L^2(Q_T)$, whence the Rellich Theorem implies that this sequence is confined to a compact subset of $H^{-1}(Q_T)$. A direct application of the $\dv$-$\crl$ Lemma (\textit{cf.}~\cite[Theorem 1]{murat}) to the product $\{\hat{v}^n_j \cdot v^n_i\}_{n\in\mathbb{N}}$ yields that 
\begin{equation*}
    G^1_{ij}(t,x,\hat{z}^n) \nabla \hat{z}_j^n
\cdot\nabla z_i^n  
\to G^1_{ij}(t,x,\hat{z}) \nabla \hat{z}_j
\cdot\nabla z_i \qquad \text{in } \mathcal{D}'(Q_T). 
\end{equation*}
We conclude that the term in \eqref{eq:term} vanishes in the as $n\to\infty$. Note that the term 
\begin{equation*}
	\int_{t_0}^t \int_\Omega \big( G^0_i(t,x,\hat{z}^n) \cdot \nabla z^n_i - G^0_i(t,x,\hat{z})\cdot \nabla z_i \big) \d x \d \tau 
\end{equation*}
also vanishes in the limit as $n\to\infty$, since $\nabla z^n_i \rightharpoonup \nabla z_i$ weakly in $L^2(Q_T)$ and $G^0_i(\cdot,\cdot,\hat{z}^n) \to G^0_i(\cdot,\cdot,\hat{z})$ strongly in $L^2(Q_T)$, as per the estimate \eqref{eq:in weak cpctness strong L2 conv for G0}. 

Returning to \eqref{eq:difference} and using the fact that $0<t_0<t<T$ were arbitrary, it follows (by possibly taking a further subsequence to let $t_0 \to 0$ and $t \to T$) that 
\begin{equation}\label{eq:before using fisher}
\int_0^T \int_\Omega 
\left[\left(1+\frac{\eps}{z_i^n}\right)
|\nabla z_i^n|^2 
-
\left(1+\frac{\eps}{z_i}\right)
|\nabla z_i|^2 
\right]
\d x \d \tau
\to 0 \qquad \text{ as } n\to\infty.
\end{equation}
Observe that the strong convergence $z^n_i \to z_i$ in $L^2(Q_T)$ and the lower semicontinuity result Lemma \ref{lem:lowersemi the one we use} imply, after passing to a further subsequence if necessary, 
\begin{equation*}
\liminf_{n\to\infty}
\int_0^T \int_\Omega 
\left[\frac{1}{z_i^n}
|\nabla z_i^n|^2 
-
\frac{1}{z_i}
|\nabla z_i|^2 
\right]
\d x \d \tau
\geq 0,
\end{equation*}
which, combining with \eqref{eq:before using fisher}, implies 
\begin{equation*}\limsup_{n\to\infty} \big(\nm{\nabla z_i^n}_{L^2(Q_T)}^2 - \nm{\nabla z_i}_{L^2(Q_T)}^2 \big) \leq 0. 
\end{equation*}
Combining the above with $\nm{\nabla z_i}_{L^2(Q_T)} \leq \liminf_{n \to \infty}\nm{\nabla z_i^n}_{L^2(Q_T)}$, which holds true because of the weak lower semicontinuity of the norm and $\nabla z^n_i \rightharpoonup \nabla z_i$ weakly in $L^2(Q_T)$, we deduce 
\begin{equation}\label{eq:strong conv L2H1}
\nm{\nabla z_i^n}_{L^2(Q_T)} \to \nm{\nabla z_i}_{L^2(Q_T)}
 \qquad \text{as } n\to\infty. 
\end{equation}
The combination of weak convergence of $\{z_i^n\}_{n\in\NN}$ in $L^2(0,T;H^1(\Omega))$ with the convergence of the norm establishes strong convergence in $L^2(0,T;H^1(\Omega))$.

Finally, we verify the strong convergence in the dual space $X'$ for the sequence of time derivative $\{\de_t z^n_i\}_{n\in\mathbb{N}}$.
Taking the difference of the weak formulations we obtain, for any $\theta \in C^1_c(Q_T)$, 
\begin{equation*}
\begin{aligned}
\langle \de_t (z^n_i - z_i) , \theta \rangle
&=
- \int_{Q_T} \big(
\eps\nabla (z^n_i - z_i) +
z^n_i\nabla z^n_i  - z_i\nabla z_i
\big) \cdot\nabla\theta 
\d x \d t
\\
& \quad - \int_{Q_T} \big( z_i^n \nabla L_i(t,x,z^n_i) - z_i \nabla L_i(t,x,z_i) \big) \cdot \nabla \theta \d x \d t 
\\ & \quad
- \int_{Q_T} \delta\big(
 z^n_i F_i(t,x,\hat{z}^n,\nabla \hat{z}^n) - 
z_i F_i(t,x,\hat{z},\nabla \hat{z}) 
 \big)\cdot\nabla\theta  \d x \d t, 
 \end{aligned}
\end{equation*}
and, as per the convergences identified in Remark \ref{rem:about Rnu in convergence fixed point}, the right-hand side vanishes in the limit as $n\to\infty$. However, we must study this limit quantitatively. To this end, observe from this previous equation that 
\begin{equation*}
\begin{aligned}
|\langle \de_t (z^n_i - z_i) , \theta \rangle|
\leq \big( \varepsilon &\Vert \nabla z^n_i - \nabla z_i \Vert_{L^2(Q_T)} + \Vert z^n_i \nabla z^n_i - z_i \nabla z_i \Vert_{L^2(Q_T)} \\ 
&+ \Vert z_i^n \nabla L_i(\cdot,\cdot,z^n_i) - z_i \nabla L_i(\cdot,\cdot,z_i) \Vert_{L^2(Q_T)} \\ 
&+ |\delta| \Vert z_i^n F_i(\cdot,\cdot,\hat{z}^n,\nabla \hat{z}^n) - z_i F_i(\cdot,\cdot,\hat{z},\nabla \hat{z}) \Vert_{L^2(Q_T)} \big) \Vert \nabla \theta \Vert_{L^2(Q_T)}, 
 \end{aligned}
\end{equation*}
after which an application of the H\"{o}lder inequality yields 
\begin{equation*}
\begin{aligned}
|\langle \de_t (z^n_i - z_i) , \theta \rangle|
\leq \big( \varepsilon &\Vert \nabla z^n_i - \nabla z_i \Vert_{L^2(Q_T)} + \Vert z^n_i \nabla z^n_i - z_i \nabla z_i \Vert_{L^2(Q_T)} \\ 
&+ \Vert z_i^n \nabla L_i(\cdot,\cdot,z^n_i) - z_i \nabla L_i(\cdot,\cdot,z_i) \Vert_{L^2(Q_T)} \\ 
&+ |\delta| \Vert z_i^n F_i(\cdot,\cdot,\hat{z}^n,\nabla \hat{z}^n) - z_i F_i(\cdot,\cdot,\hat{z},\nabla \hat{z}) \Vert_{L^2(Q_T)} \big) (|\Omega| T)^{\frac{d}{2(d+1)}} \Vert \nabla \theta \Vert_{L^r(Q_T)}. 
 \end{aligned}
\end{equation*}
Taking the supremum over all $\theta\in C^1_c(Q_T)$ with $\Vert \theta \Vert_X \leq 1$ and using the density of $C^1_c(Q_T)$ in $X$, we get 
\begin{equation*}
\begin{aligned}
\Vert \de_t (z^n_i - z_i) \Vert_{X'}
\leq (|\Omega| T)^{\frac{d}{2(d+1)}}\big( \varepsilon &\Vert \nabla z^n_i - \nabla z_i \Vert_{L^2(Q_T)} + \Vert z^n_i \nabla z^n_i - z_i \nabla z_i \Vert_{L^2(Q_T)} \\ 
&+ \Vert z_i^n \nabla L_i(\cdot,\cdot,z^n_i) - z_i \nabla L_i(\cdot,\cdot,z_i) \Vert_{L^2(Q_T)} \\ 
&+ |\delta| \Vert z_i^n F_i(\cdot,\cdot,\hat{z}^n,\nabla \hat{z}^n) - z_i F_i(\cdot,\cdot,\hat{z},\nabla \hat{z}) \Vert_{L^2(Q_T)} \big). 
 \end{aligned}
\end{equation*}
Using the strong convergence in $L^2(Q_T)$ of $\nabla z^n_i \to \nabla z_i$, $z^n_i \to z_i$, and $\nabla L_i(\cdot,\cdot,z^n_i) \to \nabla L_i(\cdot,\cdot,z_i)$, we immediately deduce that the first three terms on the right-hand side of the previous equation vanish in the limit as $n\to\infty$, since the product of two strongly convergent sequences in $L^2(Q_T)$ is itself strongly convergent in $L^2(Q_T)$. Finally, recall that $\{\hat{z}^n\}_{n\in\mathbb{N}}$ is precisely the sequence $\{R_\mu\bar{z}^n\}_{n\in\mathbb{N}}$, and that the regularisation parameter $\mu>0$ is fixed throughout this procedure. Thus, since $\Vert \bar{z}^n_i \Vert_{L^2(0,T;H^1(\Omega))}$ is bounded independently of $n$, the estimate \eqref{eq:not unif L2H2 bound on Rnu} gives the uniform boundedness in $L^2(Q_T)$ of the sequences of higher derivatives $\{\partial_t \nabla \hat{z}^n_i\}_{n\in\mathbb{N}}$ and $\{\nabla^2 \hat{z}^n_i\}_{n\in\mathbb{N}}$. It then follows from the the Aubin--Lions Lemma (\textit{cf.}~\cite[Theorem II.5.16]{boyer2012mathematical}) that $\nabla \hat{z}^n_i \to \nabla \hat{z}_i$ strongly in $L^2(Q_T)$, since we already knew from \eqref{eq:convergences table weak compactness} that the sequence converged weakly in $L^2(0,T;H^1(\Omega))$ to this same limit. Combining with the strong convergence $\hat{z}^n_i \to \hat{z}_i$, we get 
\begin{equation*}
    \hat{z}^n_i \to \hat{z}_i \qquad \text{strongly in } L^2(0,T;H^1(\Omega)). 
\end{equation*}
It is then immediate from the structure of $F_i$ given in \eqref{eq:F structure 0} that we have the strong convergence 
\begin{equation*}
    \Vert F_i(\cdot,\cdot,\hat{z}^n,\nabla \hat{z}^n) - F_i(\cdot,\cdot,\hat{z},\nabla \hat{z}) \Vert_{L^2(Q_T)} \to 0 \qquad \text{as } n \to \infty, 
\end{equation*}
whence $\Vert \de_t (z^n_i - z_i) \Vert_{X'}$ vanishes in the limit as $n\to\infty$, again using the fact that the product of two sequences converging strongly in $L^2(Q_T)$ is itself strongly convergent in $L^2(Q_T)$. The proof is complete. 
\end{proof}

We now arrive at the existence of solutions to the regularised coupled system \eqref{pdecoupled 0}.

\begin{prop}[Existence for regularised coupled system]\label{cor:existence coupled sys}
Fix $\varepsilon>0$ and $z_{i,0} \in C^\infty_c(\Omega)$ to be non-negative functions such that $\int_\Omega z_{i,0} \d x = \int_\Omega u_{i,0} \d x$ for $i \in \{1,\dots,M\}$. There exists $z=(z_i)_{i=1}^M$, belonging to the space $\Xi$, which solves the regularised coupled system \eqref{pdecoupled 0}, i.e., 
\begin{equation*}
\left\lbrace\begin{aligned}
& \partial_t z_i = \dv [z_i
(\nabla z_i + \nabla L_i(t,x,z_i) + \delta F_i(t,x,z,\nabla z)) + \eps \nabla z_i] \qquad &&\text{in } Q_T, \\ 
& 0 = \nu \cdot [z_i
(\nabla z_i + \nabla L_i(t,x,z_i) + \delta F_i(t,x,z,\nabla z)) + \eps \nabla z_i] \qquad && \text{on } \Sigma_T, \\ 
& z_i(0,\cdot) = z_{i,0} \qquad &&\text{on } \Omega, 
\end{aligned}\right.
\end{equation*}
in the weak sense: for any test function $\phi \in C^1(\bar{Q}_T)$, for $i \in \{1,\dots,M\}$, 
\begin{equation}\label{eq:weak sol coupled sys 0}
    \begin{aligned}
        \langle \partial_t z_i , \phi \rangle_{X' \times X} +\int_{Q_T} \big[ z_i(\nabla z_i + \nabla L_i(t,x,z_i)  + \delta F_i(t,x,z,\nabla z)) + \varepsilon \nabla z_i \big]\cdot\nabla\phi \d x \d t = 0, 
    \end{aligned}
\end{equation}
with $z_i(0,\cdot) = z_{i,0}$ in $(W^{1,r}(\Omega))'$. Moreover, each $z_i$ is non-negative and conserves its initial mass, and there exists a positive constant $C=C(\Omega,T,d,\delta)$, which is independent of $\varepsilon$ and $z_0=(z_{i,0})_{i=1}^M$, such that, for $i\in\{1,\dots,M\}$, 
\begin{equation}\label{eq:estimates for coupled sys 0}
\Vert z_i \Vert_{L^2(0,T;H^1(\Omega))}^2 + \nm{\de_t z_i}_{X'}
\leq
C \bigg( 1 + \Vert z_{i,0} \Vert_{L^1(\Omega)}^2 + \int_\Omega z_{i,0}\log z_{i,0} \d  x \bigg). 
\end{equation}
\end{prop}

\begin{proof}

\textbf{Step I: } Begin by showing that for each $\mu>0$, there exists $z=(z_i)_{i=1}^M$, belonging to the space $C^{2,1}(\bar{Q}_T)$, which solves the system 
\begin{equation}\label{pdecoupled}
\left\lbrace\begin{aligned}
& \partial_t z_i = \dv [z_i
(\nabla z_i + L_i(t,x,z_i) + \delta F_i(t,x,R_\mu z,\nabla R_\mu z)) + \eps \nabla z_i] \qquad &&\text{in } Q_T, \\ 
& 0 = \nu \cdot [z_i
(\nabla z_i + L_i(t,x,z_i) + \delta F_i(t,x,R_\mu z,\nabla R_\mu z)) + \eps \nabla z_i] \qquad && \text{on } \Sigma_T, \\ 
& z_i(0,\cdot) = z_{i,0} \qquad &&\text{on } \Omega, 
\end{aligned}\right.
\end{equation}
in the classical sense, and that there exists a positive constant $C=C(\Omega,T,d,C_{reg},\delta)$, independent of $\mu,\varepsilon$, and $z_{i,0}$, such that, for $i\in\{1,\dots,M\}$, 
\begin{equation}\label{eq:estimates for coupled sys}
\Vert z_i \Vert_{L^2(0,T;H^1(\Omega))}^2 + \nm{\de_t z_i}_{X'}
\leq
C \bigg( 1 + \Vert z_{i,0} \Vert_{L^1(\Omega)}^2 + \int_\Omega z_{i,0}\log z_{i,0} \d  x \bigg). 
\end{equation}
To this end, recall that 
$S_\varepsilon \circ R_\mu$
maps the Banach space $\Xi$ into itself compactly (\textit{cf.}~Lemmas \ref{lem:energy estimates 0} and \ref{lem:strong compactness}, and \eqref{delta_bound} the smallness condition on $\delta$). Subsequently,
an application of Corollary \ref{corfp}  shows that there exists $z \in \Xi$ which is a fixed point of $S_\varepsilon \circ R_\mu$. Since $S_\varepsilon \circ R_\mu$ maps $\Xi$ into $(C^{2,1}(\bar{Q}_T))^M$, the equality in $\Xi$  $z = S_\varepsilon (R_\mu (z))$ implies that there exists a representative of $z$ belonging to $(C^{2,1}(\bar{Q}_T))^M$ and, additionally, such representative solves \eqref{pdecoupled} in the classical sense. In view of the definition of $\Xi$, we automatically obtain that, for $i\in\{1,\dots,M\}$, $z_i$ is non-negative and conserves its initial mass. 

By integrating \eqref{pdecoupled} against any test function $\phi \in C^1(\bar{Q}_T)$, we get, for $i \in \{1,\dots,M\}$, 
\begin{equation}\label{eq:weak sol coupled sys}
    \begin{aligned}
        \langle \partial_t z_i, \phi \rangle_{X' \times X} +\int_{Q_T} \big( z_i(\nabla z_i + \nabla L_i(t,x,z_i)  + \delta F_i(t,x,R_\mu z,\nabla R_\mu z)) + \varepsilon \nabla z_i \big)\cdot\nabla\phi \d x \d t = 0, 
    \end{aligned}
\end{equation}
and we also have $z_i(0,\cdot) = z_{i,0}$ in $(W^{1,r}(\Omega))'$, since we already know that this latter equality holds in the pointwise sense. 

Meanwhile, the estimate on $\Vert z_i \Vert_{L^2(0,T;H^1(\Omega))}^2$ in \eqref{eq:estimates for coupled sys} follows from Lemma \ref{lem:energy estimates 0}, using the smallness of $\delta$ along with the bound 
\begin{equation*}
    |F_i(t,x,R_\mu z(t,x),\nabla R_\mu z(t,x)| \leq C_F(T,\Omega)(1+|\nabla R_\mu z(t,x)|), 
\end{equation*}
due to \eqref{eq:fbound0}, and the estimate \eqref{eq:L2H1 bound on Rnu} (\textit{cf.}~Lemma \ref{lem:smoothing operator}). Using this latter bound and the one for $\Vert z_i \Vert_{L^2(0,T;H^1(\Omega))}^2$, we then obtain the estimate on $\nm{\de_t z_i}_{X'}$ in \eqref{eq:estimates for coupled sys} from Lemma \ref{lem:time compactness bound 0}.

\textbf{Step II: } For each $\mu>0$, define $z^\mu = (z^\mu_i)_{i=1}^M$ to be the solution of \eqref{pdecoupled} provided by Step I, and recall that \eqref{eq:weak sol coupled sys} holds. Observe from the estimates \eqref{eq:estimates for coupled sys} that $\{z^\mu\}_{\mu>0}$ is a bounded sequence in $(L^2(0,T;H^1(\Omega)))^M$, and this bound is independent of $\mu$. Hence, as in the proof of Lemma \ref{lem:weak compactness}, an application of the theorem of Banach--Alaoglu for reflexive spaces and the Aubin--Lions Lemma (\textit{cf.}~\cite[Theorem II.5.16]{boyer2012mathematical}) implies the existence of a subsequence, which we still label as $\{z^{\mu}\}_{\mu>0}$, converging weakly in $(L^2(0,T;H^1(\Omega)))^M$ and strongly in $(L^2(Q_T))^M$ to some $z = (z_i)_{i=1}^M \in (L^2(0,T;H^1(\Omega)))^M$, and such that $\partial_t z^\mu \overset{*}{\rightharpoonup} \partial_t z$ weakly-* in $(X')^M$. This latter weak-* convergence in $(X')^M$ is manifestly enough to pass to the limit in the first term of the weak formulation \eqref{eq:weak sol coupled sys}, i.e., the duality product, and for the final term, we have 
\begin{equation*}
    \begin{aligned}
        I_\mu := \bigg| \int_{Q_T} \big[ z^\mu_i(\nabla z^\mu_i &+ \nabla L_i(t,x,z_i^\mu)  + \delta F_i(t,x,R_\mu z^\mu,\nabla R_\mu z^\mu)) + \varepsilon \nabla z^\mu_i \big]\cdot\nabla\phi \d x \d t \\ 
    &- \int_{Q_T} \big[ z_i(\nabla z_i + \nabla L_i(t,x,z_i)  + \delta F_i(t,x,z,\nabla z)) + \varepsilon \nabla z_i \big]\cdot\nabla\phi \d x \d t \bigg|,  
\end{aligned}
\end{equation*}
which can be expanded as $I_\mu \leq J_\mu + K_{1,\mu} + K_{2,\mu}$, where 
\begin{equation*}
    \begin{aligned}
        J_\mu := \bigg| \int_{Q_T} z^\mu_i(\nabla z^\mu_i   + \delta F_i & (t,x,R_\mu z^\mu,\nabla R_\mu z^\mu)) \cdot \nabla \phi \d x \d t - \int_{Q_T} z_i(\nabla z_i   + \delta F_i(t,x,z,\nabla z)) \cdot \nabla \phi \d x \d t \bigg|, 
    \end{aligned}
\end{equation*}
while 
\begin{equation*}
    K_{1,\mu} := \bigg| \int_{Q_T} \varepsilon \nabla z^\mu_i \cdot\nabla\phi \d x \d t - \int_{Q_T} \varepsilon \nabla z_i \cdot\nabla\phi \d x \d t \bigg|, 
\end{equation*}
and 
\begin{equation*}
    K_{2,\mu} := \bigg| \int_{Q_T} z^\mu_i \nabla L_i(t,x,z^\mu_i) \cdot \nabla \phi \d x \d t - \int_{Q_T} z_i \nabla L_i(t,x,z_i) \cdot \nabla \phi \d x \d t \bigg|. 
\end{equation*}
Observe that 
\begin{equation*}
    K_{1,\mu} = \varepsilon \bigg| \int_{Q_T} (\nabla z^\mu_i - \nabla z_i) \cdot \nabla \phi \d x \d t \bigg| \to 0, 
\end{equation*}
in view of the weak convergence $z^\mu_i \rightharpoonup z_i$ in $L^2(0,T;H^1(\Omega))$, which naturally implies $\nabla z^\mu_i \rightharpoonup \nabla z_i$ weakly in $(L^2(Q_T))^M$. Additionally, observe that 
\begin{equation*}
    K_{2,\mu} = \bigg| \int_{Q_T} (z^\mu_i \nabla L_i(t,x,z^\mu_i) - z_i \nabla L_i(t,x,z_i)) \cdot \nabla \phi \d x \d t \bigg|, 
\end{equation*}
and 
\begin{equation*}
    \begin{aligned}
        \Vert \nabla L_i(\cdot,\cdot,z^\mu_i) - \nabla L_i(\cdot,\cdot,z_i) \Vert_{L^2(Q_T)}^2 &= \int_0^T \int_\Omega |\nabla W_i*(z^\mu_i(t,\cdot)-z_i(t,\cdot))(x)|^2 \d x \d t \\ 
        &= \int_0^T \int_\Omega \bigg|\int_\Omega \nabla W_i(x-y) (z^\mu_i(t,y)-z_i(t,y)) \d y \bigg|^2 \d x \d t \\ 
        &\leq |\Omega| C_L^2 \int_0^T \int_\Omega \bigg( \int_\Omega  |z^\mu_i(t,y)-z_i(t,y)|^2 \d y \bigg) \d x \d t \\ 
        &= |\Omega|^2 C_L^2 \Vert z^\mu_i - z_i \Vert_{L^2(Q_T)}^2 \to 0 \qquad \text{as } \mu \to 0, 
    \end{aligned}
\end{equation*}
where we used the boundedness from \eqref{eq:drift V bound} and Jensen's inequality to obtain the third line, and the Fubini--Tonelli theorem to obtain the final equality, and the strong convergence $z^\mu_i \to z_i$ in $L^2(Q_T)$ to show that the limit vanishes. Thus, 
\begin{equation*}
    \Vert \nabla L_i(\cdot,\cdot,z^\mu_i) - \nabla L_i(\cdot,\cdot,z_i) \Vert_{L^2(Q_T)} \to 0 \qquad \text{as } \mu \to 0, 
\end{equation*}
and hence, being the product of two strongly convergent sequences, we get $z^\mu_i \nabla L_i(\cdot,\cdot,z^\mu_i) \to z_i \nabla L_i(\cdot,\cdot,z_i)$ strongly in $L^2(Q_T)$, whence the Cauchy--Schwarz integral inequality yields that $K_{2,\mu} \to 0$ in the limit as $\mu \to 0$. 

We return to $J_\mu$, which we write as $J_\mu \leq J_{1,\mu} + |\delta| J_{2,\mu}$, where 
\begin{equation*}
    J_{1,\mu} := \bigg| \int_{Q_T} (z^\mu_i \nabla z^\mu_i - z_i \nabla z_i) \cdot \nabla \phi \d x \d t  \bigg|, 
\end{equation*}
and 
\begin{equation*}
    J_{2,\mu} := \bigg| \int_{Q_T} \big( z^\mu_i F_i(t,x,R_\mu z^\mu, \nabla R_\mu z^\mu) - z_i F_i(t,x,z, \nabla z) \big) \cdot \nabla \phi \d x \d t  \bigg|. 
\end{equation*}
Recall that since $z^\mu_i \to z_i$ strongly in $L^2(Q_T)$ and $\nabla z^\mu_i \rightharpoonup \nabla z_i$ weakly in $(L^2(Q_T))^M$, the product $z^\mu_i \nabla z^\mu_i$ converges weakly to $z_i \nabla z_i$ in $(L^2(Q_T))^M$, whence $J_{1,\mu} \to 0$ in the limit as $\mu$ vanishes. It remains to control $J_{2,\mu}$. Observe that an application of the triangle inequality yields 
\begin{equation*}
    \begin{aligned}
        J_{2,\mu} \leq \bigg| \int_{Q_T} & (z^\mu_i - z_i) F_i(t,x,R_\mu z^\mu, \nabla R_\mu z^\mu) \cdot \nabla \phi \d x \d t \bigg| \\ 
        &+ \bigg| \int_{Q_T} z_i \big( F_i(t,x,R_\mu z^\mu, \nabla R_\mu z^\mu) - F_i(t,x,z, \nabla z) \big) \cdot \nabla \phi \d x \d t \bigg| =: L_{1,\mu} + L_{2,\mu}. 
    \end{aligned}
\end{equation*}
Note that 
\begin{equation*}
    \begin{aligned}
        L_{1,\mu} &\leq \Vert \nabla \phi \Vert_{L^\infty(Q_T)} \max_{i \in \{1,\dots,M\}} \Vert F_i(\cdot,\cdot,R_\mu z^\mu , \nabla R_\mu z^\mu) \Vert_{L^2(Q_T)} \Vert z^\mu_i - z_i \Vert_{L^2(Q_T)} \\ 
        & \leq \Vert \nabla \phi \Vert_{L^\infty(Q_T)} C_F \big( 1 + \max_{i \in \{1,\dots,M\}}\Vert \nabla R_\mu z^\mu \Vert_{L^2(Q_T)} \big) \Vert z^\mu_i - z_i \Vert_{L^2(Q_T)}, 
    \end{aligned}
\end{equation*}
where we used the usual bound \eqref{eq:fbound0} for $F_i$, and from which it follows, using the estimate \eqref{eq:L2H1 bound on Rnu}, that 
\begin{equation*}
    \begin{aligned}
        L_{1,\mu} &\leq \Vert \nabla \phi \Vert_{L^\infty(Q_T)} C_F \big( 1 + C_{reg} \max_{i \in \{1,\dots,M\}}\Vert z^\mu_i \Vert_{L^2(0,T;H^1(\Omega))} \big) \Vert z^\mu_i - z_i \Vert_{L^2(Q_T)}. 
    \end{aligned}
\end{equation*}
Furthermore, due to the weak convergence $z^\mu_i \rightharpoonup z_i$ in $L^2(0,T;H^1(\Omega))$ for each $i \in \{1,\dots,M\}$, we know that for each $i \in \{1,\dots,M\}$ this sequence is bounded, i.e., there exists a positive constant $C_{seq}$ (independent of $\mu$) such that 
\begin{equation*}
    \max_{i \in \{1,\dots,M\}}\Vert z^\mu_i \Vert_{L^2(0,T;H^1(\Omega))} \leq C_{seq} \qquad \forall \mu>0. 
\end{equation*}
Thus, we obtain 
\begin{equation*}
    \begin{aligned}
        L_{1,\mu} &\leq \Vert \nabla \phi \Vert_{L^\infty(Q_T)} C_F ( 1 + C_{reg} C_{seq}) \Vert z^\mu_i - z_i \Vert_{L^2(Q_T)} \to 0 \qquad \text{as } \mu \to 0. 
    \end{aligned}
\end{equation*}

The term $L_{2,\mu}$ is more delicate, and to study it we expand the terms in $F_i$. Indeed, an application of the triangle inequality yields 
\begin{equation*}
    \begin{aligned}
        L_{2,\mu} \leq \bigg| & \int_{Q_T} z_i \big( G_i^0(t,x,R_\mu z^\mu) - G_i^0(t,x,z) \big) \cdot \nabla \phi \d x \d t \bigg| \\ 
        &+ \sum_{j=1}^M \bigg| \int_{Q_T} z_i \big( G^1_{ij}(t,x,R_\mu z^\mu)\nabla R_\mu z^\mu_j - G^1_{ij}(t,x,z) \nabla z_j \big) \cdot \nabla \phi \d x \d t \bigg| =: \Lambda_{1,\mu} + \Lambda_{2,\mu}. 
    \end{aligned}
\end{equation*}
Observe then that, using H\"{o}lder's inequality along with the mean value inequality and the uniform boundedness of the derivatives of $(G^0_i)_{i=1}^M$, we control the first term by 
\begin{equation}\label{eq:first estimate M1nu}
    \Lambda_{1,\mu} \leq \Vert \nabla \phi \Vert_{L^\infty(Q_T)} \Vert z_i \Vert_{L^2(Q_T)} \Vert \nabla_z G^0_i \Vert_{L^\infty(Q_T\times \mathbb{R}^M)} \max_{i \in \{1,\dots,M\}} \Vert R_\mu z^\mu_i - z_i \Vert_{L^2(Q_T)}. 
\end{equation}
Then, for each fixed $i \in \{1,\dots,M\}$, we have 
\begin{equation}\label{eq:Rnu znu converges strongly to z in L2}
    \begin{aligned}
        \Vert R_\mu z^\mu_i - z_i \Vert_{L^2(Q_T)} &\leq \Vert R_\mu z^\mu_i - R_\mu z_i \Vert_{L^2(Q_T)} + \Vert R_\mu z_i - z_i \Vert_{L^2(Q_T)} \\ 
        &\leq C_{reg}\Vert z^\mu_i - z_i \Vert_{L^2(Q_T)} + \Vert R_\mu z_i - z_i \Vert_{L^2(Q_T)} \to 0 \qquad \text{as } \mu \to 0, 
    \end{aligned}
\end{equation}
by virtue of the the strong convergence $z^\mu_i \to z_i$ in $L^2(Q_T)$ and Corollary \ref{cor:smoothing operator L2 cylinder conv}, where we also used the linearity of the operator $R_\mu$ to obtain the second inequality. Thus, by returning to \eqref{eq:first estimate M1nu}, we conclude that $\Lambda_{1,\mu}$ vanishes in the limit as $\mu \to 0$. For the term $\Lambda_{2,\mu}$, observe that, for each fixed $i,j \in \{1,\dots,M\}$, we have firstly that 
\begin{equation*}
    \Vert G^1_{ij}(\cdot,\cdot,R_\mu z^\mu) - G^1_{ij}(\cdot,\cdot,z) \Vert_{L^2(Q_T)} \leq \Vert \nabla_z G^1_{ij} \Vert_{L^\infty(Q_T \times \mathbb{R}^M)} \max_{i\in\{1,\dots,M\}} \Vert R_\mu z^\mu_i - z_i \Vert_{L^2(Q_T)}, 
\end{equation*}
where, as before, we used the mean value inequality and the uniform boundedness of the derivatives of $(G^1_{ij})_{i,j=1}^M$, so that we deduce from \eqref{eq:Rnu znu converges strongly to z in L2} that 
\begin{equation}\label{eq:strong convergence of easier mixed term in M2nu}
    G^1_{ij}(\cdot,\cdot,R_\mu z^\mu) \to G^1_{ij}(\cdot,\cdot,z) \qquad \text{strongly in } L^2(Q_T). 
\end{equation}
Meanwhile, for the first term, we have the following claim, to be proved later. 
\begin{claim}\label{claim:annoying diagonal}
\begin{equation*}
    \nabla R_\mu z^\mu_j \rightharpoonup \nabla z_j \qquad \text{weakly in } (L^2(Q_T))^M. 
\end{equation*}
\end{claim}
Using the above, we then deduce from \eqref{eq:strong convergence of easier mixed term in M2nu} that we have the weak convergence of the product, i.e., 
\begin{equation}\label{eq:final weak convergence for M2nu}
    G^1_{ij}(\cdot,\cdot,R_\mu z^\mu)\nabla R_\mu z^\mu_j \rightharpoonup G^1_{ij}(\cdot,\cdot,z)\nabla z_j \qquad \text{weakly in } (L^2(Q_T))^M. 
\end{equation}
Thus, using the H\"{o}lder inequality to verify that $z_i \nabla \phi \in (L^2(Q_T))^M$, the weak convergence of \eqref{eq:final weak convergence for M2nu} implies that $\Lambda_{2,\mu} \to 0$ in the limit as $\mu \to 0$. The limit function $z \in L^2(0,T;H^1(\Omega))$ thereby satisfies the weak formulation \eqref{eq:weak sol coupled sys 0}, and the estimates \eqref{eq:estimates for coupled sys 0} follow from the weak (and weak-* for the bound in $X'$) lower semi-continuity of the norms considered in the estimate \eqref{eq:estimates for coupled sys}. The non-negativity and mass conservation follow again directly from Lemma \ref{lem:preserve sign and mass}. The attainment of the initial data in the dual Sobolev space follows from a direct application of Lemma \ref{lem:preserve initial data}. The proof is complete. 
\end{proof}

\begin{proof}[Proof of Claim \ref{claim:annoying diagonal}]
Recall from \eqref{eq:Rnu znu converges strongly to z in L2} that we have 
\begin{equation*}
    \Vert R_\mu z^\mu_i - z_i \Vert_{L^2(Q_T)} \to 0 \qquad \text{as } \mu \to 0. 
\end{equation*}
It therefore follows that, given any test function $\psi \in (C^\infty_c(Q_T))^d$, by definition of weak derivative, there holds the integration by parts relation 
\begin{equation*}
    \int_{Q_T} \nabla(R_\mu z^\mu_i - z_i) \cdot \psi \d x \d t = - \int_{Q_T} (R_\mu z^\mu_i - z_i)  \dv \psi \d x \d t, 
\end{equation*}
where we note that there are no boundary terms due to the compact support of $\psi$. Thus, we get from the above, and the H\"{o}lder inequality, that 
\begin{equation*}
    \bigg| \int_{Q_T} \nabla(R_\mu z^\mu_i - z_i) \cdot \psi \d x \d t \bigg| \leq \Vert R_\mu z^\mu_i - z_i \Vert_{L^2(Q_T)} \Vert \nabla \psi \Vert_{L^2(Q_T)} \to 0 \qquad \text{as } \mu \to 0. 
\end{equation*}
Finally, due to the density of $C^\infty_c(Q_T)$ in $L^2(Q_T)$, it follows that given any $\psi \in (L^2(Q_T))^d$, there exists a sequence $\{\psi_m\}_{m\in\mathbb{N}}$ of elements of $C^\infty_c(Q_T)$ such that $\Vert \psi_m - \psi \Vert_{L^2(Q_T)}$ vanishes as $m \to \infty$. Then, we get 
\begin{equation}\label{eq:splitting claim}
    \begin{aligned}
        \bigg| \int_{Q_T} \nabla (R_\mu z^\mu_i - z_i) \cdot \psi \d x \d t \bigg| &\leq \bigg| \int_{Q_T} \nabla(R_\mu z^\mu_i - z_i) \cdot (\psi-\psi_m) \d x \d t \bigg| + \bigg| \int_{Q_T} \nabla(R_\mu z^\mu_i - z_i) \cdot \psi_m \d x \d t \bigg| \\ 
    &\leq \Vert \nabla R_\mu z^\mu_i - \nabla z_i \Vert_{L^2(Q_T)} \Vert \psi - \psi_m \Vert_{L^2(Q_T)} + \Vert R_\mu z^\mu_i - z_i \Vert_{L^2(Q_T)} \Vert \nabla \psi_m \Vert_{L^2(Q_T)}, 
    \end{aligned}
\end{equation}
where, for the second term in the final line, we integrated by parts and then used the H\"{o}lder inequality. Now observe that, for the first term, the triangle inequality yields 
\begin{equation*}
    \begin{aligned}
        \Vert \nabla R_\mu z^\mu_i - \nabla z_i \Vert_{L^2(Q_T)} &\leq \Vert \nabla R_\mu z^\mu_i \Vert_{L^2(Q_T)} + \Vert \nabla z_i \Vert_{L^2(Q_T)} \\ 
        &\leq \Vert R_\mu z^\mu_i \Vert_{L^2(0,T;H^1(\Omega))} + \Vert \nabla z_i \Vert_{L^2(0,T;H^1(\Omega))}. 
    \end{aligned}
\end{equation*}
However, since $z^\mu_i \in L^2(0,T;H^1(\Omega))$ for each fixed $\mu>0$, it follows from estimate \eqref{eq:L2H1 bound on Rnu} (\textit{cf.}~\eqref{eq:L2H1 unif bound Rnu smoother} of Lemma \ref{lem:smoothing operator}) that, for each fixed $\mu>0$, 
\begin{equation*}
    \Vert R_\mu z^\mu_i \Vert_{L^2(0,T;H^1(\Omega))} \leq C_{reg}\Vert z^\mu_i \Vert_{L^2(0,T;H^1(\Omega))}, 
\end{equation*}
whence, since $\{z^\mu_i\}_{\mu>0}$ is a bounded sequence in $L^2(0,T;H^1(\Omega))$ on account of being a weakly convergent sequence therein, there exists a positive constant $C_{seq}$ independent of $m$ and $\mu$, but depending on $\Vert \nabla z_i \Vert_{L^2(0,T;H^1(\Omega))}$, such that 
\begin{equation*}
    \begin{aligned}
        \Vert \nabla R_\mu z^\mu_i - \nabla z_i \Vert_{L^2(Q_T)} \leq C_{seq} \qquad \forall \mu > 0. 
    \end{aligned}
\end{equation*}
In turn, by first taking the limit as $m \to \infty$ and then making $\mu$ vanish, it follows that the right-hand side of \eqref{eq:splitting claim} tends to zero as $\mu \to 0$. To summarise, given any  $\psi \in (L^2(Q_T))^d$, 
\begin{equation*}
    \int_{Q_T} \nabla(R_\mu z^\mu_i - z_i)\cdot \psi \d x \d t \to 0 \qquad \text{as } \mu \to 0, 
\end{equation*}
i.e., $\nabla R_\mu z^\mu_i \rightharpoonup \nabla z_i$ weakly in $L^2(Q_T)$, as required. 
\end{proof}

\section{Vanishing diffusivity and proof of main result}\label{sec:vanishing}

In this section, we take the vanishing diffusivity limit as $\varepsilon \to 0$ in the regularised coupled system \eqref{pdecoupled 0}.

\begin{lem}[Existence for the degenerate system]\label{lem:inviscid sys}
Fix $z_{i,0} \in C^\infty_c(\Omega)$, for $i\in\{1,\dots,M\}$, to be non-negative functions such that $\int_\Omega z_{i,0} \d x = \int_\Omega u_{i,0} \d x$ for $i \in \{1,\dots,M\}$. There exists $z=(z_i)_{i=1}^M$, belonging to the space $\Xi$, which solves the system 
\begin{equation}\label{pde inviscid}
\left\lbrace\begin{aligned}
& \partial_t z_i = \dv [z_i
(\nabla z_i + \nabla L_i(t,x,z_i) + \delta F_i(t,x,z,\nabla z))] \qquad &&\text{in } Q_T, \\ 
& 0 = \nu \cdot [z_i
(\nabla z_i + \nabla L_i(t,x,z_i) + \delta F_i(t,x,z,\nabla z))] \qquad && \text{on } \Sigma_T, \\ 
& z_i(0,\cdot) = z_{i,0} \qquad &&\text{on } \Omega, 
\end{aligned}\right.
\end{equation}
in the weak sense prescribed by Definition \ref{defi:weak sol 0}, with $z_i(0,\cdot) = z_{i,0}$ in $(W^{1,r}(\Omega))'$. Moreover, each $z_i$ is non-negative and conserves its initial mass, and there exists a positive constant $C=C(\Omega,T,d,\delta)$ such that, for $i\in\{1,\dots,M\}$, 
\begin{equation}\label{eq:estimates for inviscid sys}
\Vert z_i \Vert_{L^2(0,T;H^1(\Omega))}^2 + \nm{\de_t z_i}_{X'}
\leq
C \bigg( 1 + \Vert z_{i,0} \Vert_{L^1(\Omega)}^2 + \int_\Omega z_{i,0}\log z_{i,0} \d  x \bigg). 
\end{equation}
\end{lem}

\begin{proof}
The proof is similar to that of Proposition \ref{cor:existence coupled sys}. Recall the weak formulation prescribed by Definition \ref{defi:weak sol 0}, i.e., for any test function $\phi \in C^1(\bar{Q}_T)$, for $i \in \{1,\dots,M\}$, 
\begin{equation}\label{eq:weak sol inviscid sys}
    \begin{aligned}
        \langle \partial_t z_i , \phi \rangle_{X' \times X} +\int_{Q_T} \big[ z_i(\nabla z_i + \nabla L_i(t,x,z_i) + \delta F_i(t,x,z,\nabla z)) \big]\cdot\nabla\phi \d x \d t = 0. 
    \end{aligned}
\end{equation}
For this proof, we define $z^\eps = (z^\eps_i)_{i=1}^M$ to be the weak solution of \eqref{eq:weak sol coupled sys 0} for each $\eps>0$, provided by Proposition \ref{cor:existence coupled sys}. 

Observe from the estimates \eqref{eq:estimates for coupled sys 0} that $\{z^\eps\}_{\eps>0}$ is a bounded sequence in $(L^2(0,T;H^1(\Omega)))^M$, and this bound is independent of $\eps$. Hence, as in the proofs of Lemma \ref{lem:weak compactness} and Proposition \ref{cor:existence coupled sys}, an application of the theorem of Banach--Alaoglu for reflexive spaces and the Aubin--Lions Lemma (\textit{cf.}~\cite[Theorem II.5.16]{boyer2012mathematical}) implies the existence of a subsequence, which we still label as $\{z^{\eps}\}_{\eps>0}$, converging weakly in $(L^2(0,T;H^1(\Omega)))^M$ and strongly in $(L^2(Q_T))^M$ to some $z = (z_i)_{i=1}^M \in (L^2(0,T;H^1(\Omega)))^M$, and such that $\partial_t z^\varepsilon \overset{*}{\rightharpoonup} \partial_t z$ weakly-* in $(X')^M$; see Lemma \ref{lem:time compactness bound 0} for the estimate independent of $\varepsilon$. This latter weak-* convergence in $(X')^M$ is manifestly enough to pass to the limit in the first term of the weak formulation \eqref{eq:weak sol coupled sys 0}, i.e., the duality product. For the remaining term we have 
\begin{equation*}
    \begin{aligned}
        I_\eps := \bigg| \int_{Q_T} \big[ z^\eps_i(\nabla z^\eps_i + \nabla L_i(t,x&,z^\varepsilon_i) + \delta F_i(t,x,z^\eps,\nabla z^\eps)) + \varepsilon \nabla z^\eps_i \big]\cdot\nabla\phi \d x \d t \\ 
    &- \int_{Q_T} \big[ z_i(\nabla z_i + \nabla L_i(t,x,z_i) + \delta F_i(t,x,z,\nabla z)) \big]\cdot\nabla\phi \d x \d t \bigg|,  
\end{aligned}
\end{equation*}
which can be expanded as $I_\eps \leq J_\eps + K_{1,\eps}+K_{2,\eps}$, where 
\begin{equation*}
    \begin{aligned}
        J_\eps := \bigg| \int_{Q_T} z^\eps_i(\nabla z^\eps_i   + \delta F_i & (t,x,z^\eps,\nabla z^\eps)) \cdot \nabla \phi \d x \d t - \int_{Q_T} z_i(\nabla z_i   + \delta F_i(t,x,z,\nabla z)) \cdot \nabla \phi \d x \d t \bigg|, 
    \end{aligned}
\end{equation*}
while 
\begin{equation*}
    K_{1,\eps} := \eps \bigg| \int_{Q_T} \nabla z^\eps_i \cdot\nabla\phi \d x \d t \bigg|, 
\end{equation*}
and 
\begin{equation*}
    K_{2,\eps} := \bigg| \int_{Q_T} z^\eps_i \nabla L_i(t,x,z^\eps_i) \cdot \nabla \phi \d x \d t - \int_{Q_T} z_i \nabla L_i(t,x,z_i) \cdot \nabla \phi \d x \d t \bigg|. 
\end{equation*}
Observe that, using H\"{o}lder's inequality, 
\begin{equation*}
    K_{1,\eps} \leq \varepsilon \Vert \nabla z^\eps_i \Vert_{L^2(Q_T)}\Vert \nabla \phi \Vert_{L^2(Q_T)}, 
\end{equation*}
and note that $\Vert \nabla z^\eps_i \Vert_{L^2(Q_T)} \leq \Vert z^\eps_i \Vert_{L^2(0,T;H^1(\Omega))}$, which is uniformly bounded independently of $\varepsilon$ on account of $\{z^\eps_i\}_{\eps>0}$ being weakly convergent in the space $L^2(0,T;H^1(\Omega))$. It follows that there exists a positive constant $C$ independent of $\varepsilon$ such that 
\begin{equation*}
    K_{1,\eps} \leq \varepsilon C\Vert \nabla \phi \Vert_{L^2(Q_T)} \to 0 \qquad \text{as } \varepsilon \to 0. 
\end{equation*}
As for $K_{2,\varepsilon}$, we estimate 
\begin{equation*}
    \begin{aligned}
        \Vert \nabla L_i(\cdot,\cdot,z^\eps_i) - \nabla L_i(\cdot,\cdot,z_i) \Vert_{L^2(Q_T)}^2 &= \int_0^T \int_\Omega |\nabla W_i*(z^\eps_i(t,\cdot)-z_i(t,\cdot))(x)|^2 \d x \d t \\ 
        &= \int_0^T \int_\Omega \bigg|\int_\Omega \nabla W_i(x-y) (z^\eps_i(t,y)-z_i(t,y)) \d y \bigg|^2 \d x \d t \\ 
        &\leq |\Omega| C_L^2 \int_0^T \int_\Omega \bigg( \int_\Omega |z^\eps_i(t,y)-z_i(t,y)|^2 \d y \bigg) \d x \d t \\ 
        &= |\Omega|^2 C_L^2 \Vert z^\eps_i - z_i \Vert_{L^2(Q_T)}^2 \to 0 \qquad \text{as } \varepsilon \to 0, 
    \end{aligned}
\end{equation*}
where we used the boundedness from \eqref{eq:drift V bound} and Jensen's inequality to obtain the third line, and the Fubini--Tonelli theorem to obtain the final equality, and the strong convergence $z^\eps_i \to z_i$ in $L^2(Q_T)$ to show that the limit vanishes. Thus, 
\begin{equation*}
    \Vert \nabla L_i(\cdot,\cdot,z^\eps_i) - \nabla L_i(\cdot,\cdot,z_i) \Vert_{L^2(Q_T)} \to 0 \qquad \text{as } \varepsilon \to 0, 
\end{equation*}
and hence, being the product of two strongly convergent sequences, we get $z^\eps_i \nabla L_i(\cdot,\cdot,z^\eps_i) \to z_i \nabla L_i(\cdot,\cdot,z_i)$ strongly in $L^2(Q_T)$, whence the Cauchy--Schwarz integral inequality yields that $K_{2,\eps} \to 0$ in the limit as $\eps \to 0$.

We return to $J_\eps$, which we write as $J_\eps \leq J_{1,\eps} + |\delta| J_{2,\eps}$, where 
\begin{equation*}
    J_{1,\eps} := \bigg| \int_{Q_T} (z^\eps_i \nabla z^\eps_i - z_i \nabla z_i) \cdot \nabla \phi \d x \d t  \bigg|, 
\end{equation*}
and 
\begin{equation*}
    J_{2,\eps} := \bigg| \int_{Q_T} \big( z^\eps_i F_i(t,x,z^\eps, \nabla z^\eps) - z_i F_i(t,x,z, \nabla z) \big) \cdot \nabla \phi \d x \d t  \bigg|. 
\end{equation*}
Recall that since $z^\eps_i \to z_i$ strongly in $L^2(Q_T)$ and $\nabla z^\eps_i \rightharpoonup \nabla z_i$ weakly in $(L^2(Q_T))^M$, the product $z^\eps_i \nabla z^\eps_i$ converges weakly to $z_i \nabla z_i$ in $(L^2(Q_T))^M$, whence $J_{1,\eps} \to 0$ in the limit as $\eps$ vanishes. It remains to control $J_{2,\eps}$. We already remarked that $z^\varepsilon_i \to z_i$ strongly in $L^2(Q_T)$. It therefore suffices to show that 
\begin{equation}\label{eq:desired weak convergence in vanish visc}
    F_i(\cdot,\cdot,z^\varepsilon,\nabla z^\varepsilon) \rightharpoonup F_i(\cdot,\cdot,z,\nabla z) \qquad \text{weakly in } L^2(Q_T). 
\end{equation}
In order to verify this, we expand the above terms in $F_i$. Observe that 
\begin{equation*}
    F_i(t,x,z^\eps,\nabla z^\eps) = G^0_i(t,x,z^\eps) + \sum_{j=1}^M G^1_{ij}(t,x,z^\eps) \nabla z^\eps_j, 
\end{equation*}
and so, using the mean value theorem and the uniform boundedness of the derivatives of $(G_i)_{i=1}^M$, we see that 
\begin{equation}\label{eq:mean value ineq vanish visc strong conv}
    \Vert G^0_i(\cdot,\cdot,z^\eps) - G^0_i(\cdot,\cdot,z) \Vert_{L^2(Q_T)} \leq \Vert \nabla_z G^0_i \Vert_{L^\infty(Q_T \times \mathbb{R}^M)} \max_{i \in \{1,\dots,M\}} \Vert z^\eps_i - z_i \Vert_{L^2(Q_T)}, 
\end{equation}
and the right-hand side vanishes due to the strong convergence $z^\varepsilon_i \to z_i$ in $L^2(Q_T)$. For the other term, we have the following claim, to be proved later. 
\begin{claim}\label{claim:weak convergence of sum in vanishing viscosity}
\begin{equation*}
    \sum_{j=1}^M G^1_{ij}(\cdot,\cdot,z^\eps) \nabla z^\eps_j \rightharpoonup \sum_{j=1}^M G^1_{ij}(\cdot,\cdot,z) \nabla z_j \qquad \text{weakly in } L^2(Q_T). 
\end{equation*}
\end{claim}
It then follows immediately from the previous claim and \eqref{eq:mean value ineq vanish visc strong conv} that we obtain the desired weak convergence \eqref{eq:desired weak convergence in vanish visc}. Now, due to the strong convergence $z^\varepsilon_i \to z_i$ strongly in $L^2(Q_T)$, it follows that we have weak convergence of the product, i.e., 
\begin{equation*}
    z^\eps_i F_i(\cdot,\cdot,z^\eps, \nabla z^\eps) \rightharpoonup z_i F_i(\cdot,\cdot,z, \nabla z) \qquad \text{weakly in } L^2(Q_T), 
\end{equation*}
whence it follows that $J_{2,\eps} \to 0$ as $\eps \to 0$. The limit function $z \in L^2(0,T;H^1(\Omega))$ thereby satisfies the weak formulation \eqref{eq:weak sol inviscid sys}, and the estimates \eqref{eq:estimates for inviscid sys} follow from the weak (and weak-* for the bound in $X'$) lower semicontinuity of the norms considered in the estimate \eqref{eq:estimates for coupled sys 0}. The non-negativity and mass conservation follow again directly from Lemma \ref{lem:preserve sign and mass}, while the attainment of the initial data in the dual Sobolev space follows from a direct application of Lemma \ref{lem:preserve initial data}. 
\end{proof}

\begin{proof}[Proof of Claim \ref{claim:weak convergence of sum in vanishing viscosity}]
Fix $j \in \{1,\dots,M\}$ in the sum, and by testing against any $\psi \in (C^\infty_c(Q_T))^d$, we obtain 
\begin{equation*}
    \begin{aligned}
        \bigg| \int_{Q_T} \big( G^1_{ij}(t,x,z^\varepsilon)\nabla z^\eps_j - G^1_{ij}(t,x,z) \nabla z_j \big) \cdot \psi \d x \d t \bigg|         \leq \bigg| \int_{Q_T} \big( G^1_{ij}(t,x,z^\varepsilon) - G^1_{ij}(t,x,z) \big)\nabla z^\eps_j \cdot \psi \d x \d t \bigg|& \\ 
        + \bigg| \int_{Q_T} G^1_{ij}(t,x,z) ( \nabla z^\eps_j - \nabla z_j ) \cdot \psi \d x \d t \bigg|&, 
    \end{aligned}
\end{equation*}
and the right-hand side of the above is bounded by 
\begin{equation*}
    \Vert G^0_i(\cdot,\cdot,z^\eps) - G^0_i(\cdot,\cdot,z) \Vert_{L^2(Q_T)} \Vert \nabla z^\eps_j \Vert_{L^2(Q_T)} \Vert \psi \Vert_{L^\infty(Q_T)} + \bigg| \int_{Q_T} ( \nabla z^\eps_j - \nabla z_j ) \cdot  (G^1_{ij})^T(t,x,z) \psi \d x \d t \bigg|, 
\end{equation*}
where $(G^1_{ij})^T$ is the transpose of the matrix $G^1_{ij}$. The first term on the right-hand side of the above vanishes as $\eps \to 0$ using the strong convergence $z^\eps_i \to z_i$ strongly in $L^2(Q_T)$ and an estimate analogous to the one in \eqref{eq:mean value ineq vanish visc strong conv} and the uniform boundedness of $\Vert \nabla z^\eps_j \Vert_{L^2(Q_T)}$ independently of $\eps$, on account of the weak convergence $z^\eps_i \rightharpoonup z_i$ in $L^2(0,T;H^1(\Omega))$. The second term on the right-hand side of the above also vanishes, on account of $(G^1_{ij})^T \psi \in L^\infty(Q_T) \subset L^2(Q_T)$ and the weak convergence $z^\eps_j \rightharpoonup z_j$ in $L^2(0,T;H^1(\Omega))$. 

We have therefore shown that, given any $\psi \in (C^\infty_c(Q_T))^d$, we have 
\begin{equation}\label{eq:for second claim checkpoint}
    \begin{aligned}
        \bigg| \int_{Q_T} & \big( G^1_{ij}(t,x,z^\varepsilon)\nabla z^\eps_j - G^1_{ij}(t,x,z) \nabla z_j \big) \cdot \psi \d x \d t \bigg| \to 0 \qquad \text{as } \eps \to 0. 
    \end{aligned}
\end{equation}
Now fix any $\psi \in (L^2(Q_T))^d$. By density of $C^\infty_c(Q_T)$ in $L^2(Q_T)$, there exists a sequence $(\psi_m)_{m\in\mathbb{N}}$ of elements of $(C^\infty_c(Q_T))^d$ converging strongly to $\psi$ in the sense of $(L^2(Q_T))^d$. Then, by splitting the term $| \int_{Q_T} \big( G^1_{ij}(t,x,z^\varepsilon)\nabla z^\eps_j - G^1_{ij}(t,x,z) \nabla z_j \big) \cdot \psi \d x \d t|$ as 
\begin{equation*}
    \begin{aligned}
        \bigg| \int_{Q_T} \big( G^1_{ij}(t,x,z^\varepsilon)\nabla z^\eps_j &- G^1_{ij}(t,x,z) \nabla z_j \big) \cdot (\psi - \psi_m) \d x \d t \bigg| \\ 
        &+ \bigg| \int_{Q_T} \big( G^1_{ij}(t,x,z^\varepsilon)\nabla z^\eps_j - G^1_{ij}(t,x,z) \nabla z_j \big) \cdot \psi_m \d x \d t \bigg|, 
        \end{aligned}
        \end{equation*}
        which is itself bounded by 
        \begin{equation*}
        \begin{aligned}
        \Vert G^1_{ij} \Vert_{L^\infty(Q_T\times\mathbb{R}^M)} (\Vert \nabla z^\eps_j \Vert_{L^2(Q_T)} &+ \Vert \nabla z_j \Vert_{L^2(Q_T)}) \Vert \psi - \psi_m \Vert_{L^2(Q_T)} \\ 
        &+ \bigg| \int_{Q_T} \big( G^1_{ij}(t,x,z^\varepsilon)\nabla z^\eps_j - G^1_{ij}(t,x,z) \nabla z_j \big) \cdot \psi_m \d x \d t \bigg|. 
    \end{aligned}
\end{equation*}
Using the weak convergence in $L^2(0,T;H^1(\Omega))$ to deduce the boundedness of $\Vert \nabla z^\eps_j \Vert_{L^2(Q_T)}$, independently of $\eps$, we now use an argument identical to the one in the proof of Claim \ref{claim:annoying diagonal} to deduce that 
\begin{equation*}
    \begin{aligned}
        \bigg| \int_{Q_T} & \big( G^1_{ij}(t,x,z^\varepsilon)\nabla z^\eps_j - G^1_{ij}(t,x,z) \nabla z_j \big) \cdot \psi \d x \d t \bigg| \to 0 \qquad \text{as } \eps \to 0, 
    \end{aligned}
\end{equation*}
for any $\psi \in (L^2(Q_T))^d$, as required. 
\end{proof}

This completes the vanishing diffusivity procedure. It remains to relax the assumption on the initial data, which we also do by a limiting strategy. This is contained below, which is the proof of the main result. 

\begin{proof}[Proof of Theorem \ref{theor:main}]
Begin by assuming $p \in (1,\infty)$. In view of the density of $C^\infty_c(\Omega)$ in $L^p(\Omega)$, given $u_{i,0} \in L^p(\Omega)$ as in the statement of the theorem, there exists a sequence $\{z^m_{i,0}\}_{m\in\mathbb{N}}$ of elements of $C^\infty_c(\Omega)$, all of which are non-negative and satisfy the initial mass assumption $\int_\Omega z_{i,0}^m \d x = \int_\Omega u_{i,0} \d x$ as per Remark \ref{rem:fixing initial data for pdefrozen}, such that 
\begin{equation*}
    \Vert u_{i,0} - z^m_{i,0} \Vert_{L^p(\Omega)} \to 0 \qquad \text{as } m \to \infty, 
\end{equation*}
for each $i\in\{1,\dots,M\}$. An explicit construction for such an approximating sequence is to set 
\begin{equation*}
    z_{i,0}^m(x) := \big((u_{i,0}\mathds{1}_{\Omega})*\rho_m(x)\big)\eta_{\sigma(m)}(x) \frac{\Vert u_{i,0}\Vert_{L^1(\Omega)}}{\Vert \big((u_{i,0}\mathds{1}_\Omega)*\rho_m\big)\eta_{\sigma(m)} \Vert_{L^1(\Omega)}} \qquad \forall x \in \Omega,
\end{equation*}
where $\rho_m$ is the usual Friedrichs mollifier, $\eta_m$ is a smooth non-negative cutoff function chosen such that $\eta \equiv 1$ on $\{x \in \Omega : d(x,\partial\Omega) \geq 1/m \}$ and $\eta \equiv 0$ outside $\{x \in \Omega : d(x,\partial\Omega) \geq 1/2m \}$, and $\sigma$ is an appropriate scaling function depending on $p,d$ (i.e.~$\sigma(m)=m^q$ for some exponent $q=q(p,d)$ suitably chosen).

For the rest of this proof, for each $m\in\mathbb{N}$, we define $z^m = (z^m_i)_{i=1}^M$ to be the weak solution of \eqref{eq:weak sol inviscid sys} with initial data $z^m_{i,0}$, for $i\in\{1,\dots,M\}$, provided by Lemma \ref{lem:inviscid sys}. Observe from the estimates \eqref{eq:estimates for inviscid sys} that $\{z^m\}_{m\in\mathbb{N}}$ is a bounded sequence in $(L^2(0,T;H^1(\Omega)))^M$, and this bound is independent of $m$. Hence, as in the proofs of Lemma \ref{lem:weak compactness}, Proposition \ref{cor:existence coupled sys}, and Lemma \ref{lem:inviscid sys}, an application of the theorem of Banach--Alaoglu for reflexive spaces and the Aubin--Lions Lemma (\textit{cf.}~\cite[Theorem II.5.16]{boyer2012mathematical}) implies the existence of a subsequence, which we still label as $\{z^{m}\}_{m\in\mathbb{N}}$, converging weakly in $(L^2(0,T;H^1(\Omega)))^M$ and strongly in $(L^2(Q_T))^M$ to some $u = (u_i)_{i=1}^M \in (L^2(0,T;H^1(\Omega)))^M$, and such that $\partial_t z^m \overset{*}{\rightharpoonup} \partial_t u$ weakly-* in $(X')^M$. This latter weak-* convergence in $(X')^M$ is manifestly enough to pass to the limit in the first term of the weak formulation \eqref{eq:weak sol inviscid sys}, i.e., the duality product, and for the final term, we have 
\begin{equation*}
    \begin{aligned}
        I_m := \bigg| \int_{Q_T} \big[ z^m_i(\nabla z^m_i + \nabla L_i(t,x&,z^m_i) + \delta F_i(t,x,z^m,\nabla z^m)) \big]\cdot\nabla\phi \d x \d t \\ 
    &- \int_{Q_T} \big[ u_i(\nabla u_i + \nabla L_i(t,x,z_i) + \delta F_i(t,x,u,\nabla u)) \big]\cdot\nabla\phi \d x \d t \bigg|,  
\end{aligned}
\end{equation*}
Following a procedure identical to that of the proof of Lemma \ref{lem:inviscid sys}, we deduce that $I_m \to 0$ as $m\to\infty$. It follows that the limit function $z \in (L^2(0,T;H^1(\Omega)))^M$ thereby satisfies the weak formulation \eqref{eq:weak sol inviscid sys}. 

The estimates \eqref{eq:estimates for main theorem i}-\eqref{eq:estimates for main theorem ii} follow from the weak (and weak-* for the bound in $X'$) lower semicontinuity of the norms considered in the estimate \eqref{eq:estimates for inviscid sys}. Indeed, note that we chose $z^m_{i,0}$ at the start of the proof such that $\Vert z_{i,0}^m \Vert_{L^1(\Omega)} = \Vert u_{i,0} \Vert_{L^1(\Omega)}$ for all $m\in\mathbb{N}$. Meanwhile, the function $f:x \mapsto (x \log x) \mathds{1}_{[0,\infty)}(x)$ is continuous and satisfies the global bound $|f(x)| \leq C(1 + |x|^{p})$ for some positive constant $C$ depending only on $p$. As a result, by a consequence of the Generalised Dominated Convergence Theorem, $f$ maps $L^p(\Omega)$ continuously into $L^1(\Omega)$. It therefore follows that, since $z^m_{i,0} \to u_{i,0}$ strongly in $L^p(\Omega)$, 
\begin{equation*}
    z^m_{i,0} \log z^m_{i,0} \to u_{i,0} \log u_{i,0} \qquad \text{in } L^1(\Omega), 
\end{equation*}
and hence $\lim_{m\to\infty} \int_\Omega z_{i,0}^m \log z_{i,0}^m \d x = \int_\Omega u_{i,0} \log u_{i,0} \d x$, as required. 

The non-negativity and mass conservation follow again from Lemma \ref{lem:preserve sign and mass}. The convergence to the initial data in the sense of the third point of Definition \ref{defi:weak sol 0} follows from a direct application of Lemma \ref{lem:preserve initial data}. 

In the case $p=\infty$, we have $u_{i,0} \in L^q(\Omega)$ for any finite $q \in (1,\infty)$, and we can follow the same argument as before, approximating the initial data in $L^q(\Omega)$---as opposed to in $L^\infty(\Omega)$, which (in general) cannot be done using smooth compactly supported functions. The proof is complete. 
\end{proof}

\section*{Acknowledgements}
M.~Bruna was supported by Royal Society University Research Fellowship (grant no.~URF/R1/180040). S.~Schulz was supported by the Royal Society Award (RGF/EA/181043). This work is the result of a collaboration that started at the workshop ``Recent Advances in Degenerate Parabolic Systems with Applications to Mathematical Biology'' (jointly funded by the BOUM programme of SMAI and by the support of the ERC advanced grant Adora), which took place at the Laboratoire Jacques-Louis Lions (UPMC Paris VI) in February 2020. 

\begin{appendices}
\section{Appendices}\label{sec:appendix united}

\subsection{Proof of Lemma \ref{lem:existence of viscous approximates 0}: Existence of solutions to the regularised frozen system}\label{sec:appendix proof russian}

\begin{proof}[Proof of Lemma \ref{lem:existence of viscous approximates 0}]
We begin by recasting the problem as one with homogeneous initial condition, i.e., defining $v_i(t,x) := z_i(t,x) - z_{i,0}(x)$ for each $i \in \{1,\dots,M\}$, problem \eqref{eq:pdefrozen} now reads as 
\begin{equation}\label{eq:neumann zero initial 0}
\left\lbrace\begin{aligned}
&\de_t v_i = \dv\left[
(v_i+z_{i,0})(\nabla (v_i+z_{i,0}) + \nabla L_i(t,x,v_i+z_{i,0}) + \delta \bar{F}_i) + \varepsilon\nabla (v_i+z_{i,0}) \right]
\quad && \text{in } Q_T,
\\
&0 = \nu \cdot\left[
(v_i+z_{i,0})(\nabla (v_i+z_{i,0}) + \nabla L_i(t,x,v_i+z_{i,0}) + \delta \bar{F}_i) + \varepsilon\nabla (v_i+z_{i,0})
\right]
\quad && \text{on } \Sigma_T,
\\
&v_i(0,\cdot) = 0 \quad
&& \text{on } \Omega. 
\end{aligned}\right.
\end{equation}
By expanding the flux term, the above system may be rewritten as 
\begin{equation*}
\left\lbrace\begin{aligned}
\de_t v_i = \dv\big[&(v_i + z_{i,0} + \varepsilon) \nabla v_i + (\nabla z_{i,0} + \nabla L_i(t,x,v_i+z_{i,0}) + \delta \bar{F}_i) v_i \\ 
&+ [z_{i,0}(\nabla z_{i,0} + \nabla L_i(t,x,v_i+z_{i,0}) + \delta \bar{F}_i) + \varepsilon \nabla z_{i,0}) ] \big] 
 \qquad &&\text{in } Q_T, \\
0 = \nu \cdot\big[&(v_i + z_{i,0} + \varepsilon) \nabla v_i + (\nabla z_{i,0} + \nabla L_i(t,x,v_i+z_{i,0}) + \delta \bar{F}_i) v_i \\ 
&+ [z_{i,0}(\nabla z_{i,0} + \nabla L_i(t,x,v_i+z_{i,0}) + \delta \bar{F}_i) + \varepsilon \nabla z_{i,0}) ]
\big] \qquad &&\text{on } \Sigma_T 
\\
v_i(0,\cdot) = 0& \qquad && \text{on } \Omega, 
\end{aligned}\right.
\end{equation*}
from which the requirement of the compatibility condition \eqref{eq:compatibility 0} is apparent by considering the no-flux boundary condition at the initial time $t=0$. 

Observe that the problem \eqref{eq:neumann zero initial 0} may be rewritten in the prototypical diagonal non-divergence form 
\begin{equation}\label{eq:as in LSU 0}
    \left\lbrace\begin{aligned}
    &\partial_t v_i + Lv_i = 0 \qquad &&\text{in } Q_T, \\ 
    &0 = \nu\cdot \big[ (v_i + z_{i,0}+\varepsilon) \nabla v_i + (\nabla z_{i,0} + \nabla L_i(t,x,v_i+z_{i,0}) +\delta \bar{F}_i) v_i \big] \qquad && \text{on } \Sigma_T, \\ 
    &v_i(0,\cdot) = 0 \qquad &&\text{on } \Omega, 
    \end{aligned}\right.
\end{equation}
where we define the operator $L$ by 
\begin{equation}\label{eq:unif elliptic op 0}
    Lv = - \sum\nolimits_{j,k=1}^d a^{jk}(v) \partial^2_{jk} v + b(t,x,v,\nabla v), 
\end{equation}
with 
\begin{equation}\label{eq:coefficients of elliptic op 0}
    \left\lbrace\begin{aligned}
    a^{jk}(v) = &(v + z_{i,0}  + \varepsilon) \delta^{jk}  , \\ 
    b(t,x,v,p) = &-|p|^2 - (2\nabla z_{i,0} + \nabla L_i(t,x,v+z_{i,0}) + \delta \bar{F}_i)\cdot p \\ 
    &- \dv[\nabla z_{i,0} + \nabla L_i(t,x,v+z_{i,0}) + \delta \bar{F}_i]v - \dv[z_{i,0} \nabla L_i(t,x,v+z_{i,0})] \\ 
    &- \dv[z_{i,0}(\nabla z_{i,0} + \delta \bar{F}_i) + \varepsilon \nabla z_{i,0}]. 
    \end{aligned} \right.
\end{equation}
We verify that the conditions (a), (b), and (c) of \cite[Theorem 7.4 in Section 7 of Chapter 5]{ladyzhenskaia1988linear} (which rely on estimates (7.4)--(7.6), (7.15), (7.34), and (7.36) therein) are satisfied by the scalar equation for each $i$ in question (since the system is diagonal) in the form with homogeneous initial condition.

To begin with, notice from Lemma \ref{lem:sign preserved regularised} that the operator \eqref{eq:unif elliptic op 0} is uniformly elliptic, since for $|v| \leq \Lambda$, 
\begin{equation}\label{eq:ellipticity}
    \varepsilon|\xi|^2 \leq a^{jk}(v) \xi_j \xi_k \leq (\Lambda+ \Lambda_0 +\varepsilon)|\xi|^2 \qquad \text{for } \xi \in \mathbb{R}^d, 
\end{equation}
where $\Lambda_0$ was given in \eqref{eq:M zero}. Additionally, all derivatives of $a^{jk}$ are uniformly bounded and $b$ is subquadratic in its final argument, by which we mean: for $(t,x) \in \bar{Q}_T$ and $|v| \leq \Lambda$ and arbitrary $p\in\mathbb{R}$, 
\begin{equation*}
    |\partial_v a^{jk}| + |\partial_{vv} a^{jk}| \leq 1, 
\end{equation*}
and, using a combination of the Cauchy--Young inequality along with \eqref{eq:boundedness of bar F z}-\eqref{eq:M zero} and \eqref{eq:drift V bound}, along with \eqref{eq:coefficients of elliptic op 0} and $\varepsilon \in (0,1)$, 
\begin{equation*}
    |b(t,x,v,p)| \leq C \big( 1+|p|^2 \big), 
\end{equation*}
for some positive constant $C = C(\delta , C_L , \Lambda , \Lambda_0 , \Lambda_{\bar{z}})$. Furthermore, making $C$ larger if necessary, we have, again for $|v| \leq \Lambda$, 
\begin{equation*}
    |\partial_p b|(1+|p|) + |\partial_v b|  \leq C|p|(1+|p|) + C(1+|p|) \leq C \big(1+|p|^2 \big), 
\end{equation*}
while, again using \eqref{eq:boundedness of bar F z}-\eqref{eq:M zero} and \eqref{eq:drift V bound}, 
\begin{equation*}
    |\partial_t b| \leq C \big(1+|p|^2\big). 
\end{equation*}

In view of \eqref{eq:boundedness of bar F z}, for each fixed $(t,z,p) \in [0,T]\times \mathbb{R}^M \times\mathbb{R}^{d \times M}$, we have that $b(t,\cdot,z,p) \in C^1(\bar{\Omega})$, and so it is Lipschitz with respect to the spatial variable. Hence, all the conditions of \cite[Theorem 7.4 in Section 7 of Chapter 5]{ladyzhenskaia1988linear} are satisfied, and so there exists a unique weak solution (in the sense of Definition \ref{defi:weak sol regularised}) of the no-flux problem with homogeneous initial condition \eqref{eq:as in LSU 0}, and thus also \eqref{eq:neumann zero initial 0}, living in the space of functions with continuous derivatives of the type $\partial_t^{n_1}\partial_{x_j}^{n_2}$ for $2n_1 + n_2 < 3$; denoted by $H^{3,\frac{3}{2}}(\bar{Q}_T)$ in \cite{ladyzhenskaia1988linear}. Note in particular that this class of functions is contained in $C^{2,1}(\bar{Q}_T)$. It therefore follows that there exists $v = (v_i)_{i=1}^M \in C^{2,1}(\bar{Q}_T)$ solving the no-flux problem with homogeneous initial condition \eqref{eq:neumann zero initial 0} as a pointwise equality between continuous functions. Uniqueness in the class $C^{2,1}(\bar{Q}_T)$ can be proved by standard methods. 
\end{proof}

\begin{rem}\label{rem:unpleasant technicalities}
Notice that the requirement that $F_i$ be $C^2$ in each argument is clear from the previous estimates on the derivatives of $b$, since, for instance, one needs to bound $\partial_t \dv \big( F_i(t,x,\bar{z}(t,x),\nabla\bar{z}(t,x)) \big)$. 
\end{rem}

\subsection{Proof of Lemma \ref{lem:quantitativeH2}: Quantitative second derivative estimate for the regularised frozen system}\label{appendix:H2bounds}

Throughout this section, it will be used that the regularised frozen system \eqref{eq:pdefrozen} is diagonal. We therefore consider the single equation 
\begin{equation}\label{eq:diag sys to single eq}
\left\lbrace\begin{aligned}
&\partial_t w = \dv\left[
w(\nabla w + \nabla L + \delta \bar{F})+\varepsilon \nabla w
\right]
\qquad && \text{in } Q_T, \\ 
&0 = \nu \cdot\left[
w(\nabla w + \nabla L + \delta \bar{F})
+\varepsilon \nabla w \right]
\qquad && \text{on } \Sigma_T, 
\\
&w(0,\cdot) = w_{0} 
\qquad && \text{on } \Omega,
\end{aligned}\right.
\end{equation}
where we omitted the $i$ subscripts for clarity of presentation, and use the notation established in \eqref{eq:F bar notation}. Recall that (as per Remark \ref{rem:Mbarz bound}) $\bar{F}$ is bounded in $C^2$-norm and (as per Remark \ref{rem:fixing initial data for pdefrozen}) $w_0 \in C^\infty_c(\Omega)$. We already know from Lemmas \ref{lem:existence of viscous approximates 0} that there exists a non-negative $w \in C^{2,1}(\bar{Q}_T)$ solving the above equation in the classical sense, which satisfies the estimates of Lemma \ref{lem:energy estimates 0}. 

To begin with, we prove a quantitative $L^\infty$-bound on the solution of the regularised frozen system \eqref{eq:pdefrozen}. 

\begin{lem}[$L^\infty$-bound]\label{cor:Linfty bound}
Suppose that $z=(z_i)_{i=1}^M$ is a $C^{2,1}(\bar{Q}_T)$ solution of the regularised frozen system \eqref{eq:pdefrozen}. Then, there holds, for each $i \in \{1,\dots,M\}$, 
\begin{equation}\label{eq:w Linfty bound}
	\Vert z_i \Vert_{L^\infty(Q_T)} \leq C(\varepsilon,\delta,T,\Omega,C_L,\Vert z_{i,0} \Vert_{L^\infty(\Omega)},\Vert \bar{F}_i \Vert_{L^\infty(Q_T)}), 
\end{equation}
where the right-hand side is a positive quantity depending only on the parameters in its parentheses. 
\end{lem}
\begin{proof}
We only write the proof for $d \geq 2$. By multiplying \eqref{eq:diag sys to single eq} by the continuously differentiable function $q w^{q-1}$ (for any finite $q \geq 2$) and by closely following the argument in \cite[Proof of Theorem 3.1]{KimZhangII}, we obtain, for every $t \in [0,T]$, 
\begin{equation}\label{eq:before gagliardo}
	\begin{aligned}
	\frac{d}{dt}\int_\Omega w^q \d x + \varepsilon \int_\Omega |\nabla(w^{q/2})|^2 \d x &\leq \frac{q(q-1)}{\varepsilon} \int_\Omega w^q ( |\nabla L|^2 +\delta^2 |\bar{F}|^2 ) \d x \\ 
	&\leq C(\varepsilon,\delta,C_L,\Vert w_0 \Vert_{L^1(\Omega)},\Vert \bar{F} \Vert_{L^\infty(Q_T)}) \bigg( q^2 \int_\Omega w^q \d x \bigg), 
\end{aligned}
\end{equation}
where we used $|\nabla L(\cdot,\cdot,w)| \leq \Vert V \Vert_{C^1(\mathbb{R}^{d+1})} + |(\nabla W * w(t,\cdot))| \leq \Vert V \Vert_{C^1(\mathbb{R}^{d+1})} + \Vert W \Vert_{C^1(\mathbb{R}^{d})}\Vert w_0 \Vert_{L^1(\Omega)}$, using H\"{o}lder's inequality for the convolution, to get the bound $\Vert \nabla L \Vert_{L^\infty(Q_T)} \leq C_L(1+\Vert w_0 \Vert_{L^1(\Omega)})$. We now use the Gagliardo--Nirenberg inequality to write 
\begin{equation*}
	\int_\Omega w^q \d x = \Vert w^{q/2} \Vert_{L^2(\Omega)}^2 \leq C(\Omega) \big( \Vert \nabla(w^{q/2}) \Vert_{L^2(\Omega)}^{2\gamma} \Vert w^{q/2} \Vert_{L^1(\Omega)}^{2(1-\gamma)} + \Vert w^{q/2} \Vert_{L^1(\Omega)}^2 \big) 
\end{equation*}
for $d \geq 3$, where $\gamma = d/(d+2)$. The case $d=2$ can be dealt with analogously using Ladyzhenskaya's inequality. Then, using the above and the weighted Young inequality, the right-hand side of \eqref{eq:before gagliardo} is bounded by 
\begin{equation*}
	\begin{aligned}
	\frac{\varepsilon}{2}\Vert \nabla(w^{q/2}) \Vert_{L^2(\Omega)}^2 + C(\varepsilon,\delta,\Omega,d,C_L,\Vert w_0 \Vert_{L^1(\Omega)},\Vert \bar{F} \Vert_{L^\infty(Q_T)}) q^{\frac{2\gamma}{1-\gamma}}\Vert w^{q/2} \Vert_{L^1(\Omega)}^2. 
	\end{aligned}
\end{equation*}
Using the Gagliardo--Nirenberg and weighted Young inequalities again to bound $\Vert \nabla(w^{q/2}) \Vert_{L^2(\Omega)}$ from below by $C_1(\Omega,d)\Vert w^{q/2} \Vert_{L^2(\Omega)}^2 - C_2(\Omega,d)\Vert w^{q/2} \Vert_{L^1(\Omega)}^2$ and returning to \eqref{eq:before gagliardo}, we obtain, for every $t \in [0,T]$, 
\begin{equation*}
	\begin{aligned}
	\frac{d}{dt}\int_\Omega w^q \d x + \int_\Omega w^q \d x &\leq C(\varepsilon,\delta,\Omega,d,C_L,\Vert w_0 \Vert_{L^1(\Omega)},\Vert \bar{F} \Vert_{L^\infty(Q_T)}) \bigg( q^{2+\frac{2\gamma}{1-\gamma}}\bigg( \int_\Omega w^{q/2} \d x \bigg)^2  + q^2 \bigg), 
\end{aligned}
\end{equation*}
and the above holds for every $q \in [2,\infty)$. We now conclude using the iterative result \cite[Lemma 3.2]{KimZhangII}, as per \cite[Proof of Theorem 3.1]{KimZhangII}. 
\end{proof}

With the previous estimate in hand, we proceed to the bound on the second derivative. 

\begin{proof}[Proof of Lemma \ref{lem:quantitativeH2}]
We neglect the drift terms for the time being, for clarity of presentation, though we show how to treat them at the end of the proof. Hence, we restrict our focus to the higher regularity of the single equation 
\begin{equation*}
\left\lbrace\begin{aligned}
&\partial_t w = \dv\left[
w(\nabla w + \delta \bar{F})+\varepsilon \nabla w
\right]
\qquad && \text{in } Q_T, \\ 
&0 = \nu \cdot\left[
w(\nabla w + \delta \bar{F})
+\varepsilon \nabla w \right]
\qquad && \text{on } \Sigma_T, 
\\
&w(0,\cdot) = w_{0} 
\qquad && \text{on } \Omega,
\end{aligned}\right.
\end{equation*}
Recall that, by Lemma \ref{lem:energy estimates 0}, there holds 
\begin{equation}\label{eq:L2H1 bound that you already know}
	\Vert w \Vert^2_{L^2(0,T;H^1(\Omega))} \leq C(\delta,T,|\Omega|,\Vert w_0 \Vert_{L^\infty(\Omega)},\Vert \bar{F} \Vert_{L^\infty(Q_T)}), 
\end{equation}
where we bound the right-hand side of \eqref{eq:L2H1 bound 0} using $\Vert \bar{F} \Vert^2_{L^2(Q_T)} \leq T |\Omega| \Vert \bar{F} \Vert^2_{L^\infty(Q_T)}$, $\Vert w_0 \Vert_{L^1(\Omega)} \leq |\Omega|\Vert w_0 \Vert_{L^\infty(\Omega)}$, and $\int_\Omega w_0 \log w_0 \d x \leq \int_\Omega w_0 (\log w_0)_+ \d x \leq C|\Omega|(1+\Vert w_0 \Vert_{L^\infty(\Omega)}^2)$ for some universal constant $C$.

\textbf{Step I: } Define the new (non-negative) function 
\begin{equation}\label{eq:new psi def}
	\psi(t,x) := \frac{1}{2}w(t,x)^2 + \varepsilon w(t,x) \qquad \forall (t,x) \in \bar{Q}_T, 
\end{equation}
from which it follows that we may write, since $w$ is non-negative itself, 
\begin{equation*}
	w(t,x) = -\varepsilon + \sqrt{\varepsilon^2 + 2\psi(t,x)} \qquad \forall (t,x) \in \bar{Q}_T, 
\end{equation*}
and we note the formula 
\begin{equation}\label{eq:formula dt w}
	\partial_t w = \frac{\partial_t \psi}{\sqrt{\varepsilon^2+2\psi}} \qquad \text{in } \bar{Q}_T. 
\end{equation}
The evolution equation for $w$ may be rewritten as 
\begin{equation}\label{eq:new evolution}
\left\lbrace\begin{aligned}
&\partial_t \psi = (\varepsilon^2+2\psi)^{1/2}\dv\left[\nabla \psi + \delta w \bar{F} \right]
\qquad && \text{in } Q_T, \\ 
&0 = \nu \cdot\left[
\nabla \psi + \delta w \bar{F} \right]
\qquad && \text{on } \Sigma_T, 
\\
&\psi(0,\cdot) = \frac{1}{2}w_{0}^2 + \varepsilon w_0 
\qquad && \text{on } \Omega. 
\end{aligned}\right.
\end{equation}
We now test the above against $\dv [\nabla \psi + \delta w \bar{F}]$. An integration by parts on the left-hand side (with no boundary terms because of the no-flux condition) yields 
\begin{equation*}
	-\int_\Omega \partial_t \nabla \psi \cdot [\nabla \psi + \delta w \bar{F}] \d x = \int_\Omega (\varepsilon^2+2\psi)^{1/2}(\Delta \psi + \delta \dv(w \bar{F}) )^2 \d x \geq \varepsilon \int_\Omega (\Delta \psi + \delta \dv(w \bar{F}) )^2 \d x, 
\end{equation*}
where the final inequality follows from the non-negativity of $\psi$. Thus, 
\begin{equation}\label{eq:unpleasant but making progress}
	\frac{1}{2}\frac{d}{dt}\int_\Omega |\nabla \psi + \delta w\bar{F}|^2 \d x + \varepsilon \int_\Omega (\Delta \psi + \delta \dv(w \bar{F}) )^2 \d x \leq - \delta \int_\Omega \partial_t (w \bar{F}) \cdot [\nabla \psi + \delta w \bar{F}] \d x. 
\end{equation}
The right-hand side may be expanded as 
\begin{equation}\label{eq:split the unpleasant}
	\begin{aligned}
		-\delta \int_\Omega & (\partial_t w) \bar{F} \cdot [\nabla \psi + \delta w \bar{F}] \d x - \delta \int_\Omega w (\partial_t \bar{F}) \cdot [\nabla \psi + \delta w \bar{F}] \d x, 
		\end{aligned}
		\end{equation}
		and, by the Cauchy--Schwarz integral inequality, the second term of the above is bounded above by 
		\begin{equation*}
		\begin{aligned}
		\delta^2 |\Omega|^{1/2} \Vert w \Vert_{L^\infty({Q}_T)}^2 \Vert \partial_t \bar{F} \Vert_{L^2(\Omega)} \Vert \bar{F} \Vert_{L^\infty(Q_T)} + \delta \Vert w \Vert_{L^\infty(\Omega)}\Vert \partial_t \bar{F} \Vert_{L^2(\Omega)}\Vert \nabla \psi \Vert_{L^2(\Omega)}. 
	\end{aligned}
\end{equation*}
The first term in \eqref{eq:split the unpleasant} may be rewritten, using the formula \eqref{eq:formula dt w} and the equation \eqref{eq:new evolution}, as 
\begin{equation*}
	\begin{aligned}
		-\delta \int_\Omega & \dv[\nabla \psi + \delta w \bar{F}] \bar{F} \cdot [\nabla \psi + \delta w \bar{F}] \d x = -\delta \int_\Omega (\Delta \psi + \delta \dv(w \bar{F})) \bar{F} \cdot [\nabla \psi + \delta w \bar{F}] \d x. 
	\end{aligned}
\end{equation*}
Using the Cauchy--Schwartz integral inequality and the Young inequality, the right-hand side of the above is bounded above by 
\begin{equation*}
	\begin{aligned}
		\frac{\varepsilon}{2}\int_\Omega (\Delta \psi + \delta \dv(w \bar{F}))^2 \d x  + \frac{\delta^2}{2\varepsilon} \Vert \bar{F} \Vert^2_{L^\infty(Q_T)} \int_\Omega |\nabla \psi + \delta w \bar{F}|^2 \d x. 
	\end{aligned}
\end{equation*}
Returning to \eqref{eq:unpleasant but making progress}, we therefore have 
\begin{equation*}
	\begin{aligned}
		\frac{1}{2}\frac{d}{dt} \int_\Omega |\nabla \psi + \delta w\bar{F}|^2 \d x + \frac{\varepsilon}{2} \int_\Omega & (\Delta \psi + \delta \dv(w \bar{F}) )^2 \d x \leq \delta^2 |\Omega|^{1/2} \Vert w \Vert_{L^\infty({Q}_T)}^2 \Vert \partial_t \bar{F} \Vert_{L^2(\Omega)} \Vert \bar{F} \Vert_{L^\infty(Q_T)} \\ 
		&+ \delta \Vert w \Vert_{L^\infty(Q_T)}\Vert \partial_t \bar{F} \Vert_{L^2(\Omega)}\Vert \nabla \psi \Vert_{L^2(\Omega)} + \frac{\delta^2}{2\varepsilon} \Vert \bar{F} \Vert^2_{L^\infty(Q_T)} \int_\Omega |\nabla \psi + \delta w \bar{F}|^2 \d x. 
		\end{aligned}
\end{equation*}
Integrating in time and using the H\"{o}lder inequality in the second term on the right-hand side, we get, for every $t\in[0,T]$, 
\begin{equation*}
	\begin{aligned}
		\frac{1}{2}\int_\Omega |\nabla & \psi + \delta w\bar{F}|^2 \d x + \frac{\varepsilon}{2} \int_0^t \int_{\Omega} (\Delta \psi + \delta \dv(w \bar{F}) )^2 \d x \d t \\ 
		\leq & \frac{1}{2}\int_\Omega |\nabla \psi(0,x)+\delta w_0(x) \bar{F}(0,x)|^2 \d x + \delta^2 |\Omega|^{1/2} \Vert w \Vert_{L^\infty({Q}_T)}^2 \Vert \partial_t \bar{F} \Vert_{L^1(0,T;L^2(\Omega))} \Vert \bar{F} \Vert_{L^\infty(Q_T)} \\ 
		&  + \delta \Vert w \Vert_{L^\infty(Q_T)}\Vert \partial_t \bar{F} \Vert_{L^1(0,T;L^2(\Omega))} \sup_{\tau \in [0,t]}\Vert \nabla \psi(\tau,\cdot) \Vert_{L^2(\Omega)} + \frac{\delta^2}{2\varepsilon} \Vert \bar{F} \Vert^2_{L^\infty(Q_T)} \int_0^t \int_\Omega |\nabla \psi + \delta w \bar{F}|^2 \d x. 
		\end{aligned}
\end{equation*}
By the triangle inequality, $\Vert \nabla \psi(\tau,\cdot) \Vert_{L^2(\Omega)} \leq \Vert \nabla \psi(\tau,\cdot) + \delta w(\tau,\cdot) \bar{F}(\tau,\cdot) \Vert_{L^2(\Omega)} + \Vert \delta w \bar{F}(\tau,\cdot) \Vert_{L^2(\Omega)}$, and using the Young inequality as well, we get, for every $t\in[0,T]$, 
\begin{equation*}
	\begin{aligned}
		\frac{1}{2}\int_\Omega |\nabla \psi &+ \delta w\bar{F}|^2 \d x + \frac{\varepsilon}{2} \int_0^t \int_{\Omega} (\Delta \psi + \delta \dv(w \bar{F}) )^2 \d x \d t \leq \frac{1}{2}\int_\Omega |\nabla \psi(0,x)+\delta w_0(x) \bar{F}(0,x)|^2 \d x \\ 
		&+ \delta^2 |\Omega|^{1/2} \Vert w \Vert_{L^\infty({Q}_T)}^2 \Vert \partial_t \bar{F} \Vert_{L^1(0,T;L^2(\Omega))} \Vert \bar{F} \Vert_{L^\infty(Q_T)} + \frac{1}{4}\sup_{\tau \in [0,t]}\Vert \nabla \psi(\tau,\cdot) + \delta w(\tau,\cdot) \bar{F}(\tau,\cdot) \Vert_{L^2(\Omega)}^2 \\ 
		&  + \delta^2 \Vert w \Vert_{L^\infty(Q_T)}^2\Vert \partial_t \bar{F} \Vert_{L^1(0,T;L^2(\Omega))}^2 + \delta^2|\Omega|T \Vert w \Vert_{L^\infty(Q_T)}^2\Vert \partial_t \bar{F} \Vert_{L^1(0,T;L^2(\Omega))} \Vert \bar{F} \Vert_{L^\infty(Q_T)} \\ 
		&+ \frac{\delta^2}{2\varepsilon} \Vert \bar{F} \Vert^2_{L^\infty(Q_T)} \int_0^t \int_\Omega |\nabla \psi + \delta w \bar{F}|^2 \d x. 
		\end{aligned}
\end{equation*}
As a result, we write 
\begin{equation*}
	\begin{aligned}
		\frac{1}{4}\sup_{\tau\in[0,t]}\Vert & \nabla \psi(\tau,\cdot) + \delta w(\tau,\cdot)\bar{F}(\tau,\cdot)\Vert_{L^2(\Omega)}^2 + \frac{\varepsilon}{2} \int_0^t \int_{\Omega} (\Delta \psi + \delta \dv(w \bar{F}) )^2 \d x \d t \leq \\ 
		&\frac{1}{2}\int_\Omega |\nabla \psi(0,x)+\delta w_0(x) \bar{F}(0,x)|^2 \d x + \delta^2 |\Omega|^{1/2} \Vert w \Vert_{L^\infty({Q}_T)}^2 \Vert \partial_t \bar{F} \Vert_{L^1(0,T;L^2(\Omega))} \Vert \bar{F} \Vert_{L^\infty(Q_T)} \\ 
		& + \delta^2|\Omega|T \Vert w \Vert_{L^\infty(Q_T)}^2\Vert \partial_t \bar{F} \Vert_{L^1(0,T;L^2(\Omega))} \Vert \bar{F} \Vert_{L^\infty(Q_T)} + \frac{\delta^2}{2} \Vert w \Vert_{L^\infty(Q_T)}^2\Vert \partial_t \bar{F} \Vert_{L^1(0,T;L^2(\Omega))}^2 \\ 
		&   + \frac{\delta^2}{2\varepsilon} \Vert \bar{F} \Vert^2_{L^\infty(Q_T)} \int_0^t \sup_{y\in[0,\tau]}\Vert \nabla \psi(y,\cdot) + \delta w(y,\cdot) \bar{F}(y,\cdot) \Vert_{L^2(\Omega)}^2 \d \tau. 
		\end{aligned}
\end{equation*}
Using a combination of \eqref{eq:w Linfty bound}, the formula \eqref{eq:new psi def}, and \eqref{eq:L2H1 bound that you already know}, we deduce from the previous inequality that, for every $t \in [0,T]$, 
\begin{equation}\label{eq:unpleasant yields its fruit}
	\begin{aligned}
		 \sup_{\tau\in[0,t]}\Vert \nabla \psi(\tau,\cdot) + \delta w(\tau,\cdot)&\bar{F}(\tau,\cdot) \Vert_{L^2(\Omega)}^2 + \varepsilon \int_0^t \int_{\Omega} (\Delta \psi + \delta \dv(w \bar{F}) )^2 \d x \d t \\ 
		& \leq C(\varepsilon,\delta,T,\Omega,\Vert w_0 \Vert_{L^\infty(\Omega)}, \Vert \nabla w_0 \Vert_{L^\infty(\Omega)},\Vert \bar{F}\Vert_{L^\infty(Q_T)},\Vert \partial_t \bar{F}\Vert_{L^1(0,T;L^2(\Omega))}) \cdot \\ 
		&\qquad \bigg( 1 + \int_0^t \sup_{y\in[0,\tau]}\Vert \nabla \psi(y,\cdot) + \delta w(y,\cdot)\bar{F}(y,\cdot)\Vert_{L^2(\Omega)}^2 \d \tau \bigg), 
		\end{aligned}
\end{equation}
where $C$ on the right-hand side denotes a quantity depending only the parameters inside its brackets. Then, an application of Gr\"{o}nwall's Lemma yields 
\begin{equation*}
	\begin{aligned}
		\sup_{t\in[0,T]}\Vert \nabla \psi(t,\cdot) &+ \delta w(t,\cdot)\bar{F}(t,\cdot)\Vert_{L^2(\Omega)}^2 \leq C(\varepsilon,\delta,T,\Omega,\Vert w_0 \Vert_{C^1(\bar{\Omega})},\Vert \bar{F}\Vert_{L^\infty(Q_T)},\Vert \partial_t \bar{F}\Vert_{L^1(0,T;L^2(\Omega))}). 
		\end{aligned}
\end{equation*}
Hence, returning to \eqref{eq:unpleasant yields its fruit} and using \eqref{eq:w Linfty bound}, we have 
\begin{equation}\label{eq:where you must return when adding nonlocal drift}
	\begin{aligned}
	\varepsilon \int_{Q_T} (\Delta \psi + \delta & \dv(w\bar{F}))^2 \d x \d t \leq C(\varepsilon,\delta,T,\Omega,\Vert w_0 \Vert_{C^1(\bar{\Omega})},\Vert \bar{F}\Vert_{L^\infty(Q_T)},\Vert \partial_t \bar{F}\Vert_{L^1(0,T;L^2(\Omega))}). 
	\end{aligned}
\end{equation}
Using $\int_{Q_T} (\dv(w\bar{F}))^2 \d x \d t \leq C(\Vert \bar{F}\Vert_{L^\infty(Q_T)} , \Vert \nabla \bar{F}\Vert_{L^\infty(Q_T)})\Vert w \Vert^2_{L^2(0,T;H^1(\Omega)}$ along with the estimate \eqref{eq:L2H1 bound that you already know}, and using the triangle inequality, we obtain 
\begin{equation*}
	\begin{aligned}
	&\int_{Q_T} (\Delta \psi )^2 \d x \d t \leq  C(\varepsilon,\delta,T,\Omega,\Vert w_0 \Vert_{C^1(\bar{\Omega})},\Vert \bar{F}\Vert_{L^\infty(Q_T)},\Vert \partial_t \bar{F}\Vert_{L^1(0,T;L^2(\Omega))},\Vert \nabla \bar{F}\Vert_{L^\infty(Q_T)}). 
	\end{aligned}
\end{equation*}
Note that $\Delta \psi = (\varepsilon+w)\Delta w + |\nabla w|^2$ and that $w$ is non-negative, whence the previous estimate and \eqref{eq:L2H1 bound that you already know} yield the desired estimate \eqref{eq:H2bound quantitative}.

\textbf{Step II: } We emphasise that the addition of drift terms changes nothing to this argument---and one would still define $\psi$ as per \eqref{eq:new psi def}---since, by \eqref{eq:Li terms}, they can be bounded in $L^\infty$ as follows 
\begin{equation*}
	\begin{aligned}
		\Vert \nabla L(\cdot,\cdot,w) \Vert_{L^\infty(Q_T)} &\leq \Vert V \Vert_{C^1(\mathbb{R}^{d+1})} + \Vert \nabla W*w(t,\cdot) \Vert_{L^\infty(Q_T)} \leq \Vert V \Vert_{C^1(\mathbb{R}^{d+1})} + |\Omega| \Vert W \Vert_{C^1(\mathbb{R}^{d})} \Vert w \Vert_{L^\infty(Q_T)}, 
	\end{aligned}
\end{equation*}
in conjunction with \eqref{eq:w Linfty bound}, where we omitted the $i$ subscripts. Similarly, any additional space derivatives fall directly on $V$ and $W$, not on $w$, and so 
\begin{equation*}
	\begin{aligned}
		\Vert \nabla^2 L(\cdot,\cdot,w) \Vert_{L^\infty(Q_T)} &\leq \Vert V \Vert_{C^2(\mathbb{R}^d)} + |\Omega| \Vert W \Vert_{C^2(\mathbb{R}^d)} \Vert w \Vert_{L^\infty(Q_T)}. 
	\end{aligned}
\end{equation*}
For an additional time derivative, we have 
\begin{equation*}
	\begin{aligned}
		\Vert \partial_t \nabla L(\cdot,\cdot,w) \Vert_{L^1(0,T;L^2(\Omega))} &\leq |\Omega|^{1/2}T\Vert V \Vert_{C^2(\mathbb{R}^{d+1})} + \Vert \nabla W*\partial_t w(t,\cdot) \Vert_{L^1(0,T;L^2(\Omega))}, 
	\end{aligned}
\end{equation*}
and the second term on the right-hand side of the above is controlled---using \eqref{eq:formula dt w}, \eqref{eq:new evolution}, the Cauchy--Schwartz integral inequality, and the Young inequality---as 
\begin{equation*}
	\begin{aligned}
		\int_0^T \Vert \nabla W*\partial_t w(t,\cdot) \Vert_{L^2(\Omega)} \d t &\leq \Vert W \Vert_{C^1(\mathbb{R}^{d})} T^{1/2} \Vert \partial_t w \Vert_{L^2(Q_T)} \\ 
		&\leq \frac{\varepsilon}{4}\int_{Q_T} (\Delta \psi + \delta \dv(w \bar{F}))^2 \d x \d t + \frac{1}{\varepsilon}\Vert W \Vert_{C^1(\mathbb{R}^{d})}^2 T, 
		\end{aligned}
\end{equation*}
and the first term on the right-hand side can be absorbed into the left-hand side of \eqref{eq:where you must return when adding nonlocal drift}. Thus, we have shown that the term of the form $w \nabla L$ can be handled in exactly the same way as $\delta w \bar{F}$ --- notice that we never make use of the $\delta$ smallness assumption. 
\end{proof}

\subsection{Proofs of Lemmas \ref{lem:preserve sign and mass}, \ref{lem:preserve initial data}, and \ref{lem:lowersemi the one we use}: Quantities preserved under weak limit and lower semicontinuity}\label{sec:two technical seq lemmas}

\begin{proof}[Proof of Lemma \ref{lem:preserve sign and mass}]
Given any non-negative test function $\theta \in C^1(\bar{Q}_T)$, since from the first hypothesis $\zeta^n \geq 0$ a.e.~in $Q_T$, the weak convergence implies 
\begin{equation}\label{eq:nonneg weak cpctness}
	0 \leq \int_{Q_T} \theta \zeta^n \d x \d t \to \int_{Q_T} \theta \zeta \d x \d t, 
\end{equation}
from which it follows that $\zeta \geq 0$ a.e.~in $Q_T$. Similarly, given any $\eta \in C^1([0,T])$ and interpreting it as a function in $C^1(\bar{Q}_T)$, using the second hypothesis, there holds 
\begin{equation}\label{eq:mass consv weak cpctness}
	\int_0^T \eta(t) \Lambda \d t = \int_0^T \eta(t)\bigg( \int_\Omega \zeta^n \d x \bigg) \d t \to \int_0^T \eta(t)\bigg( \int_\Omega \zeta \d x \bigg) \d t, 
\end{equation}
by the weak convergence in $L^2(Q_T)$ and the Tonelli--Fubini theorem. Thus, combining the above with the non-negativity of the limit $\zeta$, 
\begin{equation}\label{eq:LinftyL1 weak cpctness}
	\int_\Omega |\zeta(t,x)| \d x = \int_\Omega \zeta(t,x) \d x = \Lambda \qquad \text{a.e.~}t \in (0,T), 
\end{equation}
which fulfils the requirement for $\zeta \in L^\infty(0,T;L^1(\Omega))$. 
\end{proof}

\begin{proof}[Proof of Lemma \ref{lem:preserve initial data}]
Following the computation in Remark \ref{rem:embed into cts negative sobolev space}, given $\varphi \in L^r(0,T;W^{1,r}(\Omega)) = X$ with $\Vert \varphi \Vert_X \leq 1$, using the H\"{o}lder inequality, 
\begin{equation*}
    \bigg|\int_{Q_T} \zeta^n \varphi \d x \d t \bigg| \leq \Vert \zeta^n \Vert_{L^2(Q_T)} (|\Omega|T)^{\frac{d}{2(d+1)}}. 
\end{equation*}
Thus, $\Vert \zeta^n \Vert_{X'} \leq \Vert \zeta^n \Vert_{L^2(Q_T)} (|\Omega|T)^{\frac{d}{2(d+1)}}$, where $X' = L^{r'}(0,T;(W^{1,r}(\Omega))')$. Meanwhile, we also have $\partial_t \zeta^n \in X'$ for each $n\in\mathbb{N}$, and it therefore follows that $\{\zeta^n\}_{n\in\mathbb{N}}$ is a sequence in $W^{1,r'}(0,T;(W^{1,r}(\Omega))')$; in fact, it is uniformly bounded therein. By \cite[Theorem 2 of Section 5.9.2]{evans1998partial}, it follows that $\zeta^n \in C([0,T];(W^{1,r}(\Omega))')$ for each $n\in\mathbb{N}$, so that the evaluation of $\zeta^n$ at time zero is well-defined in assumption 2 of the lemma. By the same reasoning, $\zeta \in C([0,T];(W^{1,r}(\Omega))')$. Moreover, \cite[Theorem 2 of Section 5.9.2]{evans1998partial} implies that, for each $n\in\mathbb{N}$,
\begin{equation*}
    \zeta^n(t,\cdot) = \zeta^n(s,\cdot) + \int_s^t \partial_t \zeta^n(\tau,\cdot) \d \tau \qquad \forall 0 \leq s \leq t \leq T, 
\end{equation*}
where the equality holds in the sense of $(W^{1,r}(\Omega))'$. That is, given any $\phi \in W^{1,r}(\Omega)$, there holds, for all $0 \leq s \leq t \leq T$, 
\begin{equation}\label{eq:invoke evans time deriv}
    \int_\Omega \zeta^n(t,x) \phi(x) \d x = \int_\Omega \zeta^n(s,x) \phi(x) \d x + \int_s^t \langle \partial_t \zeta^n(\tau,\cdot) , \phi \rangle_\Omega \d \tau, 
\end{equation}
where $\langle \cdot , \cdot \rangle_\Omega$ denotes the duality product between $W^{1,r}(\Omega)$ and its dual, and this duality product may be interchanged with the Bochner integral by virtue of \cite[Theorem 8 of Appendix E.5]{evans1998partial} and the summability 
\begin{equation*}
    \int_s^t \Vert \partial_t \zeta^n(\tau,\cdot) \Vert_{(W^{1,r}(\Omega))'} \d \tau \leq \int_0^T \Vert \partial_t \zeta^n(\tau,\cdot) \Vert_{(W^{1,r}(\Omega))'} \d \tau \leq T^{\frac{1}{r}} \Vert \partial_t \zeta^n \Vert_{X'}, 
\end{equation*}
where we used the H\"{o}lder inequality to obtain the final inequality. We now use the same bounding strategy as in \eqref{eq:invoke evans time deriv}. We obtain that, given any any $\phi \in W^{1,r}(\Omega)$, there holds, for all $0 \leq s \leq t \leq T$, 
\begin{equation*}
    \bigg|\int_\Omega \zeta^n(t,x) \phi(x) \d x - \int_\Omega \zeta^n(s,x) \phi(x) \d x \bigg| \leq (t-s)^{\frac{1}{r}} \Vert \partial_t \zeta^n \Vert_{X'} \Vert \phi \Vert_{W^{1,r}(\Omega)}. 
\end{equation*}
Since $\{\partial_t \zeta^n\}_{n\in\mathbb{N}}$ is converging weakly-* in $X'$, it follows that  $ \Vert \partial_t \zeta^n \Vert_{X'} \leq \Lambda$ for some $\Lambda$ (though this need not be $ \Vert \partial_t \zeta \Vert_{X'}$), independent of $n$. Thus, 
\begin{equation}\label{eq:to show cty for zeta pre}
    \bigg|\int_\Omega \zeta^n(t,x) \phi(x) \d x - \int_\Omega \zeta^n(s,x) \phi(x) \d x \bigg| \leq (t-s)^{\frac{1}{r}} \Lambda \Vert \phi \Vert_{W^{1,r}(\Omega)} \quad \forall \phi \in W^{1,r}(\Omega), 
\end{equation}
for all $0 \leq s \leq t \leq T$. Passing to the limit weakly as $n\to\infty$ in the above, we obtain \eqref{eq:to show cty for zeta}. Similarly, setting $s=0$ in \eqref{eq:to show cty for zeta pre}, using assumption 2 of the lemma, and then passing to the limit weakly as $n\to\infty$, we obtain 
\begin{equation*}
    \bigg|\int_\Omega \zeta(t,x) \phi(x) \d x - \int_\Omega \zeta_0(x) \phi(x) \d x \bigg| \leq t^{\frac{1}{r}} \Lambda \Vert \phi \Vert_{W^{1,r}(\Omega)} \quad \forall 0 < t \leq T, 
\end{equation*}
respectively. Taking the supremum over all $\phi \in W^{1,r}(\Omega)$ with unit norm then yields \eqref{eq:to show cty for zeta with initial}, which concludes the proof. 
\end{proof}

\begin{proof}[Proof of Lemma \ref{lem:lowersemi the one we use}]

\textbf{Step I: } To begin with, observe that since $\zeta^n \to \zeta$ strongly in $L^1(Q_T)$, there exists a subsequence (still indexed by $n$) such that 
\begin{equation*}
    \sum_{n\in\mathbb{N}} \Vert \zeta^n - \zeta \Vert_{L^1(Q_T)} < +\infty. 
\end{equation*}
From here until the end of this proof, we only consider this particular subsequence, and we never pass to further subsequences. In what follows we use the identification $L^1(Q_T) = L^1(0,T;L^1(\Omega))$ and write $\boldsymbol{\zeta}^n(t) = \zeta^n(t,\cdot)$ and $\boldsymbol{\zeta}(t) = \zeta(t,\cdot)$. We consider, for a.e.~$t\in(0,T)$, the convergence of $\boldsymbol{\zeta}^n(t)$ towards $\boldsymbol{\zeta}(t)$ in the norm of $L^1(\Omega)$. In particular, the Minkowski inequality for infinite sums in $L^1((0,T))$ yields 
\begin{equation*}
    \bigg\Vert \sum_{n\in\mathbb{N}} \Vert \boldsymbol{\zeta}^n(\cdot) - \boldsymbol{\zeta}(\cdot) \Vert_{L^1(\Omega)}  \bigg\Vert_{L^1((0,T))} \leq  \sum_{n\in\mathbb{N}} \Vert {\zeta}^n - {\zeta} \Vert_{L^1(Q_T)} < +\infty, 
\end{equation*}
and hence the infinite series $\sum_{n\in\mathbb{N}} \Vert \boldsymbol{\zeta}^n(\cdot) - \boldsymbol{\zeta}(\cdot) \Vert_{L^1(\Omega)}$ is well-defined as an element of $L^1((0,T))$. Since any integrable function is finite almost everywhere, it follows that, for a.e.~$t\in(0,T)$, 
\begin{equation*}
    \sum_{n\in\mathbb{N}} \Vert \boldsymbol{\zeta}^n(t) - \boldsymbol{\zeta}(t) \Vert_{L^1(\Omega)} < +\infty, 
\end{equation*}
and it immediately follows from the summability of this series that, for a.e.~$t\in(0,T)$, we have $\Vert \boldsymbol{\zeta}^n(t) - \boldsymbol{\zeta}(t) \Vert_{L^1(\Omega)} \to 0$ as $n\to\infty$. Rewriting in terms of the original functions, we have shown that, for a.e.~$t\in(0,T)$, 
\begin{equation*}
    \lim_{n\to\infty} \Vert \zeta^n(t,\cdot) - \zeta(t,\cdot) \Vert_{L^1(\Omega)} = 0. 
\end{equation*}

\textbf{Step II: } Note that $\zeta$ is also non-negative from the assumption that $\zeta^n \to \zeta$ in $L^1(Q_T)$. The convergence obtained at the end of Step I is sufficient to satisfy the hypothesis of Lemma \ref{lem:lowersemi}. An application of this latter result gives, for a.e.~$t \in (0,T)$, 
\begin{equation*}
	\int_\Omega \frac{|\nabla \zeta(t,x)|^2}{\zeta(t,x)} \d x \leq \liminf_{n\to\infty} \int_\Omega \frac{|\nabla \zeta^{n}(t,x)|^2}{\zeta^{n}(t,x)} \d x. 
\end{equation*}
By integrating the above inequality with respect to time and then applying the Fatou Lemma, we obtain 
\begin{equation*}
	\begin{aligned}
		\int_{Q_T} \frac{|\nabla \zeta|^2}{\zeta} \d x \d t &\leq \int_0^T \liminf_{n\to\infty} \bigg( \int_\Omega \frac{|\nabla \zeta^{n}(t,x)|^2}{\zeta^{n}(t,x)} \d x \bigg) \d t \leq \liminf_{n\to\infty} \int_{Q_T} \frac{|\nabla \zeta^{n}|^2}{\zeta^{n}} \d x \d t, 
		\end{aligned}
\end{equation*}
as required. 
\end{proof}

\subsection{Smoothing operator}\label{appendix:smoothing op}

The purpose of this appendix is to verify the properties of the smoothing operator $R_\mu$ in \eqref{eq:Rnu introduce} stated in Section \ref{section:weak compactness}. To this end, we make the following definitions, and we note that, throughout this section only, $\phi$ refers to the Friedrichs bump function (see below) and not to a generic test function.  

\begin{defi}
Define $\phi \in C^\infty_c(\mathbb{R})$ to be the standard Friedrichs bump function, i.e., 
\begin{equation*}
    \phi(y) = \left\lbrace \begin{aligned}
        & \exp\left(- \frac{1}{1-|y|^2}\right) \qquad && \text{for } |y| < 1, \\ 
        &0 \qquad && \text{for } |y| \geq 1, 
    \end{aligned} \right.
\end{equation*}
and define, for $\mu>0$, 
\begin{equation}
\phi_\mu(t,x) := c_\phi \mu^{-(d+1)}  \phi\left(\frac{\sqrt{t^2+|x|^2}}{\mu}\right) \qquad \forall (t,x) \in \mathbb{R}^{d+1}, 
\end{equation}
where the positive constant $c_\phi$ is chosen such that $\int_{\mathbb{R}^{d+1}} \phi_1(y) \d y =1$. 
\end{defi}
Note that $\Vert \phi \Vert_{L^\infty(\mathbb{R})} = 1$ and 
\begin{equation*}
    \supp \phi_\mu = \bar{B}(0,\mu), 
\end{equation*}
where the latter is the closed ball of radius $\mu$, and  $\int_{\mathbb{R}^{d+1}} \phi_\mu(y) \d y = 1$ for every $\mu>0$. 

\begin{defi}
We define the operator $A_\mu : L^2(\mathbb{R};H^1(\mathbb{R}^d)) \to C^\infty(\bar{Q}_T)$ to be the restriction of the mollification to the parabolic cylinder, i.e., 
\begin{equation*}
    A_\mu u := \int_\mathbb{R} \int_{\mathbb{R}^d} \phi_\mu(\cdot-\tau,\cdot-y) u(\tau,y) \d y \d \tau \bigg|_{\bar{Q}_T}. 
\end{equation*}
\end{defi}

\begin{lem}[Preliminary spatial extension]\label{lem:sob extension prelim}
There exists a bounded linear operator $E' : L^2(0,T;H^1(\Omega)) \to L^2(0,T;H^1(\mathbb{R}^d))$ such that 
\begin{equation*}
    E'f(t,x) = f(t,x) \qquad \text{a.e. } (t,x) \in Q_T, 
\end{equation*}
and, with $\Omega_1 := \{ x \in \mathbb{R}^d : d(x,\Omega) \leq 1 \}$, in which $\Omega$ is compactly contained, 
\begin{equation*}
    \supp E'f(t,\cdot) \subset \Omega_1 \qquad \text{a.e. } t \in (0,T). 
\end{equation*}
Moreover, there exists a constant $C>0$, depending only on $\Omega$, such that, given any $f \in L^2(0,T;H^1(\Omega))$, 
\begin{equation}\label{eq:Eprime f just L2}
    \Vert E' f \Vert_{L^2((0,T)\times \mathbb{R}^d)} \leq C \Vert f \Vert_{L^2(Q_T)}. 
\end{equation}
\end{lem}
\begin{proof}
By definition of $L^2(0,T;H^1(\Omega)) \ni f$, we have that for a.e.~$t \in (0,T)$, the element $f(t,\cdot)$ belongs to $H^1(\Omega)$. The result follows at once from repeating the argument of the proof Theorem 1 of \cite[Section 5.4]{evans1998partial} on the function $f(t,\cdot)$, with the time coordinate kept fixed. 
\end{proof}

\begin{lem}[Sobolev extension for spaces involving time]\label{lem:sob extension time}
There exists a bounded linear operator $E : L^2(0,T;H^1(\Omega)) \to L^2(\mathbb{R};H^1(\mathbb{R}^d))$ such that 
\begin{equation*}
    Ef(t,x) = f(t,x) \qquad \text{a.e. } (t,x) \in Q_T, 
\end{equation*}
and, with $\Omega_1 := \{ x \in \mathbb{R}^d : d(x,\Omega) \leq 1 \}$, in which $\Omega$ is compactly contained, 
\begin{equation*}
    \supp Ef \subset [0,T] \times \Omega_1. 
\end{equation*}
Moreover, there exists a constant $C>0$, depending only on $\Omega$, such that, given any $u \in L^2(0,T;H^1(\Omega))$, 
\begin{equation*}
    \Vert Eu \Vert_{L^2(\mathbb{R}^{d+1})} \leq C \Vert u \Vert_{L^2(Q_T)}. 
\end{equation*}
\end{lem}
\begin{proof}
Simply define $E$ via the explicit formula 
\begin{equation*}
    Ef(t,x) := E'f(t,x) \mathds{1}(t)_{[0,T]} \qquad \text{for a.e.~} (t,x) \in \mathbb{R}^{d+1}, 
\end{equation*}
where $E'$ is the bounded linear operator of Lemma \ref{lem:sob extension prelim}. Indeed, then $E$ is manifestly linear and satisfies the equality 
\begin{equation*}
    Ef(t,x) = f(t,x) \qquad \text{a.e. } (t,x) \in Q_T, 
\end{equation*} 
and 
\begin{equation*}
    \supp Ef \subset [0,T] \times \Omega_1. 
\end{equation*}
Moreover, given any $u \in L^2(0,T;H^1(\Omega))$, we have 
\begin{equation*}
    \begin{aligned}
        \Vert E u \Vert_{L^2(\mathbb{R};H^1(\mathbb{R}^d))}^2 = \int_\mathbb{R} \Vert E u(t,\cdot) \Vert^2_{H^1(\mathbb{R}^d)} \d t = \int_0^T \Vert E' u(t,\cdot) \Vert^2_{H^1(\mathbb{R}^d)} \d t = \Vert E' u \Vert^2_{L^2(0,T;H^1(\mathbb{R}^d))}, 
    \end{aligned}
\end{equation*}
and so the first result follows directly from the fact that $E' : L^2(0,T;H^1(\Omega)) \to L^2(0,T;H^1(\mathbb{R}^d))$ is a bounded linear operator (\textit{cf.}~Lemma \ref{lem:sob extension prelim}). Similarly, 
\begin{equation*}
    \begin{aligned}
        \Vert E u \Vert_{L^2(\mathbb{R}^{d+1}))}^2 = \int_\mathbb{R} \Vert E u(t,\cdot) \Vert^2_{L^2(\mathbb{R}^d)} \d t = \int_0^T \Vert E' u(t,\cdot) \Vert^2_{L^2(\mathbb{R}^d)} \d t &= \Vert E' u \Vert^2_{L^2((0,T)\times\mathbb{R}^d)} \\ 
        &\leq C \Vert u \Vert^2_{L^2(Q_T)}, 
    \end{aligned}
\end{equation*}
where we used \eqref{eq:Eprime f just L2} of Lemma \ref{lem:sob extension prelim} for the final inequality. 
\end{proof}

\begin{lem}[Smoothing operator]\label{lem:smoothing operator}
Fix $\mu>0$. The smoothing operator $R_{\mu}$, defined explicitly by 
\begin{equation} \label{smoothing operator def}
    R_{\mu} := A_\mu \circ E : L^2(0,T;H^1(\Omega)) \to C^\infty(\bar{Q}_T), 
\end{equation}
admits, for some positive constant $C_{reg}$ independent of $\mu$, the estimate 
\begin{equation}\label{eq:L2H1 unif bound Rnu smoother}
    \Vert R_\mu u \Vert_{L^2(0,T;H^1(\Omega))} \leq C_{reg} \Vert u \Vert_{L^2(0,T;H^1(\Omega))}    \qquad \forall u \in L^2(0,T;H^1(\Omega)). 
\end{equation}
As such, it is a bounded linear operator from $L^2(0,T;H^1(\Omega))$ to $C^\infty(\bar{Q}_T)$ equipped with the subspace norm-topology of $L^2(0,T;H^1(\Omega))$. Moreover, given any $u \in L^2(0,T;H^1(\Omega))$, we have the strong convergence 
\begin{equation}\label{eq:strong L2H1 conv smoother}
    \Vert R_\mu u - u \Vert_{L^2(0,T;H^1(\Omega))} \to 0 \qquad \text{as } \mu \to 0. 
\end{equation}
Also, for some positive constant $C_\mu$ depending on $\mu,\phi,\Omega,T$, 
\begin{equation}\label{eq:L2H2 nonunif bound Rnu smoother}
    \Vert R_\mu u \Vert_{L^\infty(0,T;W^{2,\infty}(\Omega))} \leq C_\mu \Vert u \Vert_{L^2(0,T;H^1(\Omega))}   \qquad \forall u \in L^2(0,T;H^1(\Omega)). 
\end{equation}
\end{lem}

\begin{proof}
Fix any $u \in L^2(0,T;H^1(\Omega))$. By Lemma \ref{lem:sob extension time}, we have that, for every $(t,x) \in \bar{Q}_T$ and every $\mu>0$, the integral 
\begin{equation}\label{eq:no derivs Amu}
    A_\mu (Eu)(t,x) = \int_\mathbb{R} \bigg( \int_{\mathbb{R}^d} \phi_\mu(t-\tau, x-y) Eu(\tau,y) \d y \bigg) \d \tau, 
\end{equation}
is well-defined, since the integrand 
\begin{equation*}
    (\tau,y) \mapsto \phi_\mu(t-\tau, x-y) Eu(\tau,y)
\end{equation*}
is compactly supported for every choice of $(t,x) \in Q_T$. Similarly, given any $m\in\mathbb{N}$, we have that 
\begin{equation*}
    \partial^{m}_x \phi_\mu(t-\tau,x-y) Eu(\tau,y) \qquad \text{and} \qquad \partial^{m}_t \phi_\mu(t-\tau,x-y) Eu(\tau,y)
\end{equation*}
are integrable over $(0,T)\times\mathbb{R}^d$ in view of the compact support. Moreover, using the shorthand $\bar{B}_\mu(t,x)$ to mean $\bar{B}((t,x),\mu)$, we have the uniform bounds 
\begin{equation*}
    |\partial^{m}_x \phi_\mu(t-\tau,x-y) Eu(\tau,y)| \leq C_{m,\mu} \mathds{1}_{\bar{B}_\mu(t,x)}(\tau,y) |Eu(\tau,y)|, 
    \end{equation*}
    and
    \begin{equation*}
    |\partial^{m}_t \phi_\mu(t-\tau,x-y) Eu(\tau,y)| \leq \tilde{C}_{m,\mu} \mathds{1}_{\bar{B}_\mu(t,x)}(\tau,y) |Eu(\tau,y)|, 
\end{equation*}
for positive constants $C_{m,\mu},\tilde{C}_{m,\mu}$ depending only on the supremum norms of the derivatives of $\phi_\mu$. By routine arguments involving the Dominated Convergence Theorem we justify differentiating under the integral in \eqref{eq:no derivs Amu}, and obtain the formulas 
\begin{equation}\label{eq:time deriv Amu E}
    \partial_t^{(m)} A_\mu (Eu)(t,x) = \int_\mathbb{R} \int_{\mathbb{R}^d} \partial_t^{(m)} \phi_\mu(t-\tau, x-y) Eu(\tau,y) \d y \d \tau \qquad \forall (t,x) \in \bar{Q}_T, 
\end{equation}
for any $m\in\mathbb{N}\cup\{0\}$, and 
\begin{equation}\label{eq:more space derivs Amu}
    \partial_{x_j}^{(m)} A_\mu (E u)(t,x) = \int_\mathbb{R} \int_{\mathbb{R}^d} \partial_{x_j}^{(m)} \phi_\mu(t-\tau, x-y) Eu(\tau,y) \d y \d \tau \quad \forall (t,x) \in \bar{Q}_T, \text{ for } j \in \{1,\dots,M\}. 
\end{equation}
An identical calculation yields the analogous formula for any mixed derivatives. Thus, it follows that $A_\mu$ maps elements of $L^2(\mathbb{R};H^1(\mathbb{R}^d))$ to smooth functions in $\bar{Q}_T$, as stated. 

In fact, in the case of one space derivative, we have that, for any $v \in L^2(\mathbb{R};H^1(\mathbb{R}^d))$, 
\begin{equation}\label{eq:in case of one space deriv}
    \partial_{x_j} A_\mu v(t,x) = \int_\mathbb{R} \int_{\mathbb{R}^d} \phi_\mu(t-\tau, x -y) \partial_{x_j} v(\tau,y) \d \tau \d y \qquad \forall (t,x) \in \bar{Q}_T, \text{ for } j \in \{1,\dots,M\}, 
\end{equation}
where we directly used the definition of weak derivative for $\nabla v \in L^2(\mathbb{R}\times\mathbb{R}^d)$, and there is no boundary term since we integrate over all of $\mathbb{R}^d$. Moreover, $A_\mu$ is linear, and for any $v \in L^2(\mathbb{R};H^1(\mathbb{R}^d))$, 
\begin{equation*}
    \Vert A_\mu v \Vert_{L^2(0,T;H^1(\Omega))}^2 = \int_0^T \bigg( \int_\Omega |A_\mu v(t,x) |^2 + |\nabla A_\mu v(t,x)|^2 \d x \bigg) \d t. 
\end{equation*}
Observe that 
\begin{equation*}
    \begin{aligned}
        \int_0^T \int_\Omega |A_\mu v (t,x)|^2 \d x \d t &\leq \int_{Q_T} \bigg( \int_\mathbb{R} \int_{\mathbb{R}^d} \phi_\mu(t-\tau,x-y) |v(\tau,y)| \d y \d \tau \bigg)^2 \d x \d t \\ 
        &= \int_{Q_T} \bigg( \int_\mathbb{R} \int_{\mathbb{R}^d} \sqrt{\phi_\mu}(t-\tau,x-y) \sqrt{\phi_\mu}(t-\tau,x-y)|v(\tau,y)| \d y \d \tau \bigg)^2 \d x \d t, 
    \end{aligned}
\end{equation*}
so that an application of the Cauchy--Schwarz integral inequality yields 
\begin{equation*}
    \begin{aligned}
       \int_0^T \int_\Omega |A_\mu v (t,x)|^2 \d x \d t  &\leq \int_{Q_T} \bigg( \int_{\mathbb{R}^{d+1}} \phi_\mu(t-\tau,x-y) \d y \d \tau \bigg) \bigg( \int_{ \mathbb{R}^{d+1}} \phi_\mu(t-\tau,x-y) |v(\tau,y)|^2 \d y \d \tau \bigg) \d x \d t \\ 
       &\leq \int_{Q_T} \bigg( \int_{\mathbb{R}^{d+1}} \phi_\mu(t-\tau,x-y) |v(\tau,y)|^2 \d y \d \tau \bigg) \d x \d t, 
    \end{aligned}
\end{equation*}
where we used the translation invariance of Lebesgue measure and $\int_{\mathbb{R}^{d+1}} \phi_\mu(t,x) \d x \d t =1$. Then, an application of the Tonelli--Fubini theorem yields 
\begin{equation}\label{eq:for L2 cyl unif bound}
    \begin{aligned}
        \int_{Q_T} \bigg( \int_{ \mathbb{R}^{d+1}} \phi_\mu(t-\tau,x-y) |v(\tau,y)|^2 \d y \d \tau \bigg) \d x \d t &= \int_{\mathbb{R}^{d+1}} |v(\tau,y)|^2 \bigg( \int_{Q_T} \phi_\mu(t-\tau,x-y)  \d x \d t \bigg) \d y \d \tau \\ 
        &\leq \int_{\mathbb{R}^{d+1}} |v(\tau,y)|^2 \d y \d \tau = \Vert v \Vert^2_{L^2(\mathbb{R}^{d+1})}. 
    \end{aligned}
\end{equation}
Similarly, an identical argument yields 
\begin{equation*}
    \begin{aligned}
        \int_0^T \int_\Omega |\nabla A_\mu v(t,x)|^2 \d x \d t &\leq \int_{Q_T} \bigg( \int_\mathbb{R} \int_{\mathbb{R}^d} \phi_\mu(t-\tau,x-y) |\nabla v(\tau,y)| \d y \d \tau \bigg)^2 \d x \d t \leq \Vert \nabla v \Vert^2_{L^2(\mathbb{R}^{d+1})}, 
    \end{aligned}
\end{equation*}
so that 
\begin{equation*}
    \Vert A_\mu v \Vert_{L^2(0,T;H^1(\Omega))} \leq \Vert v \Vert_{L^2(\mathbb{R};H^1(\mathbb{R}^d))} \qquad \forall v \in L^2(\mathbb{R};H^1(\mathbb{R}^d)). 
\end{equation*}
It follows that $A_\mu$ is a bounded linear operator from $L^2(\mathbb{R};H^1(\mathbb{R}^d))$ to $ C^\infty(\bar{Q}_T)$ endowed with the subspace norm-topology of $L^2(0,T;H^1(\Omega))$. Hence, in view of Lemma \ref{lem:sob extension time}, the composition $R_\mu$ is a bounded linear operator from $L^2(0,T;H^1(\Omega))$ to $ C^\infty(\bar{Q}_T)$ endowed with the subspace norm-topology of $L^2(0,T;H^1(\Omega))$, and \eqref{eq:L2H1 unif bound Rnu smoother} follows. 

Returning to \eqref{eq:no derivs Amu}, a direct application of the H\"{o}lder inequality shows that 
\begin{equation*}
    |A_\mu (Eu)(t,x)| \leq \Vert \phi_\mu \Vert_{L^2(\mathbb{R}^{d+1})} \Vert Eu \Vert_{L^2(\mathbb{R}^{d+1})} \leq C_\mu \Vert u \Vert_{L^2(Q_T)}, 
\end{equation*}
where we used Lemma \ref{lem:sob extension time} to obtain the final inequality, for some constant $C_\mu$ that also depends on $\mu$. Note in passing that $\Vert \phi_\mu \Vert_{L^2(\mathbb{R}^{d+1})}$ is well-defined for every $\mu>0$ due to the smoothness and compact support of $\phi_\mu$. Thus, 
\begin{equation*}
    \Vert A_\mu (Eu) \Vert_{L^\infty(Q_T)} \leq C_\mu \Vert u \Vert_{L^2(Q_T)}, 
\end{equation*}
and an analogous computation using formula \eqref{eq:more space derivs Amu} yields the corresponding estimates for higher derivatives. We deduce that, for some positive constant $C_\mu$ depending on $\mu,\Omega,\phi,T$, the bound \eqref{eq:L2H2 nonunif bound Rnu smoother} holds.

Finally, we verify the strong convergence in $L^2(0,T;H^1(\Omega))$. Fix $u \in L^2(0,T;H^1(\Omega))$ and any Lebesgue point $(t,x) \in Q_T$. Since, by Lemma \ref{lem:sob extension time}, $Eu = u$ a.e.~on $Q_T$ it follows that $Eu(t,x) = u(t,x)$ at this Lebesgue point. Using this latter fact and $\int_{\mathbb{R}^{d+1}} \phi_\mu(t,x) \d x \d t = 1$, and since the Lebesgue point $(t,x)$ was chosen arbitrarily, it follows that there holds 
\begin{equation}\label{eq:Rnu u minus u ae on cylinder}
    R_\mu u(t,x) - u(t,x) = \int_\mathbb{R} \bigg( \int_{\mathbb{R}^d} \phi_\mu(t-\tau, x-y) \big( Eu(\tau,y) - Eu(t,x) \big) \d y \bigg) \d \tau \quad \text{a.e.~}(t,x) \in Q_T. 
\end{equation}
Similarly, using the definition of weak derivative in $\mathbb{R}^d$ to place the space derivatives on $Eu$ instead of $\phi_\mu$ without additional boundary term, 
\begin{equation}\label{eq:gradients Rnu u minus u ae on cylinder}
    \nabla(R_\mu u - u)(t,x) = \int_\mathbb{R} \bigg( \int_{\mathbb{R}^d} \phi_\mu(t-\tau, x-y) \big( \nabla Eu(\tau,y) - \nabla Eu(t,x) \big) \d y \bigg) \d \tau \quad \text{a.e.~} (t,x) \in Q_T. 
\end{equation}
We now estimate, with respect to the norm on $L^2(Q_T)$, the terms of \eqref{eq:Rnu u minus u ae on cylinder} and \eqref{eq:gradients Rnu u minus u ae on cylinder}. To begin with, 
\begin{equation*}
    \Vert R_\mu u - u \Vert^2_{L^2(Q_T)} = \int_{Q_T} \bigg| \int_\mathbb{R} \bigg( \int_{\mathbb{R}^d} \phi_\mu(t-\tau, x-y) \big( Eu(\tau,y) - Eu(t,x) \big) \d y \bigg) \d \tau \bigg|^2 \d x \d t, 
\end{equation*}
which, by writing $z=(t,x)$ and $\xi = (\tau,y)$, can be rewritten as 
\begin{equation*}
    \Vert R_\mu u - u \Vert^2_{L^2(Q_T)} = \int_{Q_T} \bigg| \int_{\mathbb{R}^{d+1}} \phi_\mu(z-\xi) \big( Eu(\xi) - Eu(z) \big) \d \xi \bigg|^2 \d z. 
\end{equation*}
An application of the Cauchy--Schwarz integral inequality then yields 
\begin{equation*}
    \begin{aligned}
        \Vert R_\mu u - u \Vert^2_{L^2(Q_T)} &\leq \int_{Q_T} \bigg( \int_{\mathbb{R}^{d+1}} \phi_\mu(z-\xi) \d \xi \bigg) \bigg(\int_{\mathbb{R}^{d+1}} \phi_\mu(z-\xi) |Eu(\xi) - Eu(z)|^2 \d \xi \bigg) \d z \\ 
        &= \int_{Q_T} \bigg(\int_{\mathbb{R}^{d+1}} \phi_\mu(\zeta) |Eu(z+\zeta) - Eu(z)|^2 \d \zeta \bigg) \d z \\ 
        &= c_\phi \int_{Q_T} \bigg(\int_{B(0,1)} \phi(|\xi|) |Eu(z+\mu\xi) - Eu(z)|^2 \d \xi \bigg) \d z, 
    \end{aligned}
\end{equation*}
 where we used the evenness of $\phi$ in the penultimate line, and took into account the scaling with respect to $\mu$ explicitly in the final line. Hence, applying the Tonelli--Fubini theorem, we get 
 \begin{equation*}
    \begin{aligned}
        \Vert R_\mu u - u \Vert^2_{L^2(Q_T)} &\leq c_\phi \int_{Q_T} \bigg(\int_{B(0,1)} |Eu(z+\mu\xi) - Eu(z)|^2 \d \xi \bigg) \d z \\ 
        &= c_\phi \int_{B(0,1)}\bigg( \int_{Q_T} |Eu(z+\mu\xi) - Eu(z)|^2 \d z \bigg) \d \xi \leq c_\phi \int_{B(0,1)} \Vert \tau_{\mu \xi}Eu - Eu \Vert_{L^2(\mathbb{R}^{d+1})}^2 \d \xi, 
    \end{aligned}
\end{equation*}
where, with slight abuse of notation, $\tau_h$ is the usual translation operator in the final line (and is no longer the coordinate with that same name). Observe that, for any $\mu>0$ and $\xi \in \mathbb{R}^{d+1}$, 
\begin{equation}\label{eq:equality of translated norms}
    \Vert \tau_{\mu \xi} Eu \Vert^2_{L^2(\mathbb{R}^{d+1})} = \int_{\mathbb{R}^{d+1}} |Eu(z+\mu\xi)|^2 \d z = \int_{\mathbb{R}^{d+1}} |Eu(z)|^2 \d z = \Vert Eu \Vert^2_{L^2(\mathbb{R}^{d+1})}, 
\end{equation}
where we used the translation invariance of Lebesgue measure. Note that the latter is well-defined by virtue of Lemma \ref{lem:sob extension time}. Hence, we have the bound 
\begin{equation*}
    \Vert \tau_{\mu \xi}Eu - Eu \Vert_{L^2(\mathbb{R}^{d+1})}^2 \leq 4 \Vert Eu \Vert^2_{L^2(\mathbb{R}^{d+1})}
\end{equation*}
which is uniform in $\mu$ and $\xi$ and is integrable on $B(0,1)$, along with, for a.e.~fixed $\xi \in B(0,1)$, 
\begin{equation*}
    \lim_{\mu \to 0} \Vert \tau_{\mu \xi}Eu - Eu \Vert_{L^2(\mathbb{R}^{d+1})}^2 = 0, 
\end{equation*}
since the translation operator $\tau_h$ is a bounded linear operator from $L^2(\mathbb{R}^{d+1})$ to itself, and is therefore continuous. Thus, an application of the Dominated Convergence Theorem yields 
\begin{equation}\label{eq:L2 cylinder conv}
    \begin{aligned}
        \Vert R_\mu u - u \Vert^2_{L^2(Q_T)} &\leq c_\phi \int_{B(0,1)} \Vert \tau_{\mu \xi}Eu - Eu \Vert_{L^2(\mathbb{R}^{d+1})}^2 \d \xi \to 0 \qquad \text{as } \mu \to 0, 
    \end{aligned}
\end{equation}
as required. Similarly, by writing $z=(t,x)$ and $\xi = (\tau,y)$, and letting $w(t,x) := \nabla Eu(t,x)$ (where we emphasize that the derivative is only taken with respect to the space coordinates, i.e., $\nabla Eu = \nabla_x Eu$), 
\begin{equation*}
    \Vert \nabla(R_\mu u - u) \Vert^2_{L^2(Q_T)} = \int_{Q_T} \bigg| \int_{\mathbb{R}^{d+1}} \phi_\mu(z-\xi) \big( w(\xi) - w(z) \big) \d \xi \bigg|^2 \d z. 
\end{equation*}
Following the previous argument to the letter, we find $\Vert \nabla(R_\mu u - u) \Vert^2_{L^2(Q_T)} \leq c_\phi \int_{B(0,1)} \Vert \tau_{\mu \xi}w - w \Vert_{L^2(\mathbb{R}^{d+1})}^2 \d \xi$. Noting that $w \in L^2(\mathbb{R}^{d+1})$ by virtue of Lemma \ref{lem:sob extension time}, we follow the same argument involving the Dominated Convergence Theorem to find that the latter integral vanishes, whence 
\begin{equation*}
    \Vert R_\mu u - u \Vert^2_{L^2(Q_T)} + \Vert \nabla(R_\mu u - u) \Vert^2_{L^2(Q_T)} \to 0 \qquad \text{as } \mu \to 0, 
\end{equation*}
and \eqref{eq:strong L2H1 conv smoother} follows. The proof is complete. 
\end{proof}

\begin{rem}\label{rem:time deriv L2H1 depends mu reg}
Note that, due to the formula \eqref{eq:time deriv Amu E}, it follows that 
\begin{equation*}
    \begin{aligned}
        | \partial_t R_\mu u (t,x) | &= \bigg| \int_\mathbb{R} \int_{\mathbb{R}^d} \partial_t \phi_\mu (t-\tau,x-y) Eu(\tau,y) \d y \d \tau  \bigg|  \leq \Vert \partial_t \phi_\mu \Vert_{L^2(\mathbb{R}^{d+1})} \Vert Eu \Vert_{L^2(\mathbb{R}^{d+1})}, 
    \end{aligned}
\end{equation*}
by the H\"{o}lder inequality, and from which it follows from Lemma \ref{lem:sob extension time} that $\Vert \partial_t R_\mu u \Vert_{L^\infty(Q_T)} \leq C_\mu \Vert u \Vert_{L^2(Q_T)}$, where the finiteness of the constant $C_\mu$ follows from the smoothness and compact support of $\phi_\mu$. The same strategy with an additional space derivative yields $\Vert \partial_t \nabla R_\mu u \Vert_{L^\infty(Q_T)} \leq \Vert \partial_t \nabla \phi_\mu \Vert_{L^2(\mathbb{R}^{d+1})} \Vert E u \Vert_{L^2(Q_T)}$. In view of this, relabelling $C_\mu$ as the relevant positive constant depending on $\mu,\Omega,T,$ (which, incidentally, will blow up in the limit as $\mu \to 0$), we have shown that there holds 
\begin{equation}\label{eq:time deriv L2H1 depends mu reg}
    \Vert \partial_t R_\mu u \Vert_{L^\infty(0,T;W^{1,\infty}(\Omega))} \leq C_\mu\Vert u \Vert^2_{L^2(Q_T)} \qquad \forall u \in L^2(0,T;H^1(\Omega)). 
\end{equation}
\end{rem}

The proof of the previous lemma also implies the following result. 

\begin{cor}\label{cor:smoothing operator L2 cylinder conv}
Fix $\mu>0$. The smoothing operator $R_{\mu}$ of Lemma \ref{lem:smoothing operator}, defined explicitly by 
\begin{equation}
    R_{\mu} := A_\mu \circ E: L^2(0,T;H^1(\Omega)) \to C^\infty(\bar{Q}_T), 
\end{equation}
admits, for some positive constant $C$ independent of $\mu$, the estimate 
\begin{equation}\label{eq:L2 cylinder unif bound Rnu smoother}
    \Vert R_\mu u \Vert_{L^2(Q_T)} \leq C \Vert u \Vert_{L^2(Q_T)}    \qquad \forall u \in L^2(0,T;H^1(\Omega)). 
\end{equation}
As such, it is a bounded linear operator from $L^2(Q_T)$ to $C^\infty(\bar{Q}_T)$ equipped with the subspace norm-topology of $L^2(Q_T)$. Moreover, given any $u \in L^2(0,T;H^1(\Omega))$, we have the strong convergence 
\begin{equation}\label{eq:strong L2 cylinder conv smoother}
    \Vert R_\mu u - u \Vert_{L^2(Q_T)} \to 0 \qquad \text{as } \mu \to 0. 
\end{equation}
\end{cor}

\begin{proof}
Recall from the proof of Lemma \ref{lem:smoothing operator}; in particular the estimates leading up to \eqref{eq:L2 cylinder conv}, that 
\begin{equation*}
    \Vert R_\mu u - u \Vert_{L^2(Q_T)} \to 0 \qquad \text{as } \mu \to 0, 
\end{equation*}
i.e., \eqref{eq:strong L2 cylinder conv smoother} holds. Similarly, recall from the estimates leading up to and including \eqref{eq:for L2 cyl unif bound} that 
\begin{equation*}
    \Vert R_\mu u \Vert_{L^2(Q_T)} \leq \Vert E u \Vert_{L^2(\mathbb{R}^{d+1})}. 
\end{equation*}
The estimate \eqref{eq:L2 cylinder unif bound Rnu smoother} now follows immediately from the above and Lemma \ref{lem:sob extension time}. 
\end{proof}

\end{appendices}

\end{document}